\documentclass[11pt,reqno]{amsart}
\usepackage{amssymb,latexsym}
\usepackage{lipsum}
\usepackage{mathtools}
\usepackage{comment}
\usepackage{amsthm}
\usepackage{enumerate}

\usepackage[letterpaper, top = 2cm, bottom = 2cm,right = 1cm, left = 1cm]{geometry}
\date{}

\newtheorem{theorem}{Theorem}[section]
\newtheorem{lemma}[theorem]{Lemma}

\theoremstyle{definition}
\newtheorem{definition}[theorem]{Definition}
\newtheorem{example}[theorem]{Example}

\newtheorem{corollary}[theorem]{Corollary}
\newtheorem{prop}[theorem]{Proposition}

\theoremstyle{remark}
\newtheorem{rem}[theorem]{Remark}

\numberwithin{equation}{section}

\title{Potentials for some Tensor Algebras}
\author{Raymundo Bautista and Daniel L\'{o}pez-Aguayo}
\address{Centro de Ciencias Matem\'{a}ticas UNAM Campus Morelia, Apartado Postal 61-3 Xangari, Morelia, Michoac\'{a}n 58089, M\'{e}xico.}
\email{raymundo@matmor.unam.mx}

\address{Escuela de Ingenier\'{i}a y Ciencias, Tecnologico de Monterrey, Boulevard Felipe \'{A}ngeles 2003\\ Pachuca de Soto, Hidalgo 42080, M\'{e}xico.}
\email{dlopez.aguayo@itesm.mx}

\begin{document}
\maketitle
\begin{abstract}
This paper generalizes former works of Derksen, Weyman and Zelevinsky about quivers with potentials. We consider semisimple finite-dimensional algebras $E$ over a field $F$, such that $E \otimes_{F} E^{op}$ is semisimple. We assume that $E$ contains a certain type of $F$-basis which is a generalization of a multiplicative basis. We study potentials belonging to the algebra of formal power series, with coefficients in the tensor algebra over $E$, of any finite-dimensional $E$-$E$-bimodule on which $F$ acts centrally. In this case, we introduce a cyclic derivative and to each potential we associate a Jacobian ideal. Finally, we develop a mutation theory of potentials, which in the case that the bimodule is $Z$-free, it behaves as the quiver case; but allows us to obtain realizations of a certain class of skew-symmetrizable integer matrices.
\end{abstract}

\tableofcontents

\begin{section}{Introduction}
There have been distinct generalizations of the notion of a quiver with potential and mutation where the underlying $F$-algebra, $F$ being a field, is replaced by more general algebras; see, for example, \cite{3,6,8}. In this paper, instead of working with a quiver, we consider a tensor algebra over $M$, where $M$ is an $E$-$E$-bimodule. We consider the case in which $E=\displaystyle \prod_{i=1}^{n} E_{i}$, where each $E_{i}$ is a semisimple finite-dimensional $F$-algebra, and each $E_{i}$ has an $F$-basis satisfying the conditions given at the beginning of Section \ref{sec7}. Then we develop a mutation theory for potentials lying in the topological algebra $\mathcal{F}_{E}(M)$, the $\langle M \rangle$-adic completion of the tensor algebra $T_{E}(M)$, where $M$ is any $E$-$E$-bimodule where $F$ acts centrally and of finite-dimension over $F$. In order to do this, in the first five sections, we treat the case when $M$ is a direct sum of $E$-$E$-bimodules of the form $E_{i} \otimes_{F} E_{j}$. We call such bimodules $Z$-free, but we do not impose any special conditions on the $F$-basis of the $E_{i}$. In the second part of the paper, we consider semisimple algebras $E_{i}$ having $F$-bases with the above-mentioned conditions, but now $M$ is any $E$-$E$-bimodule where $F$ acts centrally, and of finite-dimension over $F$. We embed $M$ into a $Z$-free $E$-$E$-bimodule $\hat{M}$, and then using our results for the $Z$-free bimodules we develop a mutation theory for potentials.
We now describe the content of each section. In Section \ref{sec2}, we consider semisimple finite-dimensional $F$-algebras $E_{1},\ldots,E_{n}$, $E=\displaystyle \prod_{i=1}^{n} E_{i}$, and $M$ an $E$-$E$-bimodule of finite-dimension over $F$ and $F$ acts centrally on such bimodule. Then, we introduce $\mathcal{F}_{E}(M)$, this is the $\langle M \rangle$-adic completion of the tensor algebra $T_{E}(M)$, where $\langle M \rangle$ is the two-sided ideal generated by $M$. We view $\mathcal{F}_{E}(M)$ as formal power series over $M$; then we provide a description (analogous to that of \cite{4}) of the topological algebra isomorphisms $\mathcal{F}_{E}(M) \rightarrow \mathcal{F}_{E}(M)$. In Section \ref{sec3}, we assume that $E=\displaystyle \prod_{i=1}^{n} E_{i}$, where each $E_{i}$ is a semisimple algebra. We introduce the notion of a $Z$-free bimodule and throughout this section, until Section \ref{sec5}, we assume that $M$ is $Z$-free. Then, following \cite{12}, we introduce a cyclic derivation $\mathfrak{h}: \mathcal{F}_{E}(M) \rightarrow \operatorname{End}_{E}(\mathcal{F}_{E}(M))$. Using this map, we define a cyclic derivative $\delta: \mathcal{F}_{E}(M) \rightarrow \mathcal{F}_{E}(M)$ by taking $\delta(f)=\mathfrak{h}(1)$ for each $f \in \mathcal{F}_{E}(M)$. Using this cyclic derivative we introduce, for any $u$ in the dual space of $M_{E}$, a partial cyclic derivative $\delta_{u}: \mathcal{F}_{E}(M) \rightarrow \mathcal{F}_{E}(M)$. In Section \ref{sec4}, for every potential $P$ of $\mathcal{F}_{E}(M)$, we define a two-sided closed ideal $R(P)$ of $\mathcal{F}_{E}(M)$. Our definition is given in terms of a $Z$-free generating set of $M$ and $F$-bases of each $E_{i}$. We prove that $R(P)$ is invariant under algebra isomorphisms $\mathcal{F}_{E}(M) \rightarrow \mathcal{F}_{E}(M')$ fixing pointwise the algebra $E$. This implies that $R(P)$ is independent of the choice of a $Z$-free generating set of $M$ and of the choice of $F$-bases of the $E_{i}$. In Section \ref{sec5}, given two potentials $P$ and $P'$ in $\mathcal{F}_{E}(M)$ we establish a necessary condition for the existence of an algebra automorphism $\varphi$ of $\mathcal{F}_{E}(M)$ such that $\varphi(P)$ is cyclically-equivalent to $P'$, i.e there exists an element $z$ in the closure of the $F$-vector subspace of $\mathcal{F}_{E}(M)$ generated by all elements of the form $fg-gf$, such that $P'=\varphi(P)+z$. In Section \ref{sec6}, we deal with $E$-$E$-bimodules which are not necessarily $Z$-free. Here we introduce the notion of a polarization of an $E$-$E$-bimodule.  The main result of this section is Theorem \ref{theo5.3} (Splitting Theorem). This theorem is analogous to the well-known Splitting Theorem (\cite[Theorem 4.6]{4}) but in a more general setting.
In Section \ref{sec7}, we consider semisimple algebras $E=\displaystyle \prod_{i=1}^{n} E_{i}$, such that the enveloping algebra $E \otimes_{F} E^{op}$ is semisimple and each $E_{i}$ has an $F$-basis satisfying certain properties. Here we define the notion of the Jacobian ideal $J(P)$ of $P$, where $P$ is any potential in $\mathcal{F}_{E}(M)$; it is shown that this ideal coincides with the ones considered in \cite{3,6,8}. In the rest of the paper we consider $E$-$E$-bimodules which are finite direct sums of $E_{i}$-$E_{j}$-bimodules and which are $E_{i}$-free on the left and $E_{j}$-free on the right. In Section \ref{sec8}, inspired by \cite{4}, we introduce the concept of premutation of an $E$-$E$-bimodule. Also, we introduce the notion of premutation of a potential $P$. In Section \ref{sec9}, we compute the square of the premutation of a potential $P$ and relate it to $P$ via a certain algebra isomorphism. Also, to each $E$-$E$-bimodule $M$, we define a skew-symmetrizable matrix $B(M)$, which plays the role of the matrix $B(A)$ defined in \cite[p.29]{4}. In Section \ref{sec10}, we examine a procedure proposed in \cite{6,8} for constructing an $F$-algebra $E$ and an $E$-$E$-bimodule $M$, such that the matrix $B(M)$ is skew-symmetrizable. Based on this procedure, and using the construction given in \cite{3}, we give a construction using the same arithmetic given in \cite{6}, but allowing $F$ to be an algebraically closed field. In Section \ref{sec11}, we present a definition of mutation of potentials of $\mathcal{F}_{E}(M)$. The main result of this section is that mutation is an involution on the set of right-equivalence classes of all reduced potentials. Finally, in Section \ref{sec12}, we consider $Z$-free bimodules. We prove that if $F$ is an infinite field, then for each finite sequence $(k_{l},\ldots,k_{1})$ of elements of $[1,n]$, there exist potentials which are $(k_{l},\ldots,k_{1})$-non-degenerate in the sense of \cite{4}. 
\end{section}

\begin{section}{Formal Power Series} \label{sec2}

Let $F$ be any field and let $E_{1}, \ldots ,E_{n}$ be finite dimensional semisimple $F$-algebras. Let $E=\displaystyle \prod _{i=1}^{n}E_{i}$ and let $M$ be an $E$-$E$-bimodule.  For $i=1,\ldots,n$, denote by $e_{i}$ the image of the the identity of $E_{i}$ in $E$ under the canonical inclusion $E_{i} \hookrightarrow E$. Then, in $E$, we have
$$1=\sum _{i=1}^{n}e_{i}.$$

Throughout this paper we will always assume that the underlying field $F$ acts centrally in all $E$-$E$-bimodules and that all such bimodules are finite-dimensional over $F$.  

Let $M$ be  an $E$-$E$-bimodule and let $T_{E}(M)=E\oplus M\oplus M^{\otimes 2}\oplus \ldots$ be the tensor algebra of $M$ over $E$. In what follows, we are going to define a topological algebra containing $T_{E}(M)$ as a dense subalgebra. 

\begin{definition} \label{def1.1} Define the algebra of formal power series over $M$ as the set
$$\mathcal{F}_{E}(M)=\prod _{i=0}^{\infty}M^{\otimes i}$$
where $M^{\otimes 0}=E$. 
\end{definition}
Then $\mathcal{F}_{E}(M)$ consists of all the functions
$$f:\mathbb{Z}^{\geq 0}\rightarrow T_{E}(M)$$
such that $f(i)\in M^{\otimes i}$ for all $i\in \mathbb{Z}^{\geq 0}$. We endow $\mathcal{F}_{E}(M)$ with an $F$-algebra structure by means of the following operations.
For $x,y\in F$ and $f,g\in \mathcal{F}_{E}(M)$ we define
$$(xf+yg)(i)=xf(i)+yg(i)$$
and 
$$(fg)(i)=\sum _{s+t=i}f(s)g(t).$$

One can check that $\mathcal{F}_{E}(M)$, with the above operations, is an $F$-algebra. We now define a metric on $\mathcal{F}_{E}(M)$ as follows: for $f\neq 0$ in $\mathcal{F}_{E}(M)$ we put $\nu (f)=\mathrm{min}\{i | f(i)\neq 0\}$; then for $f, g$ in $\mathcal{F}_{E}(M)$, we define $d(f,g)=0$ if $f=g$; and for $f\neq g$ we let $d(f,g)=2^{-\nu (f-g)}$. One can verify that this is a metric, so it defines a topology on $\mathcal{F}_{E}(M)$ such that the sum and multiplication are continuous maps, so $\mathcal{F}_{E}(M)$ is in fact a topological algebra. To each element $\alpha =\displaystyle \sum _{i=0}^{l}\alpha _{i}\in T_{E}(M)$, with $\alpha _{i}\in M^{\otimes i}$, we associate the element $\iota (\alpha )\in \mathcal{F}_{E}(M)$ given by $\iota (\alpha )(i)=\alpha _{i}$ for $i=0,\ldots,l$ and $\iota (\alpha )(j)=0$ for $j>l$.
The function $\iota :T_{E}(M) \rightarrow \mathcal{F}_{E}(M)$ is a monomorphism of $F$-algebras. \\

Throughout this paper we identify $\alpha$ with $\iota(\alpha)$; in this way, $\mathcal{F}_{E}(M)$ contains $T_{E}(M)$ as a dense subalgebra.

We recall that if $\{f_{r}\}_{r=0}^{\infty }$ is a sequence of elements of $\mathcal{F}_{E}(M)$, then the series $\displaystyle \sum _{r=0}^{\infty }f_{r}$ exists if and only if the sequence of partial sums $\{J_{m}\}_{m=0}^{\infty }$, with $J_{m}=\displaystyle \sum _{r=0}^{m}f_{r}$ , converges and in such case
$$\sum _{r=0}^{\infty }f_{r}=\displaystyle \lim_{m\rightarrow \infty }J_{m}.$$

\begin{definition} \label{def1.2} A sequence $\{f_{s}\}_{s=0}^{\infty }$ of elements of $\mathcal{F}_{E}(M)$ is said to be summable if for each non-negative integer $m$, the set $T(m)=\{s|f_{s}(m)\neq 0 \}$ is finite.
\end{definition}

\begin{rem} \label{rem1.3} If the sequence  $\{f_{s}\}_{s=0}^{\infty }$ of elements of $\mathcal{F}_{E}(M)$  is summable
then $\displaystyle \sum _{s=0}^{\infty }f_{s}$ exists and  
$$\left(\sum _{s=0}^{\infty }f_{s}\right)(i)=\sum _{s\in T(i)}f_{s}(i)$$
for all $i\in \mathbb{Z}^{\geq 0}$.
\end{rem}

\begin{rem} \label{rem1.4} If the sequences $\{ f_{s}\}_{s=0}^{\infty }$ and $\{ g_{s}\}_{s=0}^{\infty }$ of $\mathcal{F}_{E}(M)$ are summable, then the sequence $\{h_{s}\}_{s=0}^{\infty }$, where $h_{s}=\displaystyle \sum _{u+v=s}f_{u}g_{v}$, is summable and
$$\left(\sum _{s=0}^{\infty }f_{s}\right)\left(\sum _{s=0}^{\infty }g_{s}\right)=\sum _{s=0}^{\infty }h_{s}.$$
\end{rem}

\begin{rem} \label{rem1.5} If $f\in \mathcal{F}_{E}(M)$, then the sequence $\{f(i)\}_{i=0}^{\infty}$ of elements in $\mathcal{F}_{E}(M)$ is summable and
$$f=\sum _{i=0}^{\infty }f(i).$$
\end{rem}

{\bf Notation.} For each non-negative integer $i$, we define $\mathcal{F}_{E}(M)^{\geq i}$ as the set of elements $f\in \mathcal{F}_{E}(M)$ such that $f(j)=0$ for $j<i$. Observe that $\mathcal{F}_{E}(M)^{\geq i}$ is a closed ideal of $\mathcal{F}_{E}(M)$.

\begin{lemma} \label{lem1.6} Let $M^{\prime }$ be an $E$-$E$-bimodule and $\phi: M\rightarrow \mathcal{F}_{E}(M^{\prime })^{\geq 1}$ a morphism of $E$-$E$-bimodules, then there is a morphism of topological algebras
$$\varphi :\mathcal{F}_{E}(M)\rightarrow \mathcal{F}_{E}(M^{\prime })$$
such that $\varphi |_{E}=id_{E}$ and $\varphi |_{M}=\phi $.
\end{lemma} 
 
\begin{proof} By the universal property of the tensor algebra $T_{E}(M)$ there is a morphism of $F$-algebras $\hat{\phi }:T_{E}(M)\rightarrow \mathcal{F}_{E}(M^{\prime })$ such that $\hat{\phi }|_{E}=id_{E}$ and $\hat{\phi }|_{M}=\phi $. We extend the function $\hat{\phi}$ to a function $\varphi: \mathcal{F}_{E}(M)\rightarrow \mathcal{F}_{E}(M^{\prime })$ in the following way: let $f\in \mathcal{F}_{E}(M)$, then since $\phi (M)\subseteq \mathcal{F}_{E}(M^{\prime })^{\geq 1}$, we have
  $$\hat{\phi }(f(i))\in \mathcal{F}_{E}(M^{\prime })^{\geq i}.$$ 
Therefore the sequence $\{ \hat{\phi }(f(i))\}_{i=0}^{\infty }$ is summable, and we define $\varphi (f)=\displaystyle \sum _{i=0}^{\infty }\hat{\phi }(f(i))$.
 
 Clearly, for $f$ and  $g$ in $\mathcal{F}_{E}(M)$ we have $\varphi (f+g)=\varphi (f)+\varphi (g)$. Moreover,
 by Remark \ref{rem1.4}
 $$\varphi (f)\varphi (g)=\left(\sum _{i=0}^{\infty }\hat{\phi }(f(i))\right)\left(\sum _{i=0}^{\infty }\hat{\phi }(g(i))\right)=\sum _{l=0}^{\infty }\left(\sum _{s+t=l}\hat{\phi }(f(s))\hat{\phi }(g(t))\right)$$ 
 $$=\sum _{l=0}^{\infty }\hat{\phi }\left(\sum _{s+t=l}f(s)g(t)\right)=\sum _{l=0}^{\infty }\hat{\phi}(fg(l))=\varphi (fg).$$
 Therefore $\varphi $ is an algebra morphism. Since $\varphi (\mathcal{F}_{E}(M)^{\geq i})\subseteq \mathcal{F}_{E}(M^{\prime })^{\geq i}$ for all $i>0$, then $\varphi $ is a continuous map. It follows that $\varphi $ is a morphism of topological algebras, as claimed.
 \end{proof}

It follows from Lemma \ref{lem1.6} that any algebra morphism $\varphi :\mathcal{F}_{E}(M)\rightarrow \mathcal{F}_{E}(M^{\prime })$ with $\varphi |_{E}=id_{E}$  and $\varphi (M)\subseteq \mathcal{F}_{E}(M^{\prime })^{\geq 1}$ is completely determined by the pair of morphisms of $E$-$E$-bimodules $\varphi ^{(1)}:M\rightarrow M^{\prime }, \varphi ^{(2)}:M\rightarrow \mathcal{F}_{E}(M^{\prime })^{\geq 2}$ such that for $m\in M$, $\varphi (m)=\varphi ^{(1)}(m)+\varphi ^{(2)}(m)$.
 
\begin{prop} \label{prop1.7} The morphism $\varphi :\mathcal{F}_{E}(M)\rightarrow \mathcal{F}_{E}(M^{\prime })$ determined by the pair $(\varphi ^{(1)}, \varphi ^{(2)})$, with $\varphi |_{E}=id_{E}$, is an algebra isomorphism if and only if $\varphi ^{(1)}$ is an isomorphism of $E$-$E$-bimodules.
\end{prop}
 
\begin{proof} Suppose first that $\varphi $ is an algebra isomorphism, then there is an algebra morphism 
 $$\varrho:\mathcal{F}_{E}(M^{\prime })
 \rightarrow \mathcal{F}_{E}(M)$$
such that $\varrho$ is the inverse of $\varphi $. Clearly $\varrho |_{E}=id_{E} $, then $\varrho $ is determined by a pair
  $(\varrho ^{(1)}, \varrho ^{(2)})$.  For $m\in M$, we have:
  $$m=\varrho ^{(1)}\varphi ^{(1)}(m)+\varrho ^{(2)}(\varphi^{(1)}(m))+\varrho(\varphi^{2}(m))$$
Since the last three summands in the above equality are in $\mathcal{F}_{E}(M)^{\geq 2}$, then $\varrho ^{(1)}\varphi ^{(1)}(m)=m$ and therefore $\varrho ^{(1)}\varphi ^{(1)}=id_{M}$. In the same way one proves that $\varphi ^{(1)}\varrho ^{(1)}=id_{M^{\prime }}$. This shows that $\varphi ^{(1)}$ is an isomorphism of $E$-$E$-bimodules.
  
Conversely, assume that $\varphi ^{(1)}$ is an isomorphism of $E$-$E$-bimodules. Consider first the case $M=M^{\prime }$ and $\varphi ^{(1)}=id_{M}$. Take the morphism of $E$-$E$-bimodules $\psi =id _{\mathcal{F}_{E}(M)}-\varphi $. Then
  for $\mu =m_{1} \cdots m_{l}$, with $m_{i}\in M$, we have
  $$\psi (\mu )=m_{1} \cdots m_{l}-\varphi (m_{1}) \cdots \varphi (m_{l})=m_{1} \cdots m_{l}-(m_{1}-\varphi ^{(2)}(m_{1})) \cdots (m_{l}-\varphi ^{(2)}(m_{l}))$$
  $$=m_{1} \cdots m_{l}-m_{1} \cdots m_{l}+\mu ^{\prime }=\mu ^{\prime }$$
  where $\mu ^{\prime }$ is the product of $l$ elements of the form $m_{i}$ or $\varphi ^{(2)}(m_{i})$ with at least
  one factor of the form $\varphi ^{(2)}(m_{i})\in \mathcal{F}_{E}(M)^{\geq 2}$. From here we conclude that
  $\psi (\mu )\in \mathcal{F}_{E}(M)^{\geq l+1}$, and consequently 
  $$\psi (\mathcal{F}_{E}(M)^{\geq l})\subseteq \mathcal{F}_{E}(M)^{\geq l+1}$$
for all non-negative integers $l$. Then for $f\in \mathcal{F}_{E}(M)$, the sequence $\{ \psi ^{k}(f) \}_{k=0}^{\infty}$ is summable.
Define $\varrho (f)=\displaystyle \sum _{k=0}^{\infty} \psi ^{k}(f)$.  The map $\varrho :\mathcal{F}_{E}(M)\rightarrow \mathcal{F}_{E}(M)$ is continuous.    For $f\in \mathcal{F}_{E}(M)$, we have
$$\varphi \varrho (f)=\varrho (f)-\psi \varrho (f)=\sum _{k=0}^{\infty }\psi ^{k}(f)-\sum _{k=1}^{\infty }\psi ^{k}(f)=f,$$
and 
$$\varrho \varphi (f)=\varrho (f)-\varrho \psi (f)=\sum _{k=0}^{\infty }\psi ^{k}(f)-\sum _{k=1}^{\infty }\psi^{k}(f)=f.$$

Therefore $\varrho \varphi =\varphi \varrho =id_{\mathcal{F}_{E}(M)}$, showing that $\varphi $ is an isomorphism.

Now consider the general case $\varphi :\mathcal{F}_{E}(M)\rightarrow \mathcal{F}_{E}(M^{\prime })$, given by the pair of morphisms $\varphi ^{(1)}:M\rightarrow M'$ and $\varphi ^{(2)}:M\rightarrow \mathcal{F}_{E}(M')^{\geq 2}$ with $\varphi ^{(1)}$ an isomorphism. Consider $\widehat{\varphi } :\mathcal{F}_{E}(M)\rightarrow \mathcal{F}_{E}(M^{\prime })$, the morphism determined by the pair $(\varphi ^{(1)},0)$. The morphism $\widehat{\varphi }$ has an inverse map $\varrho $ determined by the pair $(\varrho ^{(1)},0)$, with $\varrho ^{(1)}$ being the inverse of $\varphi ^{(1)}$. Note that the morphism $\varrho \varphi :\mathcal{F}_{E}(M)\rightarrow \mathcal{F}_{E}(M)$ is determined by the pair $(id_{M},\varrho \varphi ^{(2)})$, thus $\varrho \varphi $ is an isomorphism. Since $\varrho $ is an isomorphism we conclude that $\varphi $ is also an isomorphism. This completes the proof.
\end{proof}

\begin{lemma} \label{lem1.8} If $\mathcal{I}$ is an ideal of $\mathcal{F}_{E}(M)$ such that $\mathcal{I}\cap E=0$, then $\mathcal{I}\subseteq \mathcal{F}_{E}(M)^{\geq 1}$.
\end{lemma} 

\begin{proof} Suppose that $\mathcal{I}$ is not contained in $\mathcal{F}_{E}(M)^{\geq 1}$, then there exists a non-zero element $w\in \mathcal{I}$ such that $w=w_{0}+w_{1}$ with $w_{0}$ a non-zero element in $E$ and $w_{1}\in \mathcal{F}_{E}(M)^{\geq 1}$. Since $E$ is a semisimple algebra, then the identity of $E$ is a sum of orthogonal central primitive idempotents. Therefore there exists an idempotent $\epsilon $ with $\epsilon w_{0}\neq 0$; so we may assume that $w_{0}$ is a non-zero element of $\epsilon E$. Since $\epsilon E$ is a simple algebra, then there exist elements $x_{1},\dots,x_{r},y_{1},\dots,y_{r}$ of $E$ such that $\epsilon =\displaystyle \sum _{i=1}^{r}x_{i}w_{0}y_{i}$. Then $\displaystyle \sum _{i=1}^{r}x_{i}wy_{i}=\epsilon - v_{1}$ with $v_{1}\in \epsilon \mathcal{F}_{E}(M)^{\geq 1}\epsilon $.

Thus $\epsilon -v_{1}\in \mathcal{I}$ and $\epsilon =(\epsilon -v)\left(\epsilon +\displaystyle \sum _{l=1}^{\infty }v^{l}\right)=\epsilon $. It follows that $\epsilon \in \mathcal{I}\cap E$, a contradiction. This proves our claim.
\end{proof}

\begin{prop} \label{prop1.9} If $\varphi :\mathcal{F}_{E}(M)\rightarrow \mathcal{F}_{E}(M^{\prime })$ is an algebra isomorphism such that $\varphi |_{E}=id_{E}$, then $\varphi (M)\subseteq \mathcal{F}_{E}(M^{\prime })^{\geq 1}$. Therefore in this case $\varphi $ is completely determined by a
pair of morphisms of $E$-$E$-bimodules $\varphi ^{(1)}:M\rightarrow M^{\prime }$, $\varphi ^{(2)}:M\rightarrow \mathcal{F}_{E}(M^{\prime })^{\geq 2}$.
\end{prop}

\begin{proof} Take $\mathcal{I}=\mathcal{F}_{E}(M)^{\geq 1}$, since $\varphi $ is an isomorphism, then
$\varphi (\mathcal{I})$ is an ideal of $\mathcal{F}_{E}(M^{\prime })$ such that $\varphi (\mathcal{I})\cap E=0$. Then
$$\varphi (M)\subseteq \varphi (\mathcal{I})\subseteq \mathcal{F}_{E}(M^{\prime })^{\geq 1}.$$ This proves our Proposition.
\end{proof}
  
\begin{prop} \label{prop1.10} Let $\{\rho _{k}\}_{k=0}^{\infty }$ be a sequence of algebra automorphisms of $\mathcal{F}_{E}(M)$, such that each $\rho _{k}$ is determined by the pair $(id_{M},\rho _{k}^{(2)})$ with $\rho _{0}^{(2)}=0$ and $\rho _{k}^{(2)}(M)\subseteq \mathcal{F}_{E}(M)^{\geq k+1}$ for all non-negative integers $k$. 
Then
\begin{enumerate}
\item[(i)] The sequence $\{\rho _{k} \cdots \rho _{1}\rho  _{0}(f)\}_{k=0}^{\infty }$ converges for all $f\in \mathcal{F}_{E}(M)$.
\item[(ii)] The morphism $\rho: \mathcal{F}_{E}(M) \rightarrow \mathcal{F}_{E}(M)$, which sends an element $f\in \mathcal{F}_{E}(M)$ into $\displaystyle \lim_{k\to \infty }\rho _{k} \cdots \rho _{0}(f)$, is an algebra automorphism such that $\rho |_{E}=id_{E}$. 
\end{enumerate}
\end{prop}

\begin{proof} First we  define $\rho :\mathcal{F}_{E}(M)\rightarrow \mathcal{F}_{E}(M)$ as follows. For $f\in \mathcal{F}_{E}(M)$, and any non-negative integer $l$
$$\rho (f)(l)=\rho _{l}\rho  _{l-1}\cdots \rho  _{0}(f)(l).$$

Observe that $\rho (f)(l)=\rho  _{l+1}\rho  _{l}\cdots \rho  _{0}(f)(l)$. Consequently
$$\rho (f)(l)=\rho _{l+k}\rho _{l+k-1}\cdots \rho _{0}(f)(l)$$
for any non-negative integer $k$.

Clearly $\rho $ is an $F$-linear map and for $f,g\in \mathcal{F}_{E}(M)$ one has
$$\rho (fg)(l)=\rho  _{l}\rho  _{l-1} \cdots \rho _{0}(fg)(l)=\rho _{l} \cdots \rho _{0}(f)\rho _{l} \cdots \rho _{0}(g)(l)$$
$$=\sum _{k_{1}+k_{2}=l}\rho _{l} \cdots \rho _{0}(f)(k_{1})\rho _{l} \cdots \rho _{0}(g)(k_{2})=\sum _{k_{1}+k_{2}=l}
\rho _{k_{1}} \cdots \rho _{0}(f)(k_{1})\rho _{k_{2}} \cdots \rho _{0}(g)(k_{2})$$
$$=\sum _{k_{1}+k_{2}=l}\rho (f)(k_{1})\rho (g)(k_{2})=\rho (f)\rho (g)(l),$$
and therefore $\rho (fg)=\rho (f)\rho (g)$. Here  $\rho (1)=1$ and $\rho (m)-m \in \mathcal{F}_{E}(M)^{\geq 2}$,
for each $m\in M$. Therefore $\rho $ is an algebra automorphism. Moreover, for each $l^{\prime }\leq l$, the following holds
$$\left(\rho (f)-\varphi _{l} \cdots \varphi _{0}(f)\right)(l^{\prime })=\rho (f)(l^{\prime })-\varphi _{l} \cdots \varphi _{0}(f)(l^{\prime})
=$$ 
$$\rho _{l^{\prime }} \cdots \rho _{0}(f)(l^{\prime })-\rho _{l^{\prime }} \cdots \rho _{0}(f)(l^{\prime })=0.$$
Therefore $\rho (f)-\rho _{l} \cdots \rho _{0}(f)\in \mathcal{F}_{E}(M)^{\geq l+1}$. This proves that
$$\rho (f)=\displaystyle \lim_{k\to \infty }(\rho _{k} \cdots \rho _{0}(f)).$$
The proof of the proposition is now complete.
\end{proof}

\end{section}

\begin{section}{Cyclic derivatives} \label{sec3}

Let $A$ be an $F$-algebra, we recall that a cyclic derivation in the sense of Rota-Sagan-Stein \cite{12} is an $F$-linear function $\mathfrak{h}:A\rightarrow \mathrm{End}_{F}(A)$ such that

\begin{equation} \label{eq2.1}
\mathfrak{h}(f_{1}f_{2})(f)=\mathfrak{h}(f_{1})(f_{2}f)+\mathfrak{h}(f_{2})(ff_{1})
\end{equation}
for all $f,f_{1},f_{2} \in A$. From a cyclic derivation $\mathfrak{h}$, we obtain a cyclic derivative $\delta :A\rightarrow A$ given by
$$\delta (f)=\mathfrak{h}(f)(1)$$
for all $f\in A$. From $(2.1)$ one obtains

\begin{equation} \label{eq2.2}
\delta (f_{1}f_{2})=\mathfrak{h}(f_{1})(f_{2})+\mathfrak{h}(f_{2})(f_{1})
\end{equation}
for all $f_{1}, f_{2}\in A$. In particular, $\delta (f_{1}f_{2})=\delta (f_{2}f_{1})$.

In this section we will construct a cyclic derivative for the $F$-algebra $\mathcal{F}_{E}(M)$, with $E$ is as in Section \ref{sec2}
and $M\cong E\otimes _{Z}M_{0}\otimes _{Z}E$, where $Z=\displaystyle \sum _{i=1}^{n}Fe_{i}$ is a subalgebra of $E$ and
$M_{0}$ is a $Z$-$Z$-subbimodule of $M$ of finite dimension over $F$.

\begin{definition} \label{def2.1} An $E$-$E$-bimodule $M$ is called $Z$-free if there exists a $Z$-$Z$-subbimodule $M_{0}$ of $M$ such that the multiplication map $E\otimes _{Z}M_{0}\otimes _{Z}E\rightarrow M$ which sends $x\otimes m\otimes y$ to
$xmy$ with $x,y\in E$ and $m\in M_{0}$ is an isomorphism of $E$-$E$-bimodules. The $Z$-$Z$-subbimodule $M_{0}$ is called a $Z$-free generator for $M$.
\end{definition}
In order to construct a cyclic derivative on $\mathcal{F}_{E}(M)$, we first define a cyclic derivation on the
tensor algebra $A=T_{E}(M)$. In order to do this, first consider the map
$$\hat{u}:A\times A\rightarrow A$$
given by $\hat{u}(f,g)=\displaystyle \sum _{i=1}^{n}e_{i}gfe_{i}$, for $f$ and $g\in A$. This is an $F$-bilinear map which is $Z$-balanced. Indeed if $s=ze_{j}$ with $z\in F$, one has
\begin{center}
$\hat{u}(fze_{j},g)=\displaystyle \sum _{i=1}^{n}e_{i}gfze_{j}e_{i}=ze_{j}gfe_{j}=\hat{u}(f,ze_{j}g).$
\end{center}
Thus there exists $u: A \otimes_{Z} A \rightarrow A$ such that $u(a \otimes b)=\hat{u}(a,b)$.
Now we define an $F$-derivation $\Delta: A\rightarrow A\otimes _{Z}A$ as follows: for $s\in E$, we put $\Delta (s)=1\otimes s-s\otimes 1$; for $m\in M_{0}$, $\Delta (m)=1\otimes m$. Then we define $M\rightarrow  T_{E}(M)$ such that for $s_{1},s_{2}\in E$ and $m\in M_{0}$ we have
$$\Delta (s_{1}ms_{2})=\Delta (s_{1})ms_{2}+s_{1}\Delta (m)s_{2}+s_{1}m\Delta (s_{2})$$
the above map is well defined because $M\cong E\otimes _{Z}M_{0}\otimes _{Z}E$ via  the multiplication map.
Now $\Delta $ can be extended to an  $F$-derivation on $A$. We define $\mathfrak{h}:A\rightarrow \mathrm{End}_{F}(A)$ as follows 
$$\mathfrak{h}(f)(g)=u(\Delta (f)g).$$

We have
\begin{align*}
&\mathfrak{h}(f_{1}f_{2})(f)=u(\Delta (f_{1}f_{2})f) \\
&=u(\Delta (f_{1})f_{2}f)+u(f_{1}\Delta (f_{2})f) \\
&=u(\Delta (f_{1})f_{2}f)+u(\Delta (f_{2})ff_{1}) \\
&=\mathfrak{h}(f_{1})(f_{2}f)+\mathfrak{h}(f_{2})(ff_{1}).
\end{align*}

Hence we obtain a cyclic derivation $\mathfrak{h}$ on $A$. We are now going to extend $\mathfrak{h}$ to $\mathcal{F}_{E}(M)$.
Take $f,g \in \mathcal{F}_{E}(M)$, then by the definition of $\mathfrak{h}$ we see that $\mathfrak{h}(f(i))(g(j))\in M^{\otimes i+j}$, and then we define $\mathfrak{h}(f)(g)(l)=\displaystyle \sum _{i+j=l}\mathfrak{h}(f(i))(g(j))$ for all non-negative integers $l$.

\begin{prop} \label{prop2.2} The $F$-linear map $\mathfrak{h}:\mathcal{F}_{E}(M)\rightarrow \mathrm{End}_{F}(\mathcal{F}_{E}(M))$ is a
cyclic derivation.
\end{prop}

\begin{proof} Let $f, f_{1},f_{2} \in \mathcal{F}_{E}(M)$. Then for any non-negative integer $l$ we have
$$\mathfrak{h}(f_{1}f_{2})(f)(l)=\sum _{i+j=l}\mathfrak{h}((f_{1}f_{2})(i))(f(j))=\sum _{i_{1}+i_{2}+j=l}\mathfrak{h}(f_{1}(i_{1})f_{2}(i_{2}))(f(j))$$
$$=\sum _{i_{1}+i_{2}+j=l}\mathfrak{h}(f_{1}(i_{1}))(f_{2}(i_{2})f(j))+\sum _{i_{1}+i_{2}+j=l}\mathfrak{h}(f_{2}(i_{2}))(f(j)f_{1}(i_{1}))$$
$$=\sum _{i_{1}+t=l}\mathfrak{h}(f_{1}(i_{1}))((f_{2}f)(t))+\sum _{i_{2}+r=l}\mathfrak{h}(f_{2}(i_{2}))((ff_{1})(r))$$
$$=\mathfrak{h}(f_{1})(f_{2}f)(l)+\mathfrak{h}(f_{2})(ff_{1})(l)=\left(\mathfrak{h}(f_{1})(f_{2}f)+\mathfrak{h}(f_{2})(ff_{1})\right)(l).$$
The result follows.
\end{proof}

From the above we obtain a cyclic derivative $\delta $ with
$$\delta (f)=\mathfrak{h}(f)(1)$$
for $f\in \mathcal{F}_{E}(M)$.

\begin{rem} \label{rem2.3} For a fixed $f\in \mathcal{F}_{E}(M)$, $\mathfrak{h}(f)(g)$ is continuous in $g$; likewise, for a fixed $g$, $\mathfrak{h}(f)(g)$ is continuous
in $f$. In particular, if $\{u_{i}\}_{i=0}^{\infty }$ is a sequence of elements of $\mathcal{F}_{E}(M)$ such that
$\displaystyle \sum _{i=0}^{\infty }u_{i}$ exists, then $\displaystyle \sum _{i=0}^{\infty }\mathfrak{h}(f)(u_{i})$ and $\displaystyle \sum _{i=0}^{\infty }\mathfrak{h}(u_{i})(g)$ exist.
Moreover
$$\mathfrak{h}(f)\left(\sum _{i=0}^{\infty }u_{i}\right)=\sum _{i=0}^{\infty }\mathfrak{h}(f)(u_{i})$$
$$\mathfrak{h}\left(\sum _{i=0}^{\infty }u_{i}\right)(g)=\sum _{i=0}^{\infty }\mathfrak{h}(u_{i})(g).$$
\end{rem}
\begin{definition} \label{def2.4} For an $E$-$E$-bimodule $N$ we define the cyclic part of $N$ as
$$N_{cyc}=\sum _{i=1}^{n}e_{i}Ne_{i}.$$
Similarly, if $u\in \mathcal{F}_{E}(M)$, we define $u_{cyc}=\displaystyle \sum _{i=1}^{n}e_{i}ue_{i}$.
An element of $\mathcal{F}_{E}(M)_{cyc}$ is called a \emph{potential} or a \emph{cyclic} element.
\end{definition}

\begin{lemma} \label{lem2.5} If $s\in E$ and $f,g,\in \mathcal{F}_{E}(M)$ then
\begin{enumerate}[(a)]
\item $\mathfrak{h}(s)(f)=[sf-fs]_{cyc}$
\item $\mathfrak{h}(sf)(g)=\mathfrak{h}(f)(gs)+[s,fg]_{cyc}$
\item $\mathfrak{h}(fs)(g)=\mathfrak{h}(f)(sg)+[s,gf]_{cyc}$
\end{enumerate}
\end{lemma}

where $[a,b]:=ab-ba$ for $a,b \in \mathcal{F}_{E}(M)$.

\begin{proof} By Equation (\ref{eq2.1}),  $(a)$ implies $(b)$ and $(c)$, so it suffices to prove $(a)$. For any non-negative integer $l$, we have
$$\mathfrak{h}(s)(f)(l)=\mathfrak{h}(s)(f(l))=u(\Delta (s)f(l))=u((1\otimes s)f(l)-(s\otimes 1)f(l))$$
$$=\sum _{i=1}^{n}e_{i}(sf(l)-f(l)s)e_{i}=[sf-fs]_{cyc}(l).$$
From here we obtain $(a)$. This proves our claim.
\end{proof}

\begin{rem} \label{rem2.6} If $s\in Z$, and $f\in \mathcal{F}_{E}(M)$, then
$$\mathfrak{h}(fs)(g)=\mathfrak{h}(f)(sg)$$
and
$$\mathfrak{h}(sf)(g)=\mathfrak{h}(f)(gs).$$
\end{rem}

\begin{definition} \label{def2.7} An element $f\in \mathcal{F}_{E}(M)$ is called directed if there exist idempotents $e_{i}, e_{j}$ such that $f=e_{j}fe_{i}$; in this case $i=\sigma (f)$ is called the \emph{start} of $f$ and $j=\tau (f)$ the \emph{end} of $f$.
\end{definition}

\begin{prop} \label{prop2.8} The following properties hold
\begin{enumerate}[(a)]
\item If  $f, f_{1},\dotsc,f_{l} \in \mathcal{F}_{E}(M)$ then
$$\mathfrak{h}(f_{1}\dotsm f_{l})(f)=\mathfrak{h}(f_{1})(f_{2} \dotsm f_{l}f)+\mathfrak{h}(f_{2})(f_{3}\dotsm f_{l}ff_{1}) + \dotsc+\mathfrak{h}(f_{l})(ff_{1}\dotsm f_{l-1}).$$
\item If $m \in EM_{0}$ then $\Delta (m)=1\otimes m$, and for $f\in \mathcal{F}_{E}(M)$
$$\mathfrak{h}(m)(f)=(m f)_{cyc}$$
\item If $m_{1},\dots,m_{l}\in EM_{0}$ and $f\in \mathcal{F}_{E}(M)$, then
$$\mathfrak{h}(m_{1} \dotsm m_{l})(f)=(m_{1}m_{2} \dotsm m_{l}f+m_{2}m_{3}\dotsm m_{l}fm_{1}+ \dotsc+m_{l}fm_{1}\dotsm m_{l-1})_{cyc}.$$
In particular
$$\delta (m_{1}m_{2} \dotsm m_{l})=(m_{1}m_{2}\dotsm m_{l}+m_{2}m_{3}\dotsm m_{l}m_{1}+\dotsc+m_{l}m_{1}\dotsm m_{l-1})_{cyc}.$$
\item If $m_{1},\dotsc,m_{l}$ are directed elements of $EM_{0}$ and $m_{1} \dotsm m_{l}$ is a non-zero cyclic element, then
$$\delta (m_{1} \dotsm m_{l})=m_{1}\dotsm m_{l}+m_{2}m_{3} \dotsm m_{l}m_{1}+\dotsc+m_{l}m_{1}\dotsm m_{l-1}.$$
\end{enumerate}
\end{prop}

\begin{proof}
\begin{enumerate}[(a)]
\item  Follows by induction on $l$ from Equation (\ref{eq2.1}).
\item It suffices to consider the case $m=sm_{0}$ with $s\in E$ and $m_{0}\in M_{0}$. We have
 $$\Delta (sm_{0})=\Delta (s)m_{0}+s\Delta (m_{0})=1\otimes sm_{0}-s\otimes m_{0}+s(1\otimes m _{0})
=1\otimes sm_{0}.$$
Then for $m\in EM_{0}$, one has 
$$\mathfrak{h}(m)(f)=\sum _{j=0}^{\infty }\mathfrak{h}(m)(f(j))=\sum _{j=0}^{\infty }u(\Delta (m)f(j))$$
$$=\sum _{j=0}^{\infty }u(1\otimes mf(j))=\sum _{j=0}^{\infty }(mf(j))_{cyc}=(mf)_{cyc}.$$
\item Follows from (a) and (b).
\item Follows from (c).
\end{enumerate}
\end{proof}
The proof is now complete.
\begin{definition} \label{def2.9} A dual basis $\{z_{i},\psi _{i}\}_{i=1}^{k}$ of $M_{E}$ is called a $Z$-dual basis if for any $m\in M$, $m\in EM_{0}$ if and only if $\psi _{i}(m)\in Z$ for all $i=1,\dotsc, k$.
\end{definition}

For our $Z$-free $E$-$E$-bimodule $M$, we can find a $Z$-dual basis of $M_{E}$. Indeed, we have $EM_{0}\cong E\otimes _{Z}M_{0}$, then $EM_{0}$ is a projective right $Z$-module. Let $\{z_{1}, \dotsc ,z_{k};\lambda _{1}, \dotsc,\lambda _{k}\}$ be a dual basis for $(EM_{0})_{Z}$ with $z_{i}\in EM_{0}$, $\lambda _{i}\in \mathrm{Hom}_{Z}((EM_{0})_{Z},Z)$. Since $M_{E}\cong EM_{0}\otimes _{Z}E$, each map $\lambda _{i}$ can be extended to a $E$-morphism $\psi _{i}:M_{E}\rightarrow E$ such that
$\psi _{i}(ms)=\lambda _{i}(m)s$ for $m\in EM_{0}$ and $s\in E$. Now one can see that $\{z_{1}, \dotsc,z_{k};\psi _{1},\dotsc,\psi _{k}\}$ is a $Z$-dual basis of $M_{E}$.

Let $\psi \in \mathrm{Hom}_{E}(M_{E},E)$. We will extend this map to an $F$-linear endomorphism of $\mathcal{F}_{E}(M)$  denoted also by $\psi $. For $s\in E$, we define $\psi (s)=0$; for $M^{\otimes l}$, with $l>1$, we define $\psi (m_{1}\otimes \dotsm \otimes m_{l})=\psi (m_{1})m_{2} \dotsm m_{l}\in M^{\otimes (l-1)}$ for $m_{1},\dots,m_{l}\in M$; and for $f\in \mathcal{F}_{E}(M)$ we define $\psi (f)\in \mathcal{F}_{E}(M)$ such that $\psi (f)(l-1)=\psi (f(l))$ for each non-negative integer $l$. Then
$$\psi (f)=\sum _{l=0}^{\infty }\psi (f(l)).$$

\begin{definition} \label{def2.10} For $\psi \in M^{*}=\mathrm{Hom}_{E}(M_{E},E)$ we define  $\delta _{\psi }:\mathcal{F}_{E}(M)\rightarrow \mathcal{F}_{E}(M)$ as 
$$\delta _{\psi }(f)=\psi (\delta (f))=\sum _{l=0}^{\infty }\psi \delta (f(l))$$
for all $f\in \mathcal{F}_{E}(M)$.
\end{definition}

Observe that if $\displaystyle \sum _{i=0}^{\infty }f_{i}$ exists then 
$$\delta _{\psi }\left(\sum _{i=0}^{\infty }f_{i}\right)=\psi \left(\mathfrak{h}\left(\sum _{i=0}^{\infty }f_{i}\right)(1)\right)=\psi \left(\sum _{i=0}^{\infty }\mathfrak{h}(f_{i})(1)\right)=\psi \left(\sum _{i=0}^{\infty }\delta (f_{i})\right)=\sum _{i=0}^{\infty }\delta _{\psi }(f_{i}).$$

\begin{prop} \label{prop2.11} Suppose that $M^{\prime }$ is an $E$-$E$-bimodule $Z$-freely generated by $M_{0}^{\prime }$. Denote by $\mathfrak{h}^{\prime }:\mathcal{F}_{E}(M^{\prime })\rightarrow \mathrm{End}_{F}(\mathcal{F}_{E}(M^{\prime }))$ the cyclic derivation constructed using $L$ and $M_{0}^{\prime }$, and denote by $\delta ^{\prime }$ the corresponding cyclic derivative. Let $\varphi: \mathcal{F}_{E}(M)\rightarrow \mathcal{F}_{E}(M^{\prime })$ be an isomorphism
of topological algebras with $\varphi |_{E}=id_{E}$. Let $\{z_{1},\dotsc,z_{k};\psi _{1},\dotsc,\psi _{k}\}$ be a $Z$-dual basis of
$M_{E}$. Then if $P$ is a potential in $\mathcal{F}_{E}(M)$, we have
$$\delta ^{\prime }(\varphi (P))=\sum _{j=1}^{k}\mathfrak{h}^{\prime }(\varphi (z_{j}))(\varphi (\delta _{\psi _{j}}(P))).$$
\end{prop}

\begin{proof} First assume $P=m_{1}\dotsm m_{l}$ with $m_{1},\dotsc,m_{l}\in EM_{0}$. Then by (a) of Proposition \ref{prop2.8}
$$\delta ^{\prime }(\varphi (P))=$$
$$\delta^{\prime} (\varphi (m_{1})\dotsm \varphi (m_{l}))=\mathfrak{h}^{\prime }(\varphi (m_{1}))(\varphi (m_{2})\dotsm \varphi (m_{l}))+\dotsc+\mathfrak{h}^{\prime }(\varphi (m_{l}))(\varphi (m_{1})\dotsm \varphi (m_{l-1}))$$
$$=\sum _{j=1}^{k}\mathfrak{h}^{\prime }(\varphi (z_{j})\psi _{j}(m_{1}))(\varphi (m_{2}\dotsm m_{l}))+\dotsc+\sum _{j=1}^{k}\mathfrak{h}^{\prime }(\varphi (z_{j})\psi _{j}(m_{l}))(\varphi (m_{1}\dotsc m_{l-1}))$$
Since our dual basis is a $Z$-dual basis, we have that $\psi _{j}(m_{i})\in Z$ for all $j=1,\dotsc,l$ and $i=1,\dotsc,l$. Then by Remark \ref{rem2.6} we obtain
$$\delta ^{\prime }(\varphi (P))=$$
$$\sum _{j=1}^{k}\mathfrak{h}^{\prime }(\varphi (z_{j}))(\varphi (\psi _{j}(m_{1})m_{2}\dotsm m_{l}+\dotsc+\psi _{j}(m_{l})m_{1}\dotsm m_{l-1}))=$$
$$\sum _{j=1}^{k}\mathfrak{h}^{\prime }(\varphi (z_{j}))(\varphi ((\delta _{\psi_{j}}(m_{1}\dotsm m_{l})))=\sum _{j=1}^{k}\mathfrak{h}^{\prime }(\varphi (z_{j}))(\varphi (\delta _{\psi_{j}}(P))).$$

Now any potential  in $M^{\otimes l}$ is a sum of potentials  of the form $Pt$ with $P$ as above and $t\in E$. Then
$$\delta ^{\prime }(\varphi (Pt))=\delta^{\prime}(\varphi (P)t) =\delta^{\prime}(t\varphi (P))=\delta^{\prime} (\varphi (tP)).$$
Therefore
$$\delta ^{\prime }(\varphi (Pt))=\sum _{j=1}^{k}\mathfrak{h}^{\prime }(\varphi (z_{j}))(\varphi (\delta _{\psi_{j}}(tP)))=\sum _{j=1}^{k}\mathfrak{h}^{\prime }(\varphi (z_{j}))(\varphi (\delta _{\psi_{j}}(Pt))).$$

From here it follows that our Proposition holds for any potential in $M^{\otimes l}$. Finally, if $P$ is any potential in
$\mathcal{F}_{E}(M)$ we have
$$\delta ^{\prime }(\varphi (P))=\sum _{l=0}^{\infty }\delta ^{\prime }(\varphi (P(l)))$$
 $$=\sum _{l=0}^{\infty }\sum _{j=1}^{k}\mathfrak{h}^{\prime }(\varphi (z_{j}))(\varphi (\delta _{\psi_{j}}(P(l))))$$
 $$=\sum _{j=1}^{k}\mathfrak{h}^{\prime }(\varphi (z_{j}))\left(\sum _{l=0}^{\infty }\varphi (\delta _{\psi_{j}}(P(l)))\right)=\sum _{j=1}^{k}\mathfrak{h}^{\prime }(\varphi (z_{j}))(\varphi (\delta _{\psi_{j}}(P))).$$
This completes the proof.
 \end{proof}
 
{\bf Notation.} We denote by $[\mathcal{F}_{E}(M),\mathcal{F}_{E}(M)]$ the closure in $\mathcal{F}_{E}(M)$ of the
$F$-subspace generated by all the elements of the form $[f,g]=fg-gf$ with $f,g\in \mathcal{F}_{E}(M).$

\begin{definition} \label{def2.12} Two potentials $P$ and $P^{\prime }$ are called cyclically equivalent if $P-P^{\prime }\in  [\mathcal{F}_{E}(M),\mathcal{F}_{E}(M)]$.
\end{definition}

Clearly, if $P$ and $P^{\prime }$ are  cyclically equivalent potentials then $\delta (P)=\delta (P^{\prime })$.

\begin{definition} \label{def2.13} A subset $T$ of $EM_{0}$ is called a $Z$-free local basis for $M_{E}$ if for each idempotent
$e_{i}$ of $S$, and for every $i=1,\dotsc,n$, $T_{i}=T\cap Me_{i}$  is a set of $E_{i}$-free generators of $Me_{i}$ and
$T=\displaystyle \bigcup _{i=1}^{n}T_{i}$.
\end{definition}
 
A $Z$-free local basis for $M_{E}$ is obtained by choosing for $i=1,\dotsc,n$ an $F$-basis $T_{i}$ of $EM_{0}e_{i}$ and taking
$T=\displaystyle \bigcup _{i=1}^{n}T_{i}$.

Throughout this paper, if $T$ is a $Z$-free local basis for $M_{E}$ and $c \in T \cap e_{i}Me_{j}$, we define a right $E$-morphism $c^{\ast}$ as the map $c^{\ast}: M_{E} \rightarrow E$ such that $c^{\ast}(c^{\prime})=0$ if $c^{\prime} \in T \setminus \{c\}$; and $c^{\ast}(c)=e_{j}$. In case $T$ is a $Z$-free local basis for $_{E}M$ and $c \in T \cap e_{i}Me_{j}$, $^{\ast}c$ is defined similarly. \\

{\bf Notation.} If $A$ and $B$ are subsets of $\mathcal{F}_{E}(M)$ we denote by $AB$ the closure of the
$F$-subspace generated by all the elements of the form $fg$ with $f\in A$ and $g\in B$. In the case $A=B$, we put
$A^{2}=AA$.

\begin{prop} \label{prop2.14} Let $\varphi $ be an algebra automorphism of $\mathcal{F}_{E}(M)$ determined by the pair
$(id_{M}, \varphi ^{(2)})$. Let $T$ be a $Z$-free local basis for $M_{E}$ and let $I$ be the closure of the two-sided ideal of $\mathcal{F}_{E}(M)$ generated by all the elements $\varphi ^{(2)}(c)$ with $c\in T$. Then if $P$ is a potential in
$\mathcal{F}_{E}(M)^{\geq 3}$, we have
$$\varphi (P)-P=\sum _{c\in T}\varphi ^{(2)}(c)\delta _{c^{*}}(P)+\alpha +\gamma $$
with $\alpha \in \mathcal{F}_{E}(M)^{\geq 1}I^{2}$ and $\gamma \in [\mathcal{F}_{E}(M),\mathcal{F}_{E}(M)]$.
\end{prop}

\begin{proof} First assume that $P=c_{1}c_{2}\dotsm c_{l}$ with $l\geq 3$ and $c_{1},\dotsc,c_{l} \in T$. Then
\begin{align*}
\varphi (P)&=\varphi (c_{1})\varphi (c_{2})\dotsm \varphi (c_{l}) \\
&=(c_{1}+\varphi ^{(2)}(c_{1}))(c_{2}+\varphi ^{(2)}(c_{2}))\dotsm(c_{l}+\varphi ^{(2)}(c_{l})) \\
&=c_{1}c_{2}\dotsm c_{l}+\varphi ^{(2)}(c_{1})c_{2}\dotsm c_{l}+\varphi ^{(2)}(c_{2})c_{3}\dotsm c_{l}c_{1}+\dotsc+\varphi ^{(2)} (c_{l})c_{1}\dotsm c_{l-1}+\mu +z
\end{align*}
where $\mu $ is a sum of products $x_{1}\dotsm x_{l}$ and each of the elements $x_{1},\dotsc,x_{l}$ lie in the set  
$\{c_{1},\dotsc,c_{l},\varphi ^{(2)}(c_{1}),\dotsc,\varphi ^{(2)}(c_{l})\}$, and there exist $i$, $j$, $i\neq j$, with $x_{i}, x_{j}\in \{\varphi ^{(2)}(c_{1}),\dotsc,\varphi ^{(2)}(c_{l})\}$,
and $z\in [\mathcal{F}_{E}(M),\mathcal{F}_{E}(M)]$. We claim that each of these products is of the form
$$\sum _{c\in T}cy_{c}+\sum _{c\in T}\varphi ^{(2)}(c)w_{c}+z^{\prime }$$
with $y_{c},w_{c}\in I^{2}\cap \mathcal{F}_{E}(M)^{\geq l-1}$ and $z^{\prime }\in [\mathcal{F}_{E}(M),\mathcal{F}_{E}(M)]$. If all the elements $x_{i}$ are in $\{\varphi ^{(2)}(c_{1}),\dotsc,\varphi ^{(2)}(c_{l})\}$, then
$x_{1}\dotsm x_{d}=\varphi ^{(2)}(c_{i_{1}})y_{c_{i_{1}}}$, with $y_{c_{i_{1}}}\in I^{2}\cap \mathcal{F}_{E}(M)^{\geq l-1}$; so we may assume that $x_{1}=c_{i_{1}}$, then $x_{1}x_{2} \dotsm x_{l}=c_{i_{1}}y_{i_{1}}$ with $y_{i_{1}}\in \mathcal{F}_{E}(M)^{\geq l-1}\cap I^{2}$. This proves our claim.  Now take any potential $P\in \mathcal{F}_{E}(M)^{\geq 3}$. Then we may assume that each $P(l)$ is a sum of elements of the form $m_{1}m_{2} \dotsm m_{l}$ with $m_{1},\dotsc,m_{l}\in EM_{0}$. Therefore $P(l)$ is an $F$-linear combination of elements of the form $c_{1}\dotsm c_{l}$, so we have
$$\varphi (P(l))-P(l)=\sum _{c\in T}\varphi ^{(2)}(c)\delta _{c^{*}}(P(l))+\sum _{c\in T}cy^{l}_{c}+\sum _{c\in T}\varphi ^{(2)}(c)w^{l}_{c}+z_{l}$$
with $y_{c},w_{c}\in \mathcal{F}_{E}(M)^{\geq l-1}\cap I^{2}$ and $z_{l}\in \mathcal{F}_{E}(M)^{\geq l}\cap  [\mathcal{F}_{E}(M),\mathcal{F}_{E}(M)]$. Therefore the sequences $\{y^{l}_{c}\}_{l=3}^{\infty }$,
$\{w_{c}^{l}\}_{l=3}^{\infty }$ and $\{z_{l}\}_{l=3}^{\infty } $ are summable. Then
\begin{align*}
&\varphi (P)-P=\sum _{l=3}^{\infty }\sum _{c\in T}\varphi ^{(2)}(c)(\delta _{c_{*}}(P(l)))+\sum _{l=3}^{\infty}\sum _{c\in T}(cy^{l}_{c}+\varphi ^{(2)}(c)w_{c}^{l})+ \sum _{l=3}^{\infty }z_{l} \\
&=\sum _{c\in T}\varphi ^{(2)}(c)\sum _{l=3}^{\infty }(\delta _{c_{*}}(P(l)))+\sum _{c\in T}\left[\varphi ^{(2)}(c)\sum _{l=3}^{\infty }y_{c}^{l}+\varphi ^{(2)}(c)\sum _{l=3}^{\infty }w_{c}^{l}\right]+\sum _{l=3}^{\infty }z_{l} \\
&=\sum _{c\in T}\varphi ^{(2)}(c)\delta _{c^{*}}(P)+\alpha +\gamma
\end{align*}
with $\alpha \in \mathcal{F}_{E}(M)^{\geq 1}I^{2}$ and $\gamma \in [\mathcal{F}_{E}(M),\mathcal{F}_{E}(M)]$.
\end{proof}
\end{section}
\begin{section}{The ideal $R(P)$} \label{sec4}

Throughout this section $E$, $Z$, $M$ and $M_{0}$ are as in the precedent sections. Now we are going to construct a
$Z$-free local basis for $M$ consisting of directed elements as follows.
\begin{enumerate}
 \item For each $e_{i}M_{0}e_{j}\neq 0$, choose an $F$-basis $A_{i,j}$ of $e_{i}M_{0}e_{j}$. Then $A=\displaystyle \bigcup _{i,j}A_{i,j}$ is an $F$-basis of $M_{0}$.
\item For each $e_{i}E$, choose an $F$-basis $L_{i}$ of $e_{i}E$ and let $L=\displaystyle \bigcup _{i=1}^{n}L_{i}$. If $t\in L_{i}$, 
$t^{*}:e_{i}E\rightarrow e_{i}F$ is the $F$-linear map such that $t^{*}(t)=e_{i}$ and $t^{*}(s)=0$ for $s\in L_{i} \setminus \{t\}$.
\end{enumerate}

Having $A$ and $L$, we consider the set $T=\{sa\}_{s\in L, a\in A}$ which is a $Z$-free local basis for $M$, consisting of directed elements.
\begin{lemma} \label{lem3.1} For $f\in \mathcal{F}_{E}(M)_{cyc}$ and $t\in E$ we have
$$\sum _{s\in L(\tau (a))}(sa)^{*}(tf-ft)s=0.$$
\end{lemma}

\begin{proof} If $f\in E$, the statement is clearly true; so we may assume that $f\in \mathcal{F}_{E}(M)^{\geq 1}$.
In this case, $f$ is a sum of elements of the form $rbg$ with $b\in A$, $r\in L(\tau (b))$ and $g\in \mathcal{F}_{E}(M)$.
Therefore it suffices to prove the lemma for elements of this form.
$$\sum _{s\in L(\tau (a))}(sa)^{*}(trbg)s=\sum _{s,s_{1}\in L(\tau (a))}(sa)^{*}(s_{1}s_{1}^{*}(tr)bg)s$$
$$=\sum _{s,s_{1}\in L(\tau (a))}(sa)^{*}(s_{1}b)gs_{1}^{*}(tr)s=\sum _{s\in L(\tau (a))}gs^{*}(tr)s\delta _{b,a}=gtr\delta _{b,a}.$$

Moreover,
$$\sum _{s\in L(\tau (a))}(sa)^{*}(rbgt)=\sum _{s\in L(\tau (a))}(sa)^{*}(rb)gt=gtr\delta _{b,a}.$$
The result follows. 
\end{proof}

\begin{prop} \label{prop3.2} Suppose $f$ and $g\in \mathcal{F}_{E}(M)$. If $a\in A$ and $t\in E$, then we have
\begin{enumerate}
 \item[(a)] $\displaystyle \sum _{s\in L(\tau (a))}(sa)^{*}\mathfrak{h}(tf)(g)s=\displaystyle \sum _{s\in L(\tau (a))}(sa)^{*}\mathfrak{h}(gt)s.$
 \item[(b)] $\displaystyle \sum _{s\in L(\tau (a))}(sa)^{*}\mathfrak{h}(ft)(g)s=\displaystyle \sum _{s\in L(\tau (a))}(sa)^{*}\mathfrak{h}(f)(tg)$
 \end{enumerate}
 \end{prop}
 
 $ $
 
\begin{proof} We first prove the following
 \begin{enumerate}[(i)]
\item  If $fg$ is cyclic, then 
$$\sum _{s\in L(\tau (a))}(sa)^{*}\mathfrak{h}(tf)(g)s=\sum _{s\in L(\tau (a))}(sa)^{*}\mathfrak{h}(f)(gt)s.$$
\item  If $gf$ is cyclic, then
$$\sum _{s\in L(\tau (a))}(sa)^{*}\mathfrak{h}(ft)(g)s=\sum _{s\in L(\tau (a))}(sa)^{*}\mathfrak{h}(f)(tg)$$
\end{enumerate}

Indeed, we have
$$\sum _{s\in L(\tau (a))}(sa)^{*}\mathfrak{h}(tf)(g)s=\sum _{s\in L(\tau (a))}(sa)^{*}(\mathfrak{h}(f)(gt)+[t,fg]_{cyc})s.$$

Here $fg$ is cyclic, then $[t,fg]_{cyc}=[t,fg]$, so by Lemma \ref{lem3.1} we obtain (i). The proof of (ii) is done in a similar way. Now if $f$ and $g$ are any elements of $\mathcal{F}_{E}(M)$, by Remark \ref{rem2.6} we have
$$\mathfrak{h}(f)(g)=\sum _{i=1}^{n}\mathfrak{h}(fe_{i})(g)=\sum _{i=1}^{n}\mathfrak{h}(fe_{i})(ge_{i})$$
and
$$\mathfrak{h}(f)(g)=\sum _{i=1}^{n}\mathfrak{h}(e_{i}f)(g)=\sum _{i=1}^{n}\mathfrak{h}(e_{i}f)(ge_{i}).$$

Then 
$$\sum _{s\in L(\tau (a))}(sa)^{*}\mathfrak{h}(tf)(g)s=\sum _{i=1}^{n}\sum _{s\in L(\tau (a))}(sa)^{*}(\mathfrak{h}(te_{i}f)(ge_{i}))s$$
but $e_{i}fge_{i}$ is cyclic and thus by (i) we have
$$\sum _{s\in L(\tau (a))}(sa)^{*}\mathfrak{h}(tf)(g)s=\sum _{i=1}^{n}\sum _{s\in L(\tau (a))}(sa)^{*}(\mathfrak{h}(e_{i}f)(ge_{i}t))s$$
$$=\sum _{s\in L(\tau (a))}(sa)^{*}\mathfrak{h}(f)(gt)s.$$
This proves $(a)$. A similar argument using (ii)  proves $(b)$.
\end{proof}

\begin{definition} \label{def3.3} Let $P$ be a potential in $\mathcal{F}_{E}(M)$, then $R(P)$ is the closure of the two-sided ideal in
$\mathcal{F}_{E}(M)$ generated by all the elements
$$X_{a^{*}}(P)=\sum _{s\in L(\tau (a))}\delta _{(sa)^{*}}(P)s$$
where $a$ runs over all elements of $A$.
\end{definition}

\begin{prop} \label{prop3.4} Let $P$ be a potential of the form $P=t_{1}a_{1}t_{2}a_{2}\cdots t_{l}a_{l}$ with
$a_{1},\dots,a_{l}\in A$ and for $i=1,\dots,l$, $t_{i}\in E_{\tau (a_{i})}$. Then, for $a\in A$
$$X_{a^{*}}(P)=\delta _{a,a_{1}}t_{2}a_{2}\cdots t_{l}a_{l}t_{1}+\delta _{a,a_{2}}t_{3}a_{3}\cdots t_{l}a_{l}t_{1}a_{1}t_{2}+ \hdots +\delta _{a,a_{l}}t_{1}a_{1}t_{2}a_{2}\cdots t_{l-1}a_{l-1}t_{l}.$$
\end{prop}

\begin{proof} We have 
$$X_{a^{*}}(P)=\sum _{s\in L(\tau (a))}(sa)^{*}\delta (t_{1}a_{1}\cdots t_{l}a_{l})s$$ 
$$=\sum _{s\in L(\tau (a))}(sa)^{*}(t_{1}a_{1}\cdots t_{l}a_{l}+t_{2}a_{2}\cdots t_{l}a_{l}t_{1}a_{1}+\hdots+t_{l}a_{l}t_{1}a_{1}\cdots t_{l-1}a_{l-1})s.$$

Using Lemma \ref{lem3.1} we obtain
$$X_{a^{*}}(P)=\sum _{s\in L(\tau (a))}(sa)^{*}(a_{1}\cdots t_{l}a_{l}t_{1}+a_{2}\cdots t_{l}a_{l}t_{1}a_{1}t_{2}+\hdots+a_{l}t_{1}a_{1}\cdots t_{l-1}a_{l-1}t_{l})s$$
$$=\delta _{a,a_{1}}t_{2}a_{2}\cdots t_{l}a_{l}t_{1}+\delta _{a,a_{2}}t_{3}a_{3}\cdots t_{l}a_{l}t_{1}a_{1}t_{2}+\hdots+\delta _{a,a_{l}}t_{1}a_{1}t_{2}a_{2}\cdots t_{l-1}a_{l-1}t_{l}.$$
\end{proof}

\begin{rem} \label{rem3.5} It follows from Proposition \ref{prop3.4} that the definition of $X_{a*}(P)$ is independent of the choice of an $F$-basis for $e_{1}E,\dots,e_{n}E$.
\end{rem}

Now consider a $Z$-free $E$-$E$-bimodule $M^{\prime}$ and $M_{0}^{\prime }$ a $Z$-$Z$-subbimodule of $M^{\prime}$ such that the multiplication map $E\otimes _{Z}M_{0}^{\prime} \otimes _{Z}E\rightarrow M^{\prime}$ is an isomorphism of $E$-$E$-bimodules. Choose, as in the case of $M$, an $F$-basis $A'(i,j)$ for $e_{i}M^{\prime }_{0}e_{j}\neq 0$. If $A^{\prime }=\displaystyle \bigcup _{i,j}A^{\prime }(i,j)$ then the set $T^{\prime }=\{sa^{\prime }\}_{s\in L(\tau (a^{\prime })),a^{\prime}}$ is a $Z$-free basis for
$M^{\prime }_{E}$. Then for any potential $Q\in \mathcal{F}_{E}(M^{\prime })$ we have defined the ideal $R(Q)$. \\

We now establish a similar result to Proposition \ref{prop2.11}. 

\begin{lemma} \label{lem3.6} Let $\varphi :\mathcal{F}_{E}(M)\rightarrow \mathcal{F}_{E}(M^{\prime })$ be an algebra isomorphism such that $\varphi |_{E}=id_{E}$, then for any potential $P\in \mathcal{F}_{E}(M)$ and $c\in A^{\prime }$, we have

$$ X_{c^{*}}(\varphi (P))=\sum _{t\in L(\tau (c)),a\in A}(tc)^{*}\mathfrak{h}^{\prime }(\varphi (a))(\varphi (X_{a^{*}}(P)))t.$$
\end{lemma}

\begin{proof} From the definition of $X_{c^{*}}$, we have
 $$X_{c^{*}}(\varphi (P))=\sum _{t\in L(\tau{c}) }(tc)^{*}\delta ^{\prime }_{c^{*}}(\varphi (P))t.$$
 
 Using Proposition \ref{prop2.11}, we obtain
 \begin{center}
 $X_{c^{*}}(\varphi (P))=\displaystyle \sum _{t\in L(\tau (c))}(tc)^{*}\left[\displaystyle \sum _{sa\in T}\mathfrak{h}^{\prime }(\varphi (sa))\varphi (\delta _{(sa)^{*}}(P))\right]t.$
 \end{center}
 By (a) of Proposition \ref{prop3.2} we have
 \begin{align*}
X_{c^{\ast}}(\varphi(P))&=\sum_{t\in L(\tau (c))}(tc)^{\ast}\left[\sum _{sa\in T}\mathfrak{h}^{\prime}(\varphi(a))\left(\varphi (\delta _{(sa)^{\ast}}(P)s)\right)\right]t\\
&=\sum_{t\in L(\tau (c))}(tc)^{\ast}\left[\sum _{a\in A}\mathfrak{h}^{\prime }(\varphi (a))\left(\varphi \left(\sum _{s\in L(\tau (a))}\delta _{(sa)^{\ast}}(P)s\right) \right)\right]t\\
&=\sum_{t \in L(\tau (c)),a \in A}(tc)^{\ast}\mathfrak{h}^{\prime}(\varphi (a))\left(\varphi (X_{a^{\ast}}(P))\right)t
\end{align*}
This completes the proof.
\end{proof}
\begin{theorem} \label{theo3.7} Let $\varphi :\mathcal{F}_{E}(M)\rightarrow \mathcal{F}_{E}(M^{\prime })$ be an algebra isomorphism such that $\varphi |_{E}=id_{E}$. Then 
$$\varphi (R(P))=R(\varphi (P))$$
for every potential $P$ of $\mathcal{F}_{E}(M) $.
\end{theorem}

\begin{proof} By Lemma \ref{lem3.6},  for $c\in A^{\prime }$ we have 
$$X_{c^{*}}(\varphi (P))=\sum _{t\in L(\tau (c)),a\in A}(tc)^{*}(\mathfrak{h}^{\prime }(\varphi (a))(\varphi (X_{a^{*}}(P)))t.$$

Each element $\varphi (a)(l)$, with $l \geq 2$, is the sum of elements of the form
$n_{1} \cdots n_{l}r$ with $n_{1},\dots, n_{l}$ directed elements in $ EM_{0}^{\prime }$ and $r\in e_{\sigma (n_{l})}E$.
By (b) of  Proposition \ref{prop3.2}, we have
 
 \begin{center}
 $\displaystyle \sum _{t\in L(\tau (c))}(tc)^{*}(\mathfrak{h}^{\prime }(n_{1}\cdots n_{l}r)(\varphi (X _{a^{*}}(P))))t=\displaystyle \sum _{t\in L(\tau (c))}(tc)^{*}(\mathfrak{h}^{\prime }(n_{1}\cdots n_{l})(\varphi (rX_{a^{*}}(P)))t.$
 \end{center}
 
 Taking $\alpha =\varphi (rX_{a^{*}}(P))$, and applying (c) of Proposition \ref{prop2.8}, yields
 $$(tc)^{*}\mathfrak{h}^{\prime }(n_{1} \cdots n_{l})(\varphi (\alpha ))t=$$
 $$((tc)^{*}(n_{1})n_{2}\cdots n_{l}\varphi (\alpha )+(tc)^{*}(n_{2})\cdots n_{l}\varphi (\alpha )n_{1}+\hdots+(tc)^{*}(n_{l})\varphi (\alpha )n_{1}\cdots n_{l-1}))t.$$
 Observe that if $l=1$, $(tc)^{*}\mathfrak{h}^{\prime }(n_{1})(\varphi (\alpha ))t=(tc)^{*}(n_{1})\varphi (\alpha )t$. Then $\displaystyle \sum _{t\in L(\tau (c))}(tc)^{*}(\mathfrak{h}^{\prime }(a(l))(X_{a^{*}}(P))t$ is a sum of elements of the form
 $$\gamma _{1}\varphi (\alpha )\gamma _{2}=\varphi (\varphi ^{-1}(\gamma _{1})\alpha \varphi ^{-1}(\gamma _{2}))$$
 where $\varphi ^{-1}(\gamma _{1})\alpha \varphi ^{-1}(\gamma _{2})\in R(P)\cap \mathcal{F}_{E}(M)^{\geq l+1}$.
 
 Therefore
 $$\sum _{t\in L(\tau (c)),a\in A}(tc)^{*}(\mathfrak{h}^{\prime }(\varphi (a)(l))(\varphi (X_{a^{*}}(P)))t=\varphi (z_{l}) $$
 with $z_{l}\in R(P)\cap \mathcal{F}_{E}(M)^{\geq l+1}$. Since the sequence $\{z_{l}\}_{l=1}^{\infty }$ is summable,
 then $\displaystyle \sum _{l=1}^{\infty }z_{l}\in R(P)$.
 
 We obtain
 $$X_{c^{*}}(\varphi (P))=\sum _{l=1}^{\infty }\left(\sum _{t\in L(\tau (c)),a\in A}(tc)^{*}(\mathfrak{h}^{\prime }(a(l))\varphi (X_{a^{*}}(P))t\right)$$
 $$=\sum _{l=1}^{\infty }\varphi (z_{l})=\varphi \left(\sum _{l=1}^{\infty }z_{l}\right)\in \varphi (R(P))$$
 Therefore $R(\varphi (P))\subseteq \varphi (R(P))$.  Moreover
 $$R(P)=R(\varphi ^{-1}(\varphi (P)))\subseteq \varphi ^{-1}(R(\varphi (P))),$$
 hence $\varphi (R(P))\subseteq R(\varphi (P))$. The result follows.  
 \end{proof}
 
\begin{prop} \label{prop3.8} Suppose $M$ is $Z$-freely generated by the $Z$-subbimodule $M_{0}^{\prime }$ of $M$,
 take $L^{\prime }(i)$ an $F$-basis for $e_{i}E$ and $A^{\prime }(i,j)$ an $F$-basis for $e_{i}M_{0}^{\prime }e_{j}\neq 0$. Let  $P$ be any potential of $\mathcal{F}_{E}(M)$, $R^{\prime }(P)$ the ideal using $L^{\prime }=\displaystyle \bigcup _{i=1}^{n}L^{\prime }(i)$ and $A^{\prime }=\displaystyle \bigcup _{i,j}A^{\prime }(i,j)$. Then
 $R^{\prime }(P)=R(P)$.
 \end{prop}
 
\begin{proof} By Remark \ref{rem3.5}, we may assume that $L=L^{\prime }$. Taking $\varphi =id_{\mathcal{F}_{E}(M)}$ in Theorem \ref{theo3.7} gives $R(P)=R^{\prime }(P)$.
\end{proof}
\end{section}

\begin{section}{Equivalence of Potentials} \label{sec5}
In this section, given potentials $P$ and $P'$ in $\mathcal{F}_{E}(M)$ we will see that if $P-P'$ lies in $R(P)^{2}$, then there exists an algebra automorphism $\varphi$ of $\mathcal{F}_{E}(M)$ such that $P'$ is cyclically-equivalent to $\varphi(P)$.

Throughout this section, if $f\in \mathcal{F}_{E}(M)$ we denote by $ \langle f \rangle$ the closure of the two-sided ideal generated by $f$ in $\mathcal{F}_{E}(M)$.

\begin{prop} \label{prop4.1} Suppose $f,g\in \mathcal{F}_{E}(M)^{\geq 2}$. Then $X_{a^{*}}(fg)$ is an element of  
$$ \mathcal{F}_{E}(M)^{\geq 1} \langle f \rangle+\langle f \rangle \mathcal{F}_{E}(M)^{\geq 1}+\mathcal{F}_{E}(M)^{\geq 1}\langle g\rangle+\langle g\rangle \mathcal{F}_{E}(M)^{\geq 1}.$$
\end{prop}

\begin{proof} We have
$$X_{a^{*}}(fg)=\sum _{s\in L(\tau (a))}(sa)^{*}\delta (fg)s=\sum _{s\in L(\tau (a))}(sa)^{*}(\mathfrak{h}(f)(g)+
\mathfrak{h}(g)(f))s.$$
Take any integer $l\geq 2$, then $f(l)$ is a sum of elements of the form $m_{1}\cdots m_{l}t$ with $m_{1},\hdots ,m_{l}\in EM_{0}$ and $t\in E$. Then by (b) of Proposition \ref{prop3.2} and (c) of Proposition \ref{prop2.8}
$$\sum _{s\in L(\tau (a))}(sa)^{*}\mathfrak{h}(m_{1}\cdots m_{l}t)(g)s=
\sum _{s\in L(\tau (a))}(sa)^{*}(\mathfrak{h}(m_{1}\cdots m_{l})(tg))s=$$
$$\sum _{s\in L(\tau (a))}((sa)^{*}(m_{1})m_{2}\cdots m_{l}tg+(sa)^{*}(m_{2})m_{3}\cdots m_{l}tgm_{1}+\hdots+(sa)^{*}(m_{l})tgm_{1}\cdots m_{l-1})s$$
and this element lies in $\mathcal{F}_{E}(M)^{\geq 1} \langle g \rangle+\langle g \rangle \mathcal{F}_{E}(M)^{\geq 1}$. Therefore
\begin{center}
$\displaystyle \sum _{s\in L(\tau (a))}(sa)^{*}\mathfrak{h}(f)(g)s=\displaystyle \sum _{l=2}^{\infty }\sum _{s\in L(\tau (a))}(sa)^{*}\mathfrak{h}(f(l))(g)s\in \mathcal{F}_{E}(M)^{\geq 1}\langle g \rangle+\langle g \rangle \mathcal{F}_{E}(M)^{\geq 1}.$
\end{center}

In a similar way,
$$\sum _{s\in L(\tau (a))}(sa)^{*}\mathfrak{h}(g)(f)s\in  \mathcal{F}_{E}(M)^{\geq 1}\langle f \rangle+\langle f \rangle \mathcal{F}_{E}(M)^{\geq 1}.$$
and the proof is now complete.
\end{proof}

\begin{prop} \label{prop4.2} If $P$ and $P^{\prime }$ are potentials in $\mathcal{F}_{E}(M)^{\geq 3}$ such that $P-P^{\prime } \in R(P)^{2}$, then
$R(P)=R(P^{\prime })$.
\end{prop}

$  $

\begin{proof} By Proposition \ref{prop4.1}, 
$$X_{a^{*}}(P)-X_{a^{*}}(P^{\prime })\subseteq \mathcal{F}_{E}(M)^{\geq 1}R(P)+R(P)\mathcal{F}_{E}(M)^{\geq 1}$$

Therefore
\begin{equation} \label{eq4.1}
R(P)\subseteq R(P^{\prime })+ \mathcal{F}_{E}(M)^{\geq 1}R(P)+R(P)\mathcal{F}_{E}(M)^{\geq 1}
\end{equation}

In particular, $R(P)\subseteq R(P^{\prime })+\mathcal{F}_{E}(M)^{\geq 3}$. Suppose we have already proved that
$$R(P)\subseteq R(P^{\prime })+\mathcal{F}_{E}(M)^{\geq l}$$ 
Then using Equation (\ref{eq4.1}) we obtain
$$R(P)\subseteq R(P^{\prime })+\mathcal{F}_{E}(M)^{\geq l+1}.$$
Hence $R(P)\subseteq R(P^{\prime })+\mathcal{F}_{E}(M)^{\geq l}$ for all positive integers $l$. This implies
that $R(P)$ is contained in the closure of $R(P^{\prime })$. Since $R(P^{\prime })$ is closed, it follows that $R(P)\subseteq R(P^{\prime })$. Then $P^{\prime }-P\subseteq R(P)^{2}\subseteq R(P^{\prime })^{2}$, hence $R(P^{\prime })\subseteq R(P)$. This completes the proof. 
\end{proof}

\begin{rem} \label{rem4.3} Let $\psi :M\rightarrow \mathcal{F}_{E}(M)$ be a morphism of $E$-$E$-bimodules, then
$\displaystyle \sum _{c\in T}\psi (c)\delta _{c^{*}}(P)$ is cyclically equivalent to $\displaystyle \sum _{a\in A}\psi (a)X_{a^{*}}(P)$.
Indeed, we have

$$\displaystyle \sum _{c\in T}\psi (c)\delta _{c^{*}}(P)=\sum _{a\in A}\sum _{s\in L(\tau (a))}s\psi (a)\delta _{(sa)^{*}}(P)$$
and this last expression is cyclically equivalent to $\displaystyle \sum _{a\in A }\psi (a)X_{a^{*}}(P)$.
\end{rem}
We will make use of the following fact.
 \begin{lemma} \label{lem4.4} Let $M$ be any $E$-$E$-bimodule, not necessarily $Z$-free. If $I$ is a closed ideal of $\mathcal{F}_{E}(M)$ and $J$ is the closure of the  ideal generated
 by $f_{1},\hdots, f_{N} \in \mathcal{F}_{E}(M)$, then any element $h\in IJ$ is cyclically equivalent to an element of the
 form $\alpha _{1}f_{1}+\hdots+\alpha _{N}f_{N}$, with $\alpha _{1},\hdots,\alpha _{N}\in I$.
 \end{lemma}
 
\begin{proof} The algebra $\mathcal{F}_{E}(M)$ is a $D$-algebra in the sense of Derksen-Weyman-Zelevinsky \cite{4}. Then our lemma follows from Lemma $13.8$ of \cite{4}.
\end{proof}
 
\begin{prop} \label{prop4.5} Suppose $P$ and $P^{\prime }$ are potentials in $\mathcal{F}_{E}(M)^{\geq 3}$ such that
 $P^{\prime }-P\in R(P)^{2}$. Then there is an algebra automorphism $\varphi $ of $\mathcal{F}_{E}(M)$, with $\varphi |_{E}=id_{E}$, such that $\varphi (P)$ is cyclically equivalent to $P^{\prime }$. Moreover, $\varphi (f)-f\in R(P)$ for all $f \in \mathcal{F}_{E}(M)$. 
 \end{prop}
   
\begin{proof} Here we follow closely the proof of Proposition $4.10$ from \cite{4}.
  
 {\bf Claim.} There is a sequence of $E$-$E$-bimodule morphisms
 $$\varphi ^{(2)}_{l}: M\rightarrow \mathcal{F}_{E}(M)\cap R(P)$$
 with $l \in \mathbb{Z}^{\geq 0}$ and $\varphi _{0}^{(2)}=0$, such that if we denote by $\varphi _{l}$ the automorphism
 of $\mathcal{F}_{E}(M)$ determined by the pair $(id_{M}, \varphi ^{(2)}_{l})$, then we have
 \begin{enumerate}[(i)]
 \item $\varphi ^{(2)}_{l}(a)\in \mathcal{F}_{E}(M)^{\geq l+1} \cap R(P)$ for all $l\in \mathbb{Z}^{\geq 0}$.
 \item $P^{\prime }$ is cyclically equivalent to
 $$\varphi _{0}\varphi _{1} \cdots \varphi _{l-1}(P)+\sum _{a\in A}\varphi ^{(2)} _{l}(a)X_{a^{*}}(P).$$ 
 \end{enumerate}  
 
 {\bf Proof of the Claim.} We construct the $E$-$E$-bimodule morphisms $\varphi _{l}^{(2)}$ by induction on $l$.
 Suppose first that $l=1$. Then Lemma \ref{lem4.4} gives that $P^{\prime }-P$ is cyclically equivalent to $\displaystyle \sum _{a\in A}f(a)X_{a^{*}}(P)$. Here $X_{a^{*}}(P)\in e_{\tau (a)}\mathcal{F}_{E}(M)e_{\sigma (a)}$ and
 $P^{\prime}-P$ is cyclic, therefore we may assume $ f(a)=e_{\sigma (a)}f(a)e_{\tau (a)}$. Since $A$
 is a set of $Z$-free generators of $M$, then there is a morphism of $E$-$E$-bimodules
 $$\varphi _{1}^{(2)}:M\rightarrow \mathcal{F}_{E}(M)^{\geq 2}$$
 such that $\varphi ^{(2)}_{1}(a)=e_{\sigma(a)}f(a)e_{\tau(a)} \in \mathcal{F}_{E}(M)^{\geq 2}\cap R(P)$. 
 
 Here $\varphi _{0}=id_{E}$, hence conditions (i) and (ii) are true for $l=1$. Suppose that for $l\geq 1$ we have already constructed the $E$-$E$-bimodule morphisms $\varphi _{0}^{(2)},\hdots,\varphi ^{(2)}_{l}$ satisfying conditions (i) and (ii).
 
 By Proposition \ref{prop2.14} and Lemma \ref{lem4.4} there exists $w\in [\mathcal{F}_{E}(M), \mathcal{F}_{E}(M)]\cap \mathcal{F}_{E}(M)^{\geq 3}$ such that
 $$\varphi _{l}(P)-P-\sum _{a\in A}\varphi _{l}^{(2)}(a)X_{a^{*}}(P)+w \in (\mathcal{F}_{E}(M)^{\geq 1}I)I$$
 where $I$ is the closure of the two-sided ideal generated by the elements $\varphi _{l}^{(2)}(a)$. Then by condition (ii) we
 have $I\subseteq \mathcal{F}_{E}(M)^{\geq l+1}(M)\cap R(P)$. Thus
 $$\varphi _{l}(P)-P-\sum _{a\in A}\varphi _{l}^{(2)}(a)X_{a^{*}}(P)+w \in (\mathcal{F}_{E}(M)^{\geq l+2}\cap (R(P)))R(P),$$
 and therefore $\varphi _{l}(P)-(P-w)\in R(P)^{2}$. Applying Proposition \ref{prop4.2} gives:
 $$R(\varphi _{l}(P))=R(P-w)=R(P).$$
From the above we obtain $\varphi _{l}(R(P))=R(\varphi _{l}(P))=R(P)$ and $R(P)=\varphi _{l}^{-1}(R(P))$. Then the element $\varphi _{l}(P)-P-\displaystyle \sum _{a\in A}\varphi _{l}^{(2)}(a)X_{a^{*}}(P)+w$  lies in
 $$ \varphi _{l}(\varphi _{l}^{-1}(\mathcal{F}_{E}(M)^{\geq l+2}\cap (R(P)))R(P))$$
 which is contained in 
 $$ \varphi _{l}((\mathcal{F}_{E}(M)^{\geq l+2}\cap (\varphi _{l}^{-1}(R(P))))\varphi _{l}^{-1}(R(P))).$$
 
 Therefore 
 \begin{equation} \label{eq4.2}
\varphi _{l}(P)-P-\sum _{a\in A}\varphi _{l}^{(2)}(a)X_{a^{*}}(P)+w =\varphi _{l}(z)
\end{equation}
 with $z\in (\mathcal{F}_{M}(E)^{\geq l+2}\cap R(P))R(P)$. By Lemma \ref{lem4.4}:
$$z=\sum _{a\in A}-h(a)X_{a^{*}}(P)+w_{1}$$ 
with $h(a)\in \mathcal{F}_{E}(M)^{\geq l+2}\cap R(P)$ and $w_{1}\in [\mathcal{F}_{E}(M),\mathcal{F}_{E}(M)]$. 
As in the case $l=1$ we may assume $h(a)=e_{\sigma(a)}h(a)e_{\tau (a)}$. Therefore there exists a morphism of
$E$-$E$-bimodules $\varphi ^{(2)}_{l+1}:M\rightarrow \mathcal{F}_{E}(M)^{\geq 2}$ with
$\varphi _{l+1}^{(2)}(a)=h(a)\in \mathcal{F}_{E}(M)^{\geq l+2}\cap R(P)$.

By (ii), $P^{\prime }$ is cyclically equivalent to 
$$\varphi _{0} \cdots \varphi _{l-1}\left(P+\sum _{a\in A}\varphi ^{(2)}_{l}(a)X_{a^{*}}(P)\right).$$

From Equation (\ref{eq4.2}) one has that the above element is cyclically equivalent to
$$\varphi _{0} \cdots \varphi _{l-1}\varphi _{l}\left(P+\sum _{a\in A}\varphi ^{(2)}_{l+1}(a)X_{a^{*}}(P)\right)$$
proving our claim. 

$  $

We now are ready to complete the proof of our Proposition. 

For each $l$, take $\rho _{l}=(\varphi _{l})^{-1}$. Then $\rho _{l}$ is determined by the pair $(id_{M},\rho _{l}^{(2)})$ where $\rho _{l}^{2}(a)=\displaystyle \sum _{k=1}^{\infty }(-\varphi ^{(2)}_{l}(a))^{k} \in \mathcal{F}_{E}(M)^{\geq l+1}.$  Now, define $\rho: \mathcal{F}_{E}(M) \rightarrow \mathcal{F}_{E}(M)$ as follows: for each $f \in \mathcal{F}_{E}(M)$, $\rho(f)=\displaystyle \lim_{l \to \infty} \rho_{l} \cdots \rho_{0}(f)$. By Proposition \ref{prop1.10}, $\rho $ is an automorphism of algebras fixing pointwise $E$. From (ii) of our previous claim we have
$$\rho _{l} \cdots \rho _{0}(P^{\prime })=P+w_{l}+z_{l}$$
with $w_{l}\in \mathcal{F}_{E}(M)^{\geq l+3}$ and $z_{l}\in [\mathcal{F}_{E}(M),\mathcal{F}_{E}(M)]$.

Then 
$$\rho (P^{\prime })(l)=P(l)+z_{l}(l).$$
Observe that if $l^{\prime }\leq l$, we have 
$$P(l^{\prime })+z_{l^{\prime }}(l^{\prime })=\rho _{l^{\prime }} \cdots \rho _{0}(P)(l^{\prime })=\rho _{l} \cdots \rho _{0}(P^{\prime })(l^{\prime })=P(l^{\prime })+z_{l}(l^{\prime }),$$
consequently $z_{l^{\prime}}(l^{\prime })=z_{l}(l^{\prime})$. Define $z\in \mathcal{F}_{E}(M)$ such that for any non negative integer $l$,
$z(l)=z_{l}(l)$. Then given $l$ and $l^{\prime }\leq l$, $(z-z_{l})(l^{\prime })=z_{l^{\prime }}(l^{\prime })-z_{l}(l^{\prime })=0$,
so $z-z_{l}\in \mathcal{F}_{E}(M)^{\geq l+1}$. 

Therefore $z=\displaystyle \lim_{l \to \infty }z_{l} \in [\mathcal{F}_{E}(M),\mathcal{F}_{E}(M)]$ and
$\rho (P^{\prime })=P+z$, so if $\varphi =\rho ^{-1}$ we obtain that
$P^{\prime }$ is cyclically equivalent to $\varphi (P)$. This completes the proof.
\end{proof}
$  $
  
\end{section}

\begin{section}{Splitting Theorem} \label{sec6}

In this section we take $E=\displaystyle \prod_{i=1}^{n} E_{i}$, where each $E_{i}$ is a semisimple finite-dimensional $F$-algebra, such that $E\otimes _{F}E^{op}$ is a semisimple $F$-algebra. The $E$-$E$-bimodules $M$ we consider here are acyclic; that is, $M_{cyc}=0$, but they are  not necessarily $Z$-free. 

\begin{definition} \label{def5.1} If $M$ is any $E$-$E$-bimodule, a polarization of $M$ is a set $\mathcal{A}\subseteq \{1,\ldots,n\}\times \{1,\ldots,n\}$ such that
\begin{enumerate}
\item[(i)] If $(i,j)\in \mathcal{A}$, then $e_{i}Me_{j}\neq 0$ and $(j,i)$ is not in $\mathcal{A}$.
\item[(ii)] If $e_{i}Me_{j}\neq 0$, then $(i,j)\in \mathcal{A}$ or $(j,i)\in \mathcal{A}$.
\end{enumerate}
\end{definition}
Clearly, if $M$ is any acyclic $E$-$E$-bimodule then there is at least one polarization of $M$.

\begin{definition} \label{def5.2} Suppose that $N$ is an $E$-$E$-bimodule, $U$ is a set of directed generators of $N$, and $U(i,j)=U \cap e_{i}Ne_{j}$. We say that $N$ is $2$-cycle complete with respect to $U$, if there is a polarization $\mathcal{A}$ of $N$ and for each $(i,j)\in \mathcal{A}$ there is a bijection $\hat{-}:U(i,j)\rightarrow U(j,i)$. In this case we define $U(\mathcal{A})=\displaystyle \bigcup _{(i,j)\in \mathcal{A}}U(i,j)$.
\end{definition}

$ $
Observe that if $N$ is a $2$-cycle complete $E$-$E$ bimodule with respect to $U$, then
$$N=\displaystyle \bigoplus _{a\in U(\mathcal{A})}(EaE\oplus E\hat{a}E)$$

$ $

Using the above notation we associate to a $2$-cycle complete $E$-$E$-bimodule $N$ a quadratic potential
$$Q_{N}=\sum _{a \in U(\mathcal{A})}a\hat{a}.$$
 
 $  $
 
 \begin{theorem} \label{theo5.3} (Splitting Theorem) Let $M$ be an $E$-$E$-bimodule with a decomposition $M=N\oplus M^{\prime }$ as $E$-$E$-bimodules where $N$ is $2$-cycle complete with respect to $U$, a set of directed generators for $N$. Let $\mathcal{L}$ be the closure of the two-sided ideal in $\mathcal{F}_{E}(M)$ generated by $N$. Suppose that for any positive integer $l$ and any potential in $\mathcal{L}\cap \mathcal{F}_{E}(M)^{\geq l}$ of the form $\displaystyle \sum _{(i,j)\in \mathcal{A}}\displaystyle \sum _{a \in U(i,j)}(a u_{a }+v_{\hat{a}}\hat{a })+z$, with $z\in [\mathcal{F}_{E}(M),\mathcal{F}_{E}(M)]\cap \mathcal{F}_{E}(M)^{\geq l}$, there is a morphism of $E$-$E$-bimodules
 $$\rho :N\rightarrow \mathcal{F}_{E}(M)$$ 
with $\rho (a )=-v_{\hat{a}}$ and $\rho (\hat{a})=-u_{a}$.
Then for any potential of the form $P=Q_{N}+P^{\geq 3 }$ with $P^{\geq 3 }\in \mathcal{F}_{E}(M)^{\geq 3}$, there exists an algebra automorphism $\varphi $ of $\mathcal{F}_{E}(M)$, with $\varphi |_{E}=id_{E}$, such that $\varphi (P)$ is cyclically equivalent to a potential $Q_{N}+Q$ with $Q\in \mathcal{F}_{E}(M^{\prime })^{\geq 3}$. 
\end{theorem}
   
\begin{proof} We follow the proof of the Splitting Theorem given in \cite{4}. First we prove the following:
 
 {\bf Claim 1}. Suppose $P$ is a potential such that 
   $$P=Q_{N}+P^{\prime }+H$$
   with $P^{\prime }$ a potential in $\mathcal{F}_{E}(M^{\prime })$ and $H$ a potential in $\mathcal{F}_{E}(M)^{\geq d+2}$. 
   
 Then there exists an algebra automorphism $\varphi $ of $\mathcal{F}_{E}(M)$ determined by a pair of morphisms of $E$-$E$-bimodules  $(id_{M},\varphi ^{(2)})$ such that $\varphi ^{(2)}(M^{\prime })=0, \varphi ^{(2)}(M)\in \mathcal{F}_{E}(M)^{\geq d+1}$ and
   $$\varphi (P)=Q_{N}+P^{\prime }+P^{\prime \prime }+H^{\prime }+z^{\prime }$$
   where  $P^{\prime \prime }$ is a potential in  $ \mathcal{F}_{E}(M)^{\geq 2d+2}$,
   $ H^{\prime }$ is a potential in $ \mathcal{F}_{E}(M)^{\geq 2d+2}$ and $z^{\prime }\in [\mathcal{F}_{E}(M),\mathcal{F}_{E}(M)]\cap \mathcal{F}_{E}(M)^{\geq d+2}$.
   
   {\bf Proof of the Claim 1}. We recall that  $U(\mathcal{A})=\displaystyle \bigcup _{(i,j)\in \mathcal{A}}U(i,j)$. We have $H=G+P^{''}$ with $G\in \mathcal{L}\cap \mathcal{F}_{E}(M)^{\geq d+2}$, $P^{''}\in \mathcal{F}_{E}(M^{\prime })^{\geq d+2}$. By Lemma \ref{lem4.4}, $G=\displaystyle \sum _{a\in U(\mathcal{A})}(au_{a}+v_{\hat{a}}\hat{a})+z_{1}$, with $z_{1}\in [\mathcal{F}_{E}(M),\mathcal{F}_{E}(M)]\cap \mathcal{F}_{E}(M)^{\geq d+2}$. By assumption, there exists a morphism
of $E$-$E$-bimodules $\varphi ^{(2)}:M\rightarrow \mathcal{F}_{E}(M)^{\geq d+1}$, such that $\varphi ^{(2)}(a)=-v_{\hat{a}}, \varphi ^{(2)}(\hat{a})=-u_{a}$, and $\varphi ^{(2)}(M^{\prime })=0$. Since $u_{a}$ and $v_{\hat{a}}$ are in $\mathcal{F}_{E}(M)^{\geq d+1}$, then $\varphi ^{(2)}(M)\subseteq \mathcal{F}_{E}(M)^{\geq d+1}$.  Let $\varphi$ be the algebra automorphism of $\mathcal{F}_{E}(M)$ determined by the pair $(id_{M},\varphi^{(2)})$. Observe that the restriction of $\varphi$ to $M^{\prime}$ is the identity map; consequently, $\varphi $ is the identity on $\mathcal{F}_{E}(M^{\prime })$. Moreover,  $\varphi (m_{1} \cdots m_{l})=m_{1} \cdots m_{l}+\mu $ where $m_{1},\dots,m_{l}\in M$ and $\mu \in \mathcal{F}_{E}(M)^{\geq l+d}$. Then for $f\in \mathcal{F}_{E}(M)^{\geq l}$,  $\varphi (f)=f+f^{\prime }$ with  $f^{\prime }\in \mathcal{F}_{E}(M)^{\geq l+d}$. We have
   
$$\varphi (P)=\varphi (Q_{N})+P'+P^{''}+\varphi(G)$$ 
$$\varphi (Q_{N})+\varphi (G)=\sum _{a\in U(\mathcal{A})}(a-v_{\hat{a}})(\hat{a}-u_{a})+(a-v_{\hat{a}})(u_{a}+u_{a}^{\prime })+(v_{\hat{a}}+v_{\hat{a}}^{\prime })(\hat{a}-u_{a})+\varphi (z_{1})$$
$$=\sum _{a\in U(\mathcal{A})}(a\hat{a}+au_{a}^{\prime}-v_{\hat{a}}u^{\prime }_{a}-v_{\hat{a}}u_{a}+v^{\prime }_{\hat{a}}\hat{a}-v^{\prime }_{\hat{a}}u_{a})+z^{\prime }$$
where $u^{\prime }_{a}, v^{\prime }_{\hat{a}}\in \mathcal{F}_{E}(M)^{\geq 2d+1}$ and $z^{\prime }=\varphi (z_{1})\in [\mathcal{F}_{E}(M), \mathcal{F}_{E}(M)]\cap \mathcal{F}_{E}(M)^{\geq d+2}$. Therefore
$$\varphi (Q_{N})+\varphi (G)=\sum _{a\in U(\mathcal{A})}a\hat{a}+H^{\prime }+z^{\prime }$$
with $H^{\prime }\in \mathcal{F}_{E}(M)^{\geq 2d+2}$. We have
$$\varphi(P)=Q_{N}+P'+P^{''}+H'+z'$$
From here we obtain our claim.

{\bf Claim 2}. There is a sequence $\{\varphi _{k}\}_{k=0}^{\infty} $ of algebra automorphisms of $\mathcal{F}_{E}(M)$ where
each $\varphi _{k}$ is determined by a pair of morphisms $(id_{M},\varphi _{k}^{(2)})$ with the following
properties
\begin{enumerate}[(i)]
\item $\varphi _{l}^{(2)}(M^{\prime })=0$ for all non-negative integers $l$. 
\item $\varphi _{0}^{(2)}=0$ and $\varphi _{k}^{(2)}(M)\subseteq \mathcal{F}_{E}(M)^{\geq 2^{k}+1}$ for all $k$.
\item There exist potentials $P_{0},\dots,P_{k} \in \mathcal{F}_{E}(M^{\prime })$ and elements $z_{0}=0 ,z_{1},\dots,z_{k}$ lying in $ [\mathcal{F}_{E}(M), \mathcal{F}_{E}(M)]$ with $z_{i}\in \mathcal{F}_{E}(M)^{\geq 2^{i-1}+2},  P_{i}\in \mathcal{F}_{E}(M)^{\geq 2^{i-1}+2}$ such that
$$\varphi _{k}\varphi _{k-1}\cdots \varphi _{0}(P)=$$
$$Q_{N}+P_{0}+\ldots+P_{k}+H_{k}+\varphi _{k-1}\cdots \varphi _{1}(z_{1})+
\varphi _{k}\cdots \varphi _{2}(z_{1})+\ldots+\varphi _{k}(z_{k-1})+z_{k}$$
where $H_{k}$ is a potential in $\mathcal{F}_{E}(M)^{\geq 2^{k}+2}\cap \mathcal{L}$.
  \end{enumerate}
  
  $ $
  
  {\bf Proof of Claim 2}. By induction on $k$. If $k=0$, we have
  $$\varphi _{0}(P)=P=Q_{N}+P^{\geq 3}.$$
  Here $\mathcal{F}_{E}(M)=\mathcal{F}_{E}(M^{\prime })\oplus \mathcal{L}$, then
  $$P=Q_{N}+P_{0}+H_{0}$$
  with $P_{0}\in \mathcal{F}_{E}(M^{\prime })^{\geq 3}$ and $H_{0}\in \mathcal{F}_{E}(M)^{\geq 3}\cap \mathcal{L}$. Since $3=2^{0}+2$, we obtain our result for $k=0$.
  
 Suppose now our result is true for $k$, let us prove it for $k+1$. We have
$$\varphi _{k} \cdots \varphi _{0}(P)=Q_{N}+P^{\prime }+H_{k}+w_{k}$$
with $P^{\prime }=P_{0}+\ldots+P_{k}$, $H_{k}\in \mathcal{F}_{E}(M)^{\geq 2^{k}+2}\cap \mathcal{L}$, and 
$$w_{k}=  \varphi _{k} \cdots \varphi _{2}(z_{1})+ \varphi _{k} \cdots \varphi _{3}(z_{2})+\ldots+\varphi _{k}(z_{k-1})+z_{k}.$$

Applying Claim $1$ with $d=2^{k}$, yields an algebra automorphism $\varphi _{k+1}$ of $\mathcal{F}_{E}(M)$ determined by a pair of  $E$-$E$-bimodule morphisms of the form $(id_{M},\varphi _{k+1}^{(2)})$ 
with $\varphi _{k+1}^{(2)}(M^{\prime })=0$ and $\varphi _{k+1}^{(2)}(M)\subseteq \mathcal{F}_{E}(M)^{\geq 2^{k}+1}$ such that

$$ \varphi _{k+1}\varphi _{k} \cdots \varphi _{0}(P)=Q_{N}+P^{\prime \prime }+P^{\prime }+H_{k+1}+z_{k+1}+\varphi _{k+1}(w_{k}),$$
with $P^{\prime \prime }\in \mathcal{F}_{E}(M^{\prime })^{\geq 2^{k}+2}$, $H_{k+1}\in \mathcal{F}_{E}(M)^{\geq 2^{k+1}+2}\cap \mathcal{L}$ and $z_{k+1}$ lies in $[\mathcal{F}_{E}(M),\mathcal{F}_{E}(M)] \cap \mathcal{F}_{E}(M)^{\geq 2^{k}+2}$. This proves our claim.

Now we are ready to finish the proof of our proposition. We have $\varphi _{k}^{(2)}(M)\subseteq \mathcal{F}_{E}(M)^{\geq 2^{k}+1}\subseteq \mathcal{F}_{E}(M)^{\geq k+1}$. Hence we may apply Proposition \ref{prop1.10}, and construct an algebra automorphism $\varphi $ of $\mathcal{F}_{E}(M)$ defined as
$$\varphi (P)=\displaystyle \lim_{k\rightarrow \infty } \varphi _{k} \cdots \varphi _{0}(P).$$

Observe that the sequence $\{w_{k}\}_{k=0}^{\infty}$ converges. Indeed, we have
$w_{k+1}=z_{k+1}+\varphi _{k+1}(w_{k})$. Now $\varphi _{k+1}(w_{k})=w_{k}+f_{k}$, where
$f_{k}, z_{k+1}\in \mathcal{F}_{E}(M)^{\geq 2^{k}+2}$. Therefore
$$w_{k+1}-w_{k}=z_{k+1}+f_{k}\in \mathcal{F}_{E}(M)^{\geq 2^{k}+2}$$
so $\{w_{k}\}_{k=0}^{\infty }$ is a Cauchy sequence, hence converges. Clearly the sequence $\{ P_{k} \}_{k=0}^{\infty }$ is summable, so
$$\varphi (P)=Q_{N}+\sum _{k=0}^{\infty }P_{k}+\displaystyle \lim_{k\rightarrow \infty }H_{k}+\displaystyle \lim_{k\rightarrow \infty }w_{k}$$
Therefore $$\varphi(P)=Q_{N}+Q+w,$$
where $Q=\displaystyle \sum _{k=0}^{\infty }P_{k}\in \mathcal{F}_{E}(M^{\prime })$ and $w=\displaystyle \lim_{k\rightarrow \infty }w_{k}\in [\mathcal{F}_{E}(M),\mathcal{F}_{E}(M)]$. Observe that $\displaystyle \lim_{k\rightarrow \infty }H_{k}=0$.

This proves our proposition.
\end{proof}

\begin{definition} \label{def5.4} Let $N$-be a $2$-cycle complete $E$-$E$-bimodule with respect to $U$, a set of directed generators of $N$. We say that the quadratic potential $Q_{N}=\displaystyle \sum _{a\in U(\mathcal{A})}a\hat{a}$ has the splitting property if for any $E$-$E$-bimodule $M^{\prime }$ and any potential $P$ in $\mathcal{F}_{E}(N\oplus M^{\prime })^{\geq 3}$ there exists an algebra automorphism $\varphi $ of $\mathcal{F}_{E}(N\oplus M^{\prime})$ such that $\varphi (Q_{N}+P)$ is cyclically equivalent to $Q_{N}+Q$, where $Q$ is a potential in $\mathcal{F}_{E}(M^{\prime })^{\geq 3}$.
\end{definition}

\begin{corollary} \label{coro5.5} If $N$ is a $2$-cycle complete $E$-$E$-bimodule with respect to $U$, a set of directed $Z$-free generators
of $N$, then $Q_{N}$ has the splitting property.
\end{corollary}

\begin{proof} Let $M^{\prime }$ be any $E$-$E$-bimodule and let $\mathcal{L}$ be the closure of the two sided-ideal in $\mathcal{F}_{E}(N\oplus M^{\prime })$ generated by $N$. We have $\mathcal{F}_{E}(N\oplus M^{\prime })\mathcal{L}=\mathcal{L}$. The ideal
$\mathcal{L}$ is generated by the elements of $U$. Therefore, by Lemma \ref{lem4.4}, any potential
$H$ in $\mathcal{L}$ is cyclically equivalent to an element of the form $\displaystyle \sum _{a\in U(\mathcal{A})}(u_{a}a+v_{\hat{a}}\hat{a})$, and the latter is cyclically equivalent to
$\displaystyle \sum _{a\in U(\mathcal{A})}(au_{a}+ v_{\hat{a}}\hat{a})$.

We have $au_{a}=ae_{\sigma (a)}u_{a}$ is a cyclic element, so we may assume  $u_{a}$ is directed; similarly, we may assume $v_{\hat{a}}$ is directed. We have $\tau (u_{a})=\sigma (a)=\tau (\hat{a})$ and $\sigma (u_{a})=\tau (a)=\sigma (\hat{a})$. Similarly,
$\tau (v_{\hat{a}})=\tau (a)$ and $\sigma (v_{\hat{a}})=\sigma (a)$. Since the elements of $U$ are $Z$-free generators there is a morphism of $E$-$E$-bimodules $\psi :N\rightarrow \mathcal{F}_{E}(N\oplus M^{\prime })$ such
that $\psi (a)=-v_{\hat{a}}$ and $\psi (\hat{a})=-u_{a}$ for all $a\in U(\mathcal{A})$. Then our corollary follows from Theorem \ref{theo5.3}. 
\end{proof}

Here $E\otimes _{F}E^{op}$ is a semisimple algebra. Then for any pair $i,j$, $E_{i}\otimes _{F}E_{j}^{op}$ is also a semisimple algebra. For the identity  of $E_{i}\otimes _{F}E_{j}^{op}$, we have a decomposition into a sum of orthogonal primitive idempotents
$$1_{E_{i}}\otimes 1_{E_{j}^{op}}=\sum _{\epsilon\in  I(i,j)}\epsilon .$$

Each $E_{i}\otimes _{F}E_{j}^{op}$-left module $(E_{i}\otimes E_{j}^{op})\epsilon $ is a simple module. Hence if $W$ is an $E_{i}$-$E_{j}$-bimodule it is also an $E_{i}\otimes _{F}E_{j}^{op}$-left module where for $x\in E_{i},y\in E_{j}$ and  $w\in W$, $(x\otimes y)* w=xwy$.  In case $w=\epsilon *w\neq 0$, $E_{i}wE_{j}$ is a simple $E_{i}$-$E_{j}$-bimodule. Moreover, there is an isomorphism of $E_{i}$-$E_{j}$-bimodules from $E_{i}\epsilon E_{j}$ to $E_{i}wE_{j}$ sending $\epsilon $ to $w$. Therefore if $W^{\prime }$ is another $E_{i}$-$E_{j}$-bimodule and $w^{\prime }=\epsilon *w^{\prime }\neq 0$, then there is an isomorphism
of $E_{i}$-$E_{j}$-bimodules sending $w$ to $w^{\prime }$.

Observe that we have an antiisomorphism $\mathfrak{s}:E_{i}\otimes_{F} E_{j}^{op}\rightarrow E_{j}\otimes _{F}E_{i}^{op}$
with $\mathfrak{s}(x\otimes y)=y\otimes x$. Then 
$$1_{E_{j}}\otimes 1_{E_{i}}=\mathfrak{s}(1_{E_{i}}\otimes 1_{E_{j}})=\sum _{\epsilon \in I(i,j)}\mathfrak{s}(\epsilon)$$
gives a decomposition of the identity in the algebra $E_{j}\otimes _{F}E_{i}^{op}$ into a sum of orthogonal primitive idempotents.

\begin{corollary} \label{coro5.6} Let $N$ be a $2$-cycle complete $E$-$E$-bimodule with respect to a set $U$ which has the following
properties
\begin{enumerate}[(i)]
\item If $(i,j)\in \mathcal{A}$ and $a\in U(i,j)$, then $a=\epsilon _{a} *a$ for some idempotent $\epsilon _{a}\in I(i,j)$.
\item If $(i,j)\in \mathcal{A}$ and $a\in U(i,j)$, then $\hat{a}=\mathfrak{s}(\epsilon _{a})*\hat{a}$.
\end{enumerate}
Then $Q_{N}=\displaystyle \sum _{a\in U(\mathcal{A})}a\hat{a}$ has the splitting property.
\end{corollary}

\begin{proof} Take $M^{\prime }$ an $E$-$E$-bimodule and let $M=N\oplus M^{\prime}$.  Let $\mathcal{L}$ be the closure of the two-sided ideal in $\mathcal{F}_{E}(M)$ generated by $N$. Then if $H$ is a potential in
$\mathcal{L}$, we have as in the proof of Corollary \ref{coro5.5} that $H$ is cyclically equivalent to an element of the form $\displaystyle \sum _{a\in U(\mathcal{A})}(au_{a}+v_{\hat{a}}\hat{a})$.

We have $au_{a}=\epsilon _{a}*au_{a}$ and $\epsilon _{a}=\displaystyle \sum _{u}x_{u}\otimes y_{u}$ for some $x_{u}\in E_{i},
y_{u}\in E_{j}$. Then $au_{a}=\displaystyle \sum _{u}x_{u}ay_{u}u_{a}$ and this element is cyclically equivalent to
$\displaystyle \sum _{u}ay_{u}u_{a}x_{u}=a(\mathfrak{s}(\epsilon _{a})*u_{a})$. Similarly, the element
$v_{\hat{a}}\hat{a}=v_{\hat{a}}(\mathfrak{s}(\epsilon _{a})*\hat{a})$ is cyclically equivalent to the element
$(\epsilon _{a}*v_{\hat{a}})\hat{a}$. Therefore any potential $H\in \mathcal{L}$ is cyclically equivalent to an element
of the form $\displaystyle \sum _{a\in U(\mathcal{A})}(a(\mathfrak{s}(\epsilon _{a})*u_{a})+(\epsilon _{a}*v_{\hat{a}})\hat{a})$. Moreover, there is a morphism of $E$-$E$-bimodules $\psi :N\rightarrow \mathcal{F}_{E}(M)$ such that
$$\psi (a)=\psi (\epsilon _{a}*a)=\epsilon _{a}*v_{\hat{a}}$$ 
$$\psi (\hat{a})=\psi (\mathfrak{s}(\epsilon _{a}*\hat{a}))=\mathfrak{s}(\epsilon _{a})*u_{a}.$$
The result now follows from Theorem \ref{theo5.3}.
\end{proof}
\end{section}
\begin{section}{Local Bases and cyclic derivations} \label{sec7}
Here we consider semisimple algebras $E=\displaystyle \prod _{i=1}^{n}E_{i}$, such that the enveloping algebra
$E\otimes _{F}E^{op}$ is semisimple and  each $E_{i}$ has an $F$-basis satisfying
some properties given below. These properties are trivially satisfied in case each $E_{i}$ equals the underlying base field $F$. In most of the related papers about algebras with potentials (\cite{4,6,8}), these conditions are satisfied. Under these conditions on $E$, we find a cyclic derivation on $\mathcal{F}_{E}(M)$ for any $E$-$E$-bimodule $M$. Using this cyclic derivation we can define the Jacobian ideal $J(P)$ for any potential $P$ in $\mathcal{F}_{E}(M)$. This ideal coincides with the ones considered in \cite{3,6,8}. Finally, we study quadratic potentials and their Jacobian ideals.

\begin{definition} \label{def6.1} A subset $L \subseteq E$ is a local basis for $E$ if $L$ is an $F$-basis consisting of directed elements. 
\end{definition}

If $L$ is a local basis for $E$, then each $L_{i}=Ee_{i}\cap L$ is an $F$-basis for $Ee_{i}$ and $L=\displaystyle \bigcup _{i=1}^{n}L_{i}$. For each $t\in L$, we define
$t^{*}:E\rightarrow Z$ as the $Z$-linear map such that $t^{*}(t)=e_{i}$ and $t^{*}(s)=0$ for $s\in L \setminus \{t\}$.  \\

Throughout the rest of the paper we will assume that the local bases of $E$ satisfy the following properties: 
\begin{enumerate}
\item All elements of $L$ are invertible.
\item If $c(i)=\mathrm{dim}_{F}e_{i}E$, then $\operatorname{char} F\nmid c(i)$.
\end{enumerate}
and also the following condition
\begin{equation} \label{eq6.1}
e_{i} \in L_{i} \ \text{and if} \ s,t\in L_{i}, e_{i}^{*}(st^{-1})\neq 0 \ \text{implies} \ s=t; \text{likewise}, e_{i}^{*}(t^{-1}s)\neq 0 \ \text{implies} \ s=t.
\end{equation}

\begin{example} \label{ex6.2} (Demonet, \cite{3}) Take $E=\displaystyle \prod _{i=1}^{n}F[G_{i}]$, the finite direct product of group algebras of a finite group $G_{i}$, whose order
is not divisible by the characteristic of the base field $F$. For each $i \in \{1,\ldots,n\}$, take $L_{i}=G(i)$ and $L=\displaystyle \bigcup _{i=1}^{n}L_{i}$. 
\end{example}

\begin{example} \label{ex6.3} Let $E=\displaystyle \prod _{i=1}^{n}D_{i}$, where each $D_{i}$ is a division ring over $F$, and suppose that $D_{i}$ has a semi-multiplicative basis $L_{i}$: this means that if $s,t\in L_{i}$, then $st=\lambda w$ with $w\in L_{i}$ and $\lambda \in F$. Moreover, $\mathrm{dim}_{F}D_{i}$ is not divisible by the characteristic of the field $F$.
\end{example}

\begin{example} \label{ex6.4} (Zelevinsky and Labardini-Fragoso, \cite{8}) Let $G \supseteq F$ a Galois extension such that $\mathrm{Gal}(G/F)$ is a cyclic group. Take $E=\displaystyle \prod _{i=1}^{n}E_{i}$ where each $E_{i}$ is an intermediate field of $G/F$.
\end{example}

\begin{example} \label{ex6.5} Take $E=\displaystyle \prod _{i=1}^{n}F_{i}$, where each $F_{i}$ is an abelian extension of the field $F$, containing a $c(i)$-primitive root of unity.
\end{example}

\begin{rem} \label{rem6.6} Observe that in Definition $6.1$ the set consisting of all the elements $r^{-1}$, with $r \in L$, is also a $Z$-local basis for $L$. Then in the rest of the paper, for $r \in L_{i}$, $(r^{-1})^{\ast}: E \rightarrow Z$ is the $Z$-linear map such that $(r^{-1})^{\ast}(r^{-1})=e_{i}$, and for $t \in L \setminus \{r\}$, $(r^{-1})^{\ast}(t)=0$. 
\end{rem}

\begin{prop} \label{prop6.7} For $s, t, t_{1}\in L_{i}$ one has
$$\sum _{r\in L_{i}}(r^{-1})^{*}(t_{1}^{-1}s^{-1})r^{*}(st)=\delta _{t,t_{1}}.$$
\end{prop}

\begin{proof} We have
$$st=\sum _{r\in L_{i}}r^{*}(st)r,$$
and 
$$t_{1}^{-1}s^{-1}=\sum _{r_{1}\in L_{i}}(r_{1}^{-1})^{*}(t_{1}^{-1}s^{-1})r_{1}^{-1}.$$
Multiplying the above equalities yields
$$t_{1}^{-1}t=\sum _{r_{1},r\in L_{i}}(r_{1}^{-1})^{*}(t_{1}^{-1}s^{-1})r^{*}(st)r_{1}^{-1}r.$$
Applying $e_{i}^{*}$ to the above equality and using (\ref{eq6.1}) we get

$$\delta _{t,t_{1}}=\sum _{r_{1},r\in L_{i}}(r_{1}^{-1})^{*}(t_{1}^{-1}s^{-1})r^{*}(st)e_{i}^{*}(r_{1}^{-1}r)$$
$$=\sum _{r\in L_{i}}(r^{-1})^{*}(t_{1}^{-1}s^{-1})r^{*}(st).$$
\end{proof}

This proves our result.

\begin{corollary} \label{coro6.8} For each $r,r_{1},s\in L_{i}$ we have
$$\sum _{t\in L_{i}}r^{*}(st)(r_{1}^{-1})^{*}(t^{-1}s^{-1})=\delta _{r,r_{1}}.$$
\end{corollary}

\begin{proof} Consider the square matrices of order $c(i)$, $B=(b_{t_{1},r})$ and $A=(a_{r,t})$, defined as $b_{t_{1},r}=(r^{-1})^{*}(t_{1}^{-1}s^{-1})$ and $a_{r,t}=r^{*}(st)$.  Applying Proposition \ref{prop6.7} yields $BA=I$, where $I$ is the identity matrix of order $c(i)$. Therefore  $AB=I$, and the result follows.
\end{proof}

\begin{prop} \label{prop6.9} For $s_{1}, s, t \in L_{i}$ we have the following equalities
$$\sum _{r\in L_{i}}(r^{-1})^{*}(t^{-1}s_{1}^{-1})r^{*}(st)=\delta _{s_{1},s},$$
$$\sum _{s\in L_{i}}r^{*}(st)(r_{1}^{-1})^{*}(t^{-1}s^{-1})=\delta _{r,r_{1}}.$$
 \end{prop}
 
\begin{proof} Similar to the proofs of Proposition $\ref{prop6.7}$ and Corollary $\ref{coro6.8}$.
 \end{proof}
 
\begin{lemma} \label{lem6.10} If $f\in e_{i}\mathcal{F}_{E}(M)e_{i}$ and $t\in E(i)$, then
 \begin{enumerate}[(a)]
\item $t \left(\displaystyle \sum _{\omega \in L_{i}}\omega ^{-1}f\omega \right)=\left(\displaystyle \sum _{\omega \in L_{i}}\omega ^{-1}f\omega \right)t.$
\item $\displaystyle \sum _{\omega \in L_{i}}\omega ^{-1}tf\omega =\sum _{\omega \in L_{i}}\omega ^{-1}ft\omega.$ 
 \end{enumerate}
 \end{lemma}
 $  $
 
 \begin{proof} Since the set of all elements $s^{-1}$, with $s\in L_{i}$, forms an $F$-basis for $E(i)$, it suffices to prove the statement
 for $t=s^{-1}$ with $s\in L_{i}$. We have
 $$s^{-1}\left(\sum _{\omega \in L_{i}}\omega ^{-1}f\omega \right)=\sum _{r,\omega \in L_{i}}(r^{-1})^{*}(s^{-1}\omega ^{-1})r^{-1}f\omega ss^{-1}$$
 $$=\sum _{r, r_{1}, \omega \in L_{i}}r^{-1}f(r^{-1})^{*}(s^{-1}\omega ^{-1})r_{1}^{*}(\omega s)r_{1}s^{-1}$$
 $$=\sum _{r,r_{1}}f\left(\sum _{\omega \in L_{i}}(r^{-1})^{*}(s^{-1}\omega ^{-1})r_{1}^{*}(\omega s)\right)r_{1}s^{-1}$$
 $$=\sum _{r,r_{1}\in L_{i}}r^{-1}f\delta _{r,r_{1}}r_{1}s^{-1}=\left(\sum _{r\in L_{i}}r^{-1}fr\right)s^{-1}.$$
 
 This proves (a). A similar reasoning proves (b). 
 \end{proof}
 
With our assumptions on the local bases $L$ for $E$, we can now define a new cyclic derivation on $\mathcal{F}_{E}(M)$, where $M$ is a $Z$-free $E$-$E$-bimodule. Once done this, then we can define a cyclic derivation in $\mathcal{F}_{E}(M)$ for any $E$-$E$-bimodule $M$.
 
\begin{definition} \label{def6.11} Let $M$ be an $E$-$E$-bimodule which is $Z$-free. Let $f$ and $g$ be elements of $\mathcal{F}_{E}(M)$. We define
$$\widetilde{\mathfrak{h}}(f)(g)=\sum _{s\in L}s^{-1}\mathfrak{h}(f)(g)s.$$
\end{definition}

Clearly $\widetilde{\mathfrak{h}}:\mathcal{F}_{E}(M)\rightarrow \mathrm{End}_{F}(\mathcal{F}_{E}(M))$ is 
 a cyclic derivation. Then we have the corresponding cyclic derivative $\widetilde{\delta }(f)=\widetilde{\mathfrak{h}}(f)(1)$, and for $\psi \in M^{*}$ the cyclic derivative with respect to $\psi $: $\widetilde{\delta }_{\psi}=\psi \widetilde{\delta }$.
  
\begin{prop} \label{prop6.12}  Let $M$ be a $Z$-free $E$-$E$-bimodule, then 
\begin{enumerate}
\item For $f,g\in \mathcal{F}_{E}(M)$ and $t\in E$ we have
$$ \widetilde{\mathfrak{h}}(tf)(g)=\widetilde{\mathfrak{h}}(f)(gt)$$
and
$$ \widetilde{\mathfrak{h}}(ft)(g)=\widetilde{\mathfrak{h}}(f)(tg).$$
\item If $f,g$  are directed elements of $ \mathcal{F}_{E}(M)$ and $\widetilde{\mathfrak{h}}(f)(g)\neq 0$, then
$\sigma (f)=\tau (g)$ and $\tau (f)= \sigma (g)$; thus $fg$ is cyclic.
\item If $m_{1},\dots,m_{l}$ are directed elements of $M$ with $m_{1}\cdots m_{l}\neq 0$ and $g\in \mathcal{F}_{E}(M)$, then
$$\widetilde{\mathfrak{h}}(m_{1} \cdots m_{l})(g)=\sum _{s\in L}s^{-1}(m_{1} \cdots m_{l}g+m_{2} \cdots m_{l}gm_{1}+\ldots+m_{l}gm_{1}\cdots m_{l-1})s$$
\end{enumerate}  
\end{prop}
 
\begin{proof}
\begin{enumerate}
\item Follows from Lemmas \ref{lem2.5} and  \ref{lem6.10}.
\item Follows from item $1$.
\item We may assume $m_{i}=x_{i}t_{i}$ with $x_{i}\in EM_{0}, t_{i}\in E$ for $i=1,\dots,l$ and 
$g=e_{\sigma (m_{l})}ge_{\tau (m_{1})}$.
Then using (c) of Proposition \ref{prop2.8} 
\begin{align*}
&\widetilde{\mathfrak{h}}(m_{1}\cdots m_{l})(g)=\widetilde{\mathfrak{h}}(x_{1}(t_{1}x_{2})\cdots (t_{l-1}x_{l}))(t_{l}g) \\
&=\sum _{s\in L}s^{-1}(x_{1}(t_{1}x_{2}) \cdots (t_{l-1}x_{l})(t_{l}g)+(t_{1}x_{2})\cdots (t_{l-1}x_{l})(t_{l}g)x_{1}+
\ldots+(t_{l-1}x_{l})(t_{l}g)\cdots (t_{l-2}x_{l-1}))s \\
&=\sum _{s\in L}s^{-1}((x_{1}t_{1})(x_{2}t_{2})\cdots (x_{l}t_{l})g+(x_{2}t_{2})\cdots (x_{l}t_{l})g(x_{1}t_{1})+\ldots+(x_{l})(t_{l}g)(x_{1}t_{1})\cdots (x_{l-1}t_{l-1}))s \\
&=\sum _{s\in L}s^{-1}(m_{1}\cdots m_{l}g+m_{2}\cdots m_{l}gm_{1}+\ldots+m_{l}gm_{1}\cdots m_{l-1})s.
\end{align*}
 \end{enumerate}
This completes the proof.
\end{proof}

\begin{prop} \label{prop6.13} Suppose  $M$ is any $E$-$E$-bimodule, then there is a cyclic derivation $\widetilde{\mathfrak{h}}:\mathcal{F}_{E}(M)\rightarrow \mathrm{End}_{F}(\mathcal{F}_{E}(M))$ such that for $m_{1},\dots,m_{l}$ directed elements of $M$ with $m_{1}\cdots m_{l}\neq 0$, we have
$$\widetilde{\mathfrak{h}}(m_{1}\cdots m_{l})(g)=\sum _{s\in L}s^{-1}(m_{1}\cdots m_{l}g+m_{2}\cdots m_{l}gm_{1}+\ldots+m_{l}gm_{1}\cdots m_{l-1})s.$$
\end{prop}

\begin{proof} Since $E\otimes _{F}E^{op}$ is semisimple then there is an inclusion $M\rightarrow {\bf M}$ where ${\bf M}$ is a $Z$-free $E$-$E$-bimodule. The above inclusion induces an inclusion $\mathcal{F}_{E}(M)\rightarrow \mathcal{F}_{E}({\bf M})$. We have the cyclic derivation $\widetilde{\mathfrak{h}}:\mathcal{F}_{E}({\bf M})\rightarrow \mathrm{End}_{F}(\mathcal{F}_{E}({\bf M}))$. By Proposition \ref{prop6.12}, if $f$ and $g \in \mathcal{F}_{E}(M)$ then $\widetilde{\mathfrak{h}}(f)(g)$ lies in $\mathcal{F}_{E}(M)$. Therefore the restriction of $\widetilde{\mathfrak{h}}$ to $\mathcal{F}_{E}(M)$ gives the desired cyclic derivation.
\end{proof}

\begin{prop} \label{prop6.14} Let $M$ be any  $E$-$E$-bimodule and $P \in \mathcal{F}_{E}(M)_{cyc}$. Then the $F$-linear map
$$\widetilde{\delta }^{P}:M^{*}\rightarrow \mathcal{F}_{E}(M)$$
given by $\widetilde{\delta }^{P}(\psi )=\psi (\widetilde{\delta }(P))$,
is a morphism of $E$-$E$-bimodules.
\end{prop}

\begin{proof} For $t\in L_{i}$,
$$\widetilde{\delta }^{P}(t\psi )=\sum _{s\in L}t\psi \left(s^{-1}\delta (P)s\right)=t\left(\sum _{s\in L}\psi (s^{-1}\delta (P)s)\right)=t\widetilde{\delta }^{P}(\psi ),$$
and
$$\widetilde{\delta }^{P}(\psi t)=\psi \left(t\sum _{s\in L_{i}}s^{-1}e_{i}\delta (P)e_{i}s \right)=\psi \left(\sum _{s\in L_{i}}s^{-1}e_{i}\delta (P)e_{i}s\right)t=\widetilde{\delta }^{P}(\psi )t.$$
\end{proof}

\begin{definition} \label{def6.15} Let $M$ be any $E$-$E$-bimodule and $P \in \mathcal{F}_{E}(M)_{cyc}$. We define the Jacobian ideal $J(P)$ as the closure of the two-sided ideal in $\mathcal{F}_{E}(M)$ generated by $\widetilde{\delta }^{P}(M^{*})$.
\end{definition}

\begin{prop} \label{prop6.16} Suppose $M$ is an $E$-$E$-bimodule $Z$-freely generated by the set $A$, consisting of directed elements of $M$ and $P \in \mathcal{F}_{E}(M)_{cyc}$. Then  $R(P)=J(P)$.
\end{prop}

\begin{proof} First remark that for $a\in A$ and $s\in L(\tau (a))$, $(sa)^{*}=a^{*}s^{-1}$. Then
$$X_{a^{*}}(P)=\sum _{s\in L(\tau (a))}(sa)^{*}(\delta (P)s)=a^{*}\left(\sum _{s\in L(\tau (a))}s^{-1}(\delta (P)s)\right)$$
$$=a^{*}\left(\sum _{s\in L}s^{-1}\delta (P)s\right)=\widetilde{\delta }^{P}(a^{*}).$$

Moreover, if $\psi \in M^{*}$, then $\psi =\displaystyle \sum _{u,a}\lambda _{u}a^{*}\nu _{u}$ with
$\lambda _{u}\in E(\sigma (a)), \nu _{u}\in E(\tau (a))$. Therefore
$$\widetilde{\delta }^{P}(\psi )=\sum _{u,a}\lambda _{u}\widetilde{\delta }^{P}(a^{*})\nu _{u}=\sum _{u,a}\lambda _{u}X_{a^{*}}(P)\nu _{u}.$$ 
and the statement follows.
\end{proof}

\begin{prop} \label{prop6.17} Let $M$ be any $E$-$E$-bimodule and $P \in \mathcal{F}_{E}(M)_{cyc}$.  Suppose that $\varphi $ is an algebra automorphism of $\mathcal{F}_{E}(M)$ with $\varphi |_{E}=id_{E}$. Then $$\varphi (J(P))=J(\varphi (P)).$$
\end{prop}

\begin{proof} Take $M\rightarrow {\bf M}$ an inclusion of $E$-$E$-bimodules with ${\bf M}$ a $Z$-free $E$-$E$-bimodule. Then there is an $E$-$E$-bimodule $M^{\prime }$ with ${\bf M}=M\oplus M^{\prime}$. Let $\mathcal{L}$ be the closure of the two-sided ideal in $\mathcal{F}_{E}({\bf M })$ generated by $M^{\prime }$. We have
$$\mathcal{F}_{E}({\bf M})=\mathcal{L}\oplus \mathcal{F}_{E}(M).$$

Consider $J_{{\bf M}}(P)$, the Jacobian ideal of $P$ in $\mathcal{F}_{E}({\bf M})$. Then
$$ J_{{\bf M}}(P)=(\mathcal{L}J(P)+J(P)\mathcal{L}+\mathcal{L}J(P)\mathcal{L})\oplus J(P)$$
and therefore $J(P)=\mathcal{F}_{E}(M)\cap J_{{\bf M}}(P)$. Now let $\varphi$ be an algebra automorphism of $\mathcal{F}_{E}(M)$, with $\varphi |_{E}=id_{E}$. This automorphism
is determined by a pair of morphisms of $E$-$E$-bimodules $\varphi ^{(1)}:M\rightarrow M$ and
$\varphi ^{(2)}:M\rightarrow \mathcal{F}_{E}(M)^{\geq 2}$. Take now
$$\phi ^{1}=\left ( \begin{array}{cc} \varphi ^{(1)}&0\\0&id_{M^{\prime }} \end{array}\right ) :M\oplus M^{\prime }\rightarrow M\oplus M^{\prime }$$
and 
$$\phi ^{(2)}:M\oplus M^{\prime }\rightarrow \mathcal{F}_{E}({\bf M})$$
such that $\phi ^{(2)}|_{M}=\varphi ^{(2)}$ and $\phi ^{(2)}(M^{\prime })=0$.

The above pair of morphisms determine an algebra automorphism $\phi $ of $\mathcal{F}_{E}({\bf M})$ whose restriction
to $\mathcal{F}_{E}(M)$ is $\varphi $.

Then
$$\varphi (J(P))=\phi (J_{{\bf M}}(P)\cap \mathcal{F}_{E}(M))=\phi (J_{{\bf M}}(P))\cap \phi (\mathcal{F}_{E}(M))$$
$$=J_{{\bf M }}(\phi (P))\cap \mathcal{F}_{E}(M)=J_{{\bf M}}(\varphi (P))\cap \mathcal{F}_{E}(M)=J(\varphi (P)).$$

The proof is now complete.
\end{proof}

\begin{prop} \label{prop6.18} Suppose $M$ is any $E$-$E$-bimodule with polarization $\mathcal{A}$ and $P $ a potential in $M^{\otimes 2}$ such that
$\widetilde{\delta }^{P}:M^{*}\rightarrow M$ is an isomorphism of $E$-$E$-bimodules. Then $P$ is cyclically equivalent to a potential
$$Q=\sum _{x}\alpha _{x}\beta _{x}$$
where each $E\alpha _{x}E$ and $E\beta _{x}E$ are simple $E$-$E$-bimodules and
$$\bigoplus _{(i,j)\in \mathcal{A}}e_{i}Me_{j}=\bigoplus _{x}E\alpha _{x}E$$
$$\bigoplus _{(i,j)\in \mathcal{A}}e_{j}Me_{i}=\bigoplus _{x}E\beta _{x}E.$$
\end{prop}

\begin{proof} Since $\widetilde{\delta }^{P}$ is an isomorphism, it induces isomorphisms
$$\widetilde{\delta }^{P}:(e_{i}Me_{j})^{*}\rightarrow e_{j}Me_{i}.$$

We recall that we have a decomposition of $E_{i}\otimes _{F}E_{j}^{op}=\displaystyle \bigoplus _{z\in U(i,j)}E_{i}\epsilon _{z}E_{j}^{op}$ into a sum of simple $E_{i}\otimes _{F}E_{j}^{op}$- modules and
$E_{j}\otimes _{F}E_{i}^{op}=\displaystyle \bigoplus _{z\in U(i,j)}E_{j}\mathfrak{s}(\epsilon _{z})E_{i}^{op}$ is a sum of simple
$E_{j}\otimes _{F}E_{i}^{op}$-modules. Here $e_{i}\otimes e_{j}=\displaystyle \sum _{z\in U(i,j)}\epsilon _{z}$ is a sum of elements of a set of primitive orthogonal idempotents $\epsilon _{z}$.

Then 
$$\displaystyle \bigoplus _{(i,j)\in \mathcal{A}}e_{i}Me_{j}=\displaystyle \bigoplus _{(i,j)\in \mathcal{A}}\displaystyle \bigoplus _{x\in V(i,j)}E\alpha _{x}E$$
with each $E\alpha _{x}E$ a simple $E$-$E$-bimodule and $\alpha _{x}=\epsilon _{x}*\alpha _{x}$, where
$\epsilon _{x}$ is an idempotent in $U(i,j)$. Then if $V=\displaystyle \bigcup _{(i,j)\in \mathcal{A}}V_{i,j}, U=\displaystyle \bigcup _{(i,j)\in \mathcal{A}}U_{i,j}$ we have
$$\displaystyle \bigoplus _{(i,j)\in \mathcal{A}}e_{i}Me_{j}=\displaystyle \bigoplus _{x\in V}E\alpha _{x}E.$$

We have $\alpha _{x}=\epsilon _{x}*\alpha _{x}$, for some idempotent $\epsilon _{x}\in U$. The potential $P$ is cyclically equivalent to a potential
$$P_{1}=\sum _{x\in V}\alpha _{x}y_{x}=\sum _{x\in V}\epsilon _{x}*\alpha _{x}y_{x}$$
which is cyclically equivalent to
$$P_{2}=\sum _{x\in V}\alpha _{x}\beta _{x}$$
with $\beta _{x}=\mathfrak{s}(\epsilon _{x})*y_{x}$. Therefore $E\beta _{x}E$ is a simple $E$-$E$-bimodule.

Now
$$\widetilde{\delta }^{P}\left(\displaystyle \bigoplus _{x\in V}(E\alpha _{x}E)^{*}\right)=\displaystyle \bigoplus _{x\in V}\widetilde{\delta }^{P}((E\alpha _{x}E)^{*})$$

Take $\psi \in (E\alpha _{x}E)^{*}$, then
$$\widetilde{\delta }^{P}(\psi )=\psi \widetilde{\delta (P)}=\psi \widetilde{\delta }(P_{2})=\sum _{s\in L}\psi (s^{-1}\alpha _{x})\beta _{x}s.$$

Therefore $\widetilde{\delta }^{P}((E\alpha _{x}E)^{*})=E\beta _{x}E$. From here we obtain
$$\displaystyle \bigoplus _{x\in V}E\beta _{x}E=\widetilde{\delta }^{P}\left(\displaystyle \bigoplus _{(i,j)\in \mathcal{A}}e_{i}Me_{j}\right)^{*}=\displaystyle \bigoplus _{(i,j)\in \mathcal{A}}e_{j}Me_{i}.$$

Then:
$$\displaystyle \bigoplus _{(i,j)\in \mathcal{A}}e_{j}Me_{i}=\displaystyle \bigoplus _{x\in V}E\beta _{x}E.$$

The proof is now complete.
\end{proof}

\begin{prop} \label{prop6.19} Suppose $M$ is an $E$-$E$-bimodule $Z$-freely generated by the set $A$, consisting of directed elements. Let $\mathcal{A}$ be a polarization of $M$ and let $A(\mathcal{A})=\displaystyle \bigcup _{(i,j)\in \mathcal{A}}A(i,j)$, where $A(i,j)=A \cap e_{i}Me_{j}$.
Then:
\begin{enumerate}[(i)]
\item If $P$ is a quadratic potential in $\mathcal{F}_{E}(M)$, then $P$ is cyclically equivalent to
the potential:
$$\sum _{a\in A(\mathcal{A})}a\widetilde{\delta }^{P}(a^{*})$$
\item If $P$ is a quadratic potential in $\mathcal{F}_{E}(M)$ and
$\widetilde{\delta }^{P}:M^{*}\rightarrow M$ is an isomorphism of $E$-$E$-bimodules, then $M$ is a $2$-cycle complete $E$-$E$-bimodule
with respect to $A$; furthermore, there exists an algebra automorphism $\varphi $ of $\mathcal{F}_{E}(M)$ such that 
$\varphi (P)$ is cyclically equivalent to 
$$\sum _{a\in A(\mathcal{A})}a\hat{a}.$$
\item If $M$ is a $2$-cycle complete $E$-$E$-bimodule with respect to $A$  and $P=\displaystyle \sum _{a\in A(\mathcal{A})}a\hat{a}$, then
$$\widetilde{\delta }^{P}:M^{*}\rightarrow M$$
is an isomorphism of $E$-$E$-bimodules.
\end{enumerate}
\end{prop}

\begin{proof}
\begin{enumerate}[(i)]
\item The potential $P$ is cyclically equivalent to a potential
$$P^{\prime }=\sum _{a\in A(\mathcal{A})}ay_{a}$$ 

with $y_{a}\in e_{\sigma (a)}Me_{\tau (a)}$. Then, for a fixed $a_{0}\in A(\mathcal{A})$ we have
$$\widetilde{\delta }^{P}(a_{0}^{*})=\sum _{s\in L,a \in A(\mathcal{A})}a_{0}^{*}(s^{-1}(ay_{a}+y_{a}a)s)$$
$$=\sum _{s\in L)}(sa_{0})^{*}\displaystyle \left(\sum _{a\in A(\mathcal{A})}ay_{a}\right)=e_{\sigma  (a_{0})}y_{a_{0}}=y_{a_{0}}.$$
This proves (i).
\item Since $\widetilde{\delta }^{P}$ is an isomorphism, then this map induces an isomorphism of $E$-$E$-bimodules
$$\widetilde{\delta }^{P}:\displaystyle \bigoplus _{(i,j)\in \mathcal{A}}(e_{i}Me_{j})^{*}\rightarrow \displaystyle \bigoplus _{(i,j)\in \mathcal{A}}e_{j}Me_{i}.$$
The elements  $\widetilde{\delta }^{P}(a^{*})$ form a set of $Z$-free generators of
$\displaystyle \sum _{(i,j)\in \mathcal{A}}e_{j}Me_{i}$. Therefore there is an isomorphism of $E$-$E$-bimodules $\phi :M\rightarrow M$ such that $\phi (a)=a$ for $a\in A(\mathcal{A})$ and  $\phi (\widetilde{\delta }^{P}(a^{*}))=\hat{a}\in A$.
Let $\varphi $ be the algebra automorphism of $\mathcal{F}_{E}(M)$ induced by the pair $(\phi ,0)$. Then
$P$ is cyclically equivalent to $P^{\prime }=\displaystyle \sum _{a\in A(\mathcal{A})}a\widetilde{\delta }^{P}(a^{*})$.
Therefore $\varphi (P)$ is cyclically equivalent to
$$\varphi (P^{\prime })=\displaystyle \sum _{a\in A(\mathcal{A})}a\hat{a}.$$
This proves (ii).
\item For $a_{0}\in A(\mathcal{A})$ we have
$$\widetilde{\delta }^{P}(a_{0}^{*})=\sum _{s\in L, a \in A(\mathcal{A})}a_{0}^{*}s^{-1}(a\hat{a}+\hat{a}a)s$$
$$=\sum _{s\in L, a \in A(\mathcal{A})}(sa_{0})^{*}(a\hat{a})s=\hat{a}_{0}.$$
Similarly, $\widetilde{\delta }(P)((\hat{a_{0}})^{*})=a_{0}$. This completes the proof. 
\end{enumerate}
\end{proof}

\end{section}

\begin{section}{Premutations} \label{sec8}
In this section we assume that $E=\displaystyle \prod _{i=1}^{n}E_{i}$ is a product of semisimple finite dimensional $F$-algebras. We will not make any assumptions on the $F$-bases of $E$. Before providing a definition of mutation, we define the notion of premutation for $E$-$E$-bimodules satisfying the following

\begin{definition} \label{def7.1} Let $M$ be an $E$-$E$-bimodule such that  each $e_{i}Me_{j}\neq 0$ is an $E_{j}$-free right module and an $E_{i}$-free left module. Then a local basis for $_{E}M$ is a set $X=\displaystyle \bigcup _{i,j}X(i,j)$, with $X(i,j)=e_{i}Me_{j}\cap X$ a free $E_{i}$-basis for $e_{i}Me_{j}$ provided $e_{i}Me_{j} \neq 0$; otherwise, $X(i,j)=\emptyset$.

Similarly, a set $Y=\displaystyle \bigcup _{i,j}Y(i,j)$ is a local basis for $M_{E}$ if each $Y(i,j)$ is a free $E_{j}$-basis for $e_{i}Me_{j}$ provided $e_{i}Me_{j} \neq 0$; otherwise, $Y(i,j)=\emptyset$.
\end{definition}

\begin{definition} \label{def7.2} Let $M$ be an $E$-$E$-bimodule with a local basis $X$ for $_{E}M$ and a local basis $Y$ for
$M_{E}$. Let $k$ be an integer in $[1,n]$. Suppose that $e_{i}Me_{k}\neq 0$ implies $e_{k}Me_{i}=0,$ and $e_{k}Me_{i}\neq 0$ implies
$e_{i}Me_{k}=0$.  Following \cite{4}, we define the premutation of $M$, in the direction of $k$, as
$$\mu _{k}M=\overline{e}_{k}M\overline{e}_{k}\oplus Me_{k}M\oplus (e_{k}M)^{*}\oplus ^{*}(Me_{k}),$$  
where $\overline{e}_{k}=1-e_{k}$.
\end{definition}

We will use the following notation, for $m\in \overline{e}_{k}M\overline{e}_{k}$, we denote by $[m]$ the inclusion of $m$ into the first summand of $\mu _{k}M$, and for $f\in Me_{k}M$ we denote by $[f]$ the inclusion of $f$ into the second summand of $\mu _{k}M$. Observe that by our hypothesis, $e_{i}(\mu_{k}M)e_{k} \neq 0$ implies $e_{k}(\mu_{k}M)e_{i}=0$; likewise, $e_{k}(\mu_{k}M)e_{i} \neq 0$ implies $e_{i}(\mu_{k}M)e_{k}=0$. We will also make use of the following $E$-$E$-bimodule
$$\widehat{M}=M\oplus (e_{k}M)^{*}\oplus ^{*}(Me_{k})$$
We will identify $M$, $(e_{k}M)^{*}$ and $^{*}(Me_{k})$ with their respective inclusions into $\widehat{M}$ .

The inclusions $M\rightarrow \widehat{M}$ and $\mu _{k}M)\rightarrow \mathcal{F}_{E}(\widehat{M})$ induce inclusions of topological algebras
$$i_{M}:\mathcal{F}_{E}(M)\rightarrow \mathcal{F}_{E}(\widehat{M});$$
$$i_{\mu _{k}M}: \mathcal{F}_{E}(\mu _{k}M)\rightarrow \mathcal{F}_{E}(\widehat{M}).$$

\begin{prop} \label{prop7.3} There exists a morphism of topological algebras
$$[-]: \overline{e}_{k}\mathcal{F}_{E}(M)\overline{e}_{k}\rightarrow \mathcal{F}_{E}(\mu _{k}M)$$
which induces an isomorphism of topological algebras:
$$\overline{e}_{k}\mathcal{F}_{E}(M)\overline{e}_{k}\rightarrow \mathrm{Im}([-])=\mathcal{F}_{E_{\overline{k}}}(\overline{e}_{k}M\overline{e}_{k}\oplus Me_{k}M),$$
where $E_{\overline{k}}=\overline{e}_{k}E$. Moreover, for $h\in \overline{e}_{k}\mathcal{F}_{E}(M)\overline{e}_{k}$ we have $i_{M}(h)=i_{\mu _{k}M}([h])$.
\end{prop}

\begin{proof} Let $I$ be the collection of all finite sequences $\lambda =(i_{1},\ldots,i_{s})$ where $i_{1},\ldots,i_{s}\in [1,n]$. By $I_{k}$ we denote the subset of $I$ consisting of all sequences $\lambda ^{\prime}=(j_{1},\ldots,j_{s})$ where all the $j_{i}\neq k$. Given $\lambda \in I$, we define the $E$-$E$-bimodule $M_{\lambda }=E_{i_{1}}$ in case $\lambda =(i_{1})$; and for $\lambda =(i_{1},\ldots,i_{s})$, where $s\geq 2$,
$M_{\lambda }=e_{i_{1}}Me_{i_{2}}\otimes _{E}e_{i_{2}}Me_{i_{3}}\otimes _{E} \cdots \otimes _{E}e_{i_{s-1}}Me_{i_{s}}$.
We define $|\lambda |=s-1$, and $I(t)=\{\lambda \in I: |\lambda |=t\}$. Clearly
$$M^{\otimes l}=\displaystyle \bigoplus _{\lambda \in I(l)}M_{\lambda }.$$
Now let $T$ be the set that consists of all finite sequences $(\lambda _{1},\ldots,\lambda _{r})$ of elements $\lambda _{1},\ldots,\lambda _{r}\in I$ such that
$\lambda _{1}=(\lambda _{1}^{\prime },k), \lambda _{2}=(k,\lambda _{2}^{\prime }),\ldots,\lambda _{t}=(k,\lambda ^{\prime }_{t})$, where $\lambda ^{\prime }_{1},\ldots,\lambda _{t}^{\prime }\in I_{k}$.

Now consider the $E$-$E$-bimodules
$$U(\lambda _{1},\ldots,\lambda _{t})=M_{\lambda _{1}}\otimes _{E}M_{\lambda _{2}}\otimes _{E} \cdots \otimes _{E}M_{\lambda _{t}}=$$
$$M_{\lambda _{1}^{\prime }}\otimes _{E}Me_{k}\otimes _{E}e_{k}M\otimes _{E}M_{\lambda _{2}^{\prime }}\otimes _{E} \cdots \otimes _{E}Me_{k}\otimes _{E}e_{k}M\otimes _{E}M_{\lambda ^{\prime }_{t}},$$
and
$$V(\lambda _{1},\ldots,\lambda _{t})=M_{\lambda _{1}^{\prime }}\otimes _{E}Me_{k}M\otimes _{E}M_{\lambda ^{\prime}_{2}}\otimes _{E}Me_{k}M\otimes _{E} \cdots \otimes _{E}Me_{k}M\otimes _{E}M_{\lambda _{t}^{\prime }}$$.

Then
$$\overline{e}_{k}T_{E}(M)\overline{e}_{k}=E_{\overline{k}} \bigoplus \bigoplus_{(\lambda _{1},\ldots,\lambda _{t})\in T}U(\lambda _{1},\ldots,\lambda _{t}),$$
$$T_{E_{\overline{k}}}(\overline{e}_{k}M\overline{e}_{k}\oplus Me_{k}M)=E_{\overline{k}} \bigoplus \bigoplus _{(\lambda _{1},\ldots,\lambda _{t})\in T}V(\lambda _{1},\ldots,\lambda _{t}).$$

Now the inclusion of $\overline{e}_{k}M\overline{e}_{k}$ into $\overline{e}_{k}T_{E}(M)\overline{e}_{k}$ and the isomorphism of $Me_{k}M$ into $Me_{k}\otimes _{E}e_{k}M \subseteq \overline{e}_{k}T_{E}(M)\overline{e}_{k}$ induce
a morphism of $E$-$E$-bimodules $\overline{e}_{k}M\overline{e}_{k}\oplus Me_{k}M\rightarrow \overline{e}_{k}T_{E}(M)\overline{e}_{k}$, which in turn induces a morphism of $F$-algebras
$$\mathfrak{u}_{0}:T_{E_{\overline{k}}}(\overline{e}_{k}M\overline{e}_{k}\oplus Me_{k}M) \rightarrow \overline{e}_{k}T_{E}(M)\overline{e}_{k}.$$

The morphism $\mathfrak{u}_{0}$ induces an isomorphism of $E$-$E$-bimodules $V(\lambda _{1},\ldots,\lambda _{s}) \rightarrow U(\lambda _{1},\ldots,\lambda _{s})$ for all $(\lambda _{1},\ldots,\lambda _{s})\in T$; thus $\mathfrak{u}_{0}$ is an
isomorphism of $F$-algebras. Therefore, $\mathfrak{u}_{0}$ induces an isomorphism of topological algebras
$$\mathfrak{u}:\mathcal{F}_{E_{\overline{k}}}(\overline{e}_{k}M\overline{e}_{k}\oplus Me_{k}M)\rightarrow \overline{e}_{k}\mathcal{F}_{E}(M)\overline{e}_{k},$$
Define
$$[-]=i\mathfrak{u}^{-1}:\overline{e}_{k}\mathcal{F}_{E}(M)\overline{e}_{k}\rightarrow \mathcal{F}_{E}(\mu _{k}M),$$
where $i:\mathcal{F}_{E_{\overline{k}}}(\overline{e}_{k}M\overline{e}_{k}\oplus Me_{k}M)\rightarrow \mathcal{F}_{E}(\mu _{k}M)$ is the inclusion map. This proves the first part of our proposition. Now for $h\in \overline{e}_{k}M\overline{e}_{k}\oplus Me_{k}\otimes _{E}e_{k}M$,  we have $i_{M}(h)=i_{\mu_{kM}}([h])$. Since $ \overline{e}_{k}M\overline{e}_{k}\oplus Me_{k}\otimes _{E}e_{k}M $ generates $\overline{e}_{k}T_{E}(M)\overline{e}_{k}$, then $i_{M}(h)=i_{\mu _{k}M}([h])$ for all $h\in \overline{e}_{k}T_{E}(M)\overline{e}_{k}$. Now we are done since the result follows by continuity.

\end{proof}

\begin{prop} \label{prop7.5} If $e_{i}\mu _{k}Me_{j}\neq 0$, then it is an $E_{i}$-free left module and an $E_{j}$-free right module.
\end{prop}

\begin{proof} If $i$ and $j$ are not equal to $k$, we have
$$e_{i}\mu _{k}Me_{j}=e_{i}Me_{j}\oplus e_{i}Me_{k}Me_{j}$$
Here $e_{i}Me_{j}$ is free as a left $E_{i}$-module, and it is also a free right $E_{j}$-module. Now 
$$e_{i}Me_{k}Me_{j}=e_{i}Me_{k}\otimes _{E_{k}}e_{k}Me_{j}.$$
Since $e_{k}Me_{j}$ is a free left $E_{k}$-module, then $e_{i}Me_{k}Me_{j}\cong (e_{i}Me_{k})^{l}$, where $l=\mathrm{Card}X(k,j)$. Therefore $e_{i}Me_{k}Me_{j}$ is a free left $E_{i}$-module. Similarly, one can see that $e_{i}Me_{k}Me_{j}$ is a free $E_{j}$-right module. An $E_{i}$-free basis for $e_{i}Me_{k}Me_{j}$ is given by
$$\hat{X}_{k}(i,j)=\{x\otimes x_{1} \}_{x\in X(i,k), x_{1}\in X(k,j)}.$$
An $E_{j}$-free basis for $e_{i}Me_{k}Me_{j}$ is given by 
$$\hat{Y}_{k}(i,j)=\{y\otimes y_{1}\}_{y\in Y(i,k), y_{1}\in Y(k,j)}.$$

Now if $i\ne k$ and $j=k$, then $e_{i}\mu _{k}Me_{k}=(e_{k}Me_{i})^{*}$. Clearly $(e_{k}Me_{i})^{*}$ is a free left  $E_{i}$-module and $\hat{X}(i,k)=\{y^{*}\}_{y\in Y(k,i)}$ is an $E_{i}$-free basis. We recall that we can choose an isomorphism of $E$-$E$-bimodules $(e_{k}Me_{i})^{*}\rightarrow ^{*}(e_{k}Me_{i})$, therefore  $(e_{k}Me_{i})^{*}$ is also an $E_{k}$-free right module; and $\hat{Y}(i,k)$, the image of the set $\{x^{\ast}: x \in _{k}X\}$ under the above chosen isomorphism, is an $E_{k}$-free basis. Similarly, one can see that if $e_{i}Me_{k}\neq 0$, then $^{*}(e_{i}Me_{k})$ is free as a left $E_{k}$-module and as a right $E_{i}$-module.
A free basis of $^{*}(e_{i}Me_{k})$, as a right $E_{i}$-module, is given by the set $\hat{Y}(k,i)=\{ ^{*}x| x\in X(i,k)\}$ and a free basis, as a left $E_{k}$-module, is given by the set $\hat{X}(k,i)$ which consists of the images of the elements $y^{*}, y\in Y(i,k)$ under a chosen isomorphism $(e_{i}Me_{k})^{*}\rightarrow ^{*}(e_{i}Me_{k})$.  
\end{proof}

\begin{definition} \label{def7.6} Let $P$ be a potential in $\mathcal{F}_{E}(M)$ such that $e_{k}Pe_{k}=0$, then by Proposition \ref{prop7.3} there exists a potential $[P]$ in $\mathcal{F}_{E}(\mu _{k}M)$ such that $i_{\mu _{k}M}([P])=i_{M}(P)$. Following \cite{4}, we define
$$\mu _{k}(P)=[P]+\sum _{x\in X_{k},y\in  _{k}Y}[xy](y^{*})(^{*}x).$$
Here $[xy]\in Me_{k}M\subseteq \mu _{k}M$ and $X_{k}=\displaystyle \bigcup _{i}X(i,k)$, $_{k}Y=\displaystyle \bigcup _{i}Y(k,i)$.
\end{definition}

$ $ 
Observe that the expressions $\displaystyle \sum _{x\in X_{k},y\in  _{k}Y}[xy](y^{*})(^{*}x)$, $\displaystyle \sum _{x\in  X_{k}}(^{*}x)x$, and $\displaystyle \sum _{y\in _{k}Y}y(y)^{*}$ do not depend on the choice of the local bases $X$ and $Y$ for $_{E}M$ and $M_{E}$, respectively.
 
 Indeed, assume that $Y_{1}$ is another local basis for $M_{E}$, then
 $$\sum _{x\in X_{k},y\in  _{k}Y}[xy](y^{*})(^{*}x)=\sum _{x\in X_{k},y \in _{k}Y, y_{1}\in _{k}(Y_{1})}[xy_{1}]y_{1}(y_{1})^{*}(y)y^{*}(^{*}x)=$$
 $$\sum _{x\in X_{k}, y_{1}\in _{k}(Y_{1})}[xy_{1}]\left(\sum _{y\in _{k}Y}(y_{1})^{*}(y)y^{*})(^{*}x)\right)=\sum _{x\in X_{k},y_{1}\in _{k}(Y_{1})}[xy_{1}](y_{1}^{*})(^{*}x).$$
 Similarly, if $X_{1}$ is another local basis for $_{E}M$ we have
 $$\sum _{x\in X_{k},y\in  _{k}Y}[xy](y^{*})(^{*}x)=\sum _{x\in X_{k},y_{1}\in _{k}(Y_{1})}[xy_{1}](y_{1}^{*})(^{*}x)$$
 $$=\sum _{x_{1}\in (X_{1})_{k},y_{1}\in _{k}(Y_{1})}[x_{1}y_{1}](y_{1}^{*})(^{*}x_{1}).$$
 The invariance of the other expressions is proved in the same way.

\begin{definition} \label{def7.7} An algebra automorphism $\varphi $ of $\mathcal{F}_{E}(M)$ given by a pair of $E$-$E$-bimodule morphisms $(id_{M},\varphi ^{(2)})$ is called unitriangular.
\end{definition}

\begin{prop} \label{prop7.8} Let $\varphi: \mathcal{F}_{E}(M) \rightarrow \mathcal{F}_{E}(M_{1})$ be an algebra isomorphism determined by a pair of $E$-$E$-bimodule morphisms $(\varphi ^{(1)},\varphi ^{(2)})$, where $M$ and $M_{1}$ are $E$-$E$-bimodules satisfying the conditions of Definition $7.1$. Moreover, suppose $X$ and $X^{1}$ are local bases for $_{E}M$ and $_{E}M_{1}$, respectively; and $Y$, $Y^{1}$ are local bases for $M_{E}$ and $(M_{1})_{E}$, respectively. Then there exist algebra isomorphisms  $\hat{\phi }:\mathcal{F}_{E}(\widehat{M})\rightarrow \mathcal{F}_{E}(\widehat{M}_{1})$ and $\phi: \mathcal{F}_{E}(\mu _{k}M)\rightarrow \mathcal{F}_{E}(\mu_{k}M_{1})$, such that
$$\hat{\phi }i_{\mu _{k}M}=i_{\mu _{k}M_{1}}\phi $$
$$\hat{\phi }i_{M}=i_{M_{1}}\varphi .$$

Moreover,
$$\hat{\phi }\left(\sum _{x\in X_{k}}(^{*}x)x\right)=\sum _{x_{1}\in X^{1}_{k}}(^{*}\phi ^{(1)}(x))\varphi ^{(1)}(x)$$
$$\hat {\phi }\left(\sum _{y\in _{k}Y}yy^{*}\right)=\sum _{y_{1}\in _{k}Y^{1}}y_{1}y_{1}^{*}.$$
\end{prop}

\begin{proof} First we define a right action of $^{*}M$ on $\mathcal{F}_{E}(M)$ as follows: take $\nu \in ^{*}M$, then for $s\in E$ we define $s\nu =0$. For $m_{1} \cdots m_{l}$, with $m_{i}\in M$, we define $(m_{1} \cdots m_{l})\nu =m_{1} \cdots m_{l-1}\nu(m_{l})$; for $f=\displaystyle \sum _{u=0}^{\infty}f(u)$, we define
$$f\nu =\sum _{u=0}^{\infty }f(u)\nu .$$

We now consider the particular case in which $M=M_{1}$ and $\varphi $ is a unitriangular automorphism.

Take  the restriction of $\varphi $ to $Me_{k}$

$$\varphi :Me_{k}\rightarrow \mathcal{F}_{E}(M)e_{k}$$

Let $m(k)=\mathrm{Card}X_{k}$ and consider the $m(k)\times m(k)$-matrix with entries in $\mathcal{F}_{E}(M)$
$$H=(h_{x,x^{\prime }})_{x,x^{\prime }\in X_{k}}$$
where $h_{x,x^{\prime }}=\varphi (x)(^{*}x^{\prime})$. Observe that $h_{x,x^{\prime }}=e_{\tau (x)}h_{x,x^{\prime }}e_{\tau (x^{\prime })}$. We have
$$\varphi (x)=\sum _{x^{\prime }}h_{x,x^{\prime }}x^{\prime }.$$

For $s\in E$ and $x\in X_{k}$, we have
$$xs=\sum _{x^{\prime }\in X_{k}}w(x ,x^{\prime })(s)x^{\prime }$$
with $w(x,x^{\prime })(s)\in E_{\tau (x^{\prime })}$. Then
$$\varphi (xs)=\sum _{x^{\prime }, x^{\prime \prime }\in X_{k}}w(x,x^{\prime })(s)h_{x^{\prime },x^{\prime \prime }}x^{\prime \prime }; \ \varphi (x)s=\sum _{x^{\prime },x^{\prime \prime }}h_{x,x^{\prime }}w(x^{\prime },x^{\prime \prime })(s)x^{\prime \prime }.$$
Therefore for any $x_{0}\in X_{k}$, we obtain:
$$\varphi (xs)^{*}x_{0}=\varphi (x)s^{*}x_{0},$$

Taking into account that $h_{x^{\prime },x_{0}}e_{\tau (x_{0})}=h_{x, x^{\prime }}$ and that
$w(x^{\prime },x_{0})e_{\tau (x_{0})}=w(x^{\prime },x_{0})$, yields

$$\sum _{x^{\prime }\in X_{k}}w(x,x^{\prime })(s)h_{x^{\prime },x_{0}}=\sum _{x^{\prime }\in X_{k}}h_{x,x^{\prime }}w(x^{\prime },x_{0})(s).$$

Therefore if we consider the $m(k)\times m(k)$-matrix with entries in $E$, $W(s)=(w(x,x^{\prime })(s))$, we have
$$W(s)H=HW(s)$$
for all $s\in E_{k}$.  Now consider the restriction of $\varphi$ to $e_{k}M$
$$\varphi :e_{k}M\rightarrow e_{k}\mathcal{F}_{E}(M).$$

Let $n(k)=\mathrm{Card}_{k}Y$. We define a matrix $G=(g_{y^{\prime },y})$, of order $n(k)$ and with entries in $\mathcal{F}_{E}(M)$, as $g_{y^{\prime },y}=(y^{\prime })^{*}\varphi (y)$; here $e_{\sigma (y^{\prime })}g_{y^{\prime },y}e_{\sigma (y)}=g_{y^{\prime },y}$.

We have
$$\varphi (y)=\sum _{y^{\prime }\in _{k}Y}y^{\prime }g_{y^{\prime },y}.$$

For $s\in E_{k}$ we have 
$$sy=\sum _{y^{\prime }\in _{k}Y}y^{\prime }\omega (y^{\prime },y)(s)$$
with $\omega (y^{\prime },y)(s)\in E_{\sigma (y^{\prime })}$. Therefore if we take the $n(k)\times n(k)$-matrix
$\Omega (s) =(\omega (y^{\prime },y)(s))$, we have
$$G\Omega (s)=\Omega (s)G$$
for all $s\in E$.

Now we are going to define a morphism $\phi :\widehat{M}\rightarrow \mathcal{F}_{E}(\widehat{M})$. We first define a morphism of $E$-$E$-bimodules
$$\rho _{1}: ^{*}(Me_{k})
\rightarrow e_{k}\mathcal{F}_{E}(\widetilde{M}). $$
First consider the matrix $H$ associated to the morphism
$$\varphi :Me_{k}\rightarrow \mathcal{F}_{E}(M).$$

The matrix $H$ lies in the set of matrices $\mathcal{M}$ of size $m(k)\times m(k)$, with entries in $\mathcal{F}_{E}(M)$,
consisting of all the matrices $U=(u_{x,x^{\prime}})$ such that $e_{\tau (x)}u_{x,x^{\prime }}e_{\tau (x^{\prime })}=u_{x,x^{\prime }}$. The set $\mathcal{M}$ is an $F$-vector subspace of the $F$-algebra of all $n(k)\times n(k)$-matrices, with entries in $\mathcal{F}_{E}(M)$, and with unit $I=(\delta _{x,x^{\prime }}e_{\tau (x^{\prime })})$.

Since $\varphi $ is unitriangular, then $H=I+R$ where $R$ is a matrix with entries in $\mathcal{F}_{E}(M)^{\geq 1}$. Therefore the sequence $\left(\displaystyle \sum _{i=0}^{n}(-1)^{i}R^{i}\right)_{n \geq 1}$ converges to $\displaystyle \sum _{i=0}^{\infty }(-1)^{i}R^{i}$, where $R^{0}=I$; hence $H^{-1}=\displaystyle \sum _{i=0}^{\infty} (-1)^{i}R^{i}$.

Take $H^{-1}=(\hat{h}_{x,x^{\prime }})$ and define the $E$-left morphism $\rho _{1}$ such that
for $^{*}x$, with $x\in X_{k}$, we have
$$\rho _{1}(^{*}x)=\sum _{x^{\prime }\in X_{k}}(^{*}x^{\prime })\hat{h}_{x^{\prime },x}.$$
We claim that $\rho _{1}$ is also a right $E$-morphism. To show this we first calculate $s(^{*}x)$ for
$s\in E_{k}$.

We have 
$$s(^{*}x)=\sum _{x^{\prime }}(^{*}x^{\prime })s(^{*}x)(x^{\prime })=\sum _{x^{\prime }}(^{*}x^{\prime })(^{*}x)(x^{\prime }s)=$$
$$=\sum _{x^{\prime }}(^{*}x^{\prime })(^{*}x)\left(\sum _{x^{\prime \prime }}w(x^{\prime },x^{\prime \prime })(s)x^{\prime \prime }\right)=\sum _{x^{\prime }}(^{*}x^{\prime })w(x^{\prime },x)(s).$$
Then
$$\rho _{1}(s(^{*}x))=\sum _{x^{\prime \prime },x^{\prime }\in X_{k}}(^{*}x^{\prime \prime })\hat{h}_{x^{\prime \prime },x^{\prime }}w(x^{\prime },x)(s),$$
$$s\rho _{1}(^{*}x)=\sum _{x^{\prime \prime },x^{\prime }\in X_{k}}(^{*}x^{\prime \prime })w(x^{\prime \prime },x^{\prime })(s)\hat{h}_{x^{\prime },x}.$$

Now, from $W(s)H=HW(s)$ we obtain that $H^{-1}W(s)=W(s)H^{-1}$, and from this we have 
$$\rho _{1}(s(^{*}x))=s\rho _{1}(^{*}x)$$
for all $s\in E_{k}$ and $x\in X_{k}$. Therefore $\rho _{1}$ is a morphism of $E$-$E$-bimodules.

Now we define a morphism
$$\rho _{2}:(e_{k}M)^{*}\rightarrow \mathcal{F}_{E}(\widehat{M})e_{k}.$$
For this, we consider the matrix $G$ associated to the morphism
$$\varphi : e_{k}M\rightarrow e_{k}\mathcal{F}_{E}(M).$$

The matrix $G$ is in $\mathcal{M}^{\prime }$, the  $F$-vector subspace of the algebra of $n(k)\times n(k)$-matrices with entries in $\mathcal{F}_{E}(M)$ consisting of the matrices $(u_{y^{\prime },y})$ such that $e_{\sigma (y^{\prime })}u_{y^{\prime },y}e_{y}=u_{y^{\prime },y}$. Here $\mathcal{M}^{\prime }$ is an $F$-algebra with usual multiplication of matrices and with unit $I^{\prime}=(\delta _{y^{\prime },y}e_{\sigma (y)})$.  Since $\varphi $ is unitriangular it follows that $G$ has an inverse in $\mathcal{M}^{\prime }$.

Take $G^{-1}=(\tilde{g}_{y^{\prime },y})$, then for $y\in _{k}Y$ we define $\rho _{2}$ as the left $E$-morphism such that
$$\rho _{2}(y^{*})=\sum _{y^{\prime }\in _{k}Y}\tilde{g}_{y,y^{\prime }}(y^{\prime })^{*}.$$

Let us prove that $\rho _{2}$ is also a right $E$-morphism. Let $s\in E_{k}$, then
$$y^{*}s=\sum _{y^{\prime }\in _{k}Y}(y^{*}s)(y^{\prime })(y^{\prime })^{*}=\sum _{y^{\prime }\in _{k}Y}y^{*}(sy^{\prime })(y^{\prime })^{*}$$
$$=\sum _{y^{\prime },y^{\prime \prime }}y^{*}(y^{\prime \prime }\omega (y^{\prime \prime },y^{\prime })(s))(y^{\prime })^{*}=\sum _{y^{\prime }\in _{k}Y}\omega (y,y^{\prime })(s)(y^{\prime})^{*}.$$

Therefore
$$\rho _{2}(y^{*}s)=\sum _{y^{\prime },y^{\prime \prime }\in _{k}Y}\omega (y,y^{\prime })(s)\tilde{g}_{y^{\prime },y^{\prime \prime }}(y^{\prime \prime })^{*},$$
$$\rho _{2}(y^{*})s=\sum _{y^{\prime }, y^{\prime \prime }\in _{k}Y}\tilde{g}_{y,y^{\prime }}\omega (y^{\prime },y^{\prime \prime })(y^{\prime \prime })^{*}.$$

From the equality $\Omega (s)G=G\Omega (s)$ we obtain $G^{-1}\Omega (s)=\Omega (s)G^{-1}$, and from this last
equality we deduce that $\rho _{2}(y^{*}s)=\rho _{2}(y^{*})s$ for all $y\in _{k}Y$ and $s\in E_{k}$. Therefore $\rho _{2}$ is
a morphism of $E$-$E$-bimodules. Consider now the morphism of $E$-$E$-bimodules
$$\hat{\phi }=(\varphi , \rho _{1},\rho _{2}):\widehat{M}\rightarrow \mathcal{F}_{E}(\widehat{M})$$
this map can be extended to an algebra automorphism of $\mathcal{F}_{E}(\widehat{M})$, which we also denote by $\hat{\phi}$. Observe that for $w\in \widehat{M}$, $\hat{\phi }(w)=w+w^{\prime }$ with $w^{\prime }\in \mathcal{F}_{E}(\widehat{M})^{\geq 2 }$, hence the automorphism $\hat{\phi }$ is unitriangular.
 
Now we obtain
$$\hat{\phi }\left(\sum _{x\in X_{k}}(^{*}x)x\right)=\sum _{x\in X_{k}}\hat{\phi }(^{*}x)\hat{\phi }(x)=\sum _{x\in X_{k}}\rho _{1}(^{*}x)\varphi (x)=$$
$$\sum _{x^{\prime }, x^{\prime \prime }\in X_{k}}(^{*}x^{\prime })\hat{h}_{x^{\prime },x}h_{x.x^{\prime \prime }}x^{\prime \prime }=\sum _{x^{\prime },x^{\prime \prime }\in X_{k}}(^{*}x^{\prime })\delta _{x^{\prime },x^{\prime \prime }}e_{\tau (x^{\prime \prime })}x^{\prime \prime }=\sum _{x^{\prime }\in X_{k}}(^{*}x^{\prime })x^{\prime }.$$

Moreover,
$$\hat{\phi }\left(\sum _{y \in _{k}Y }yy^{*}\right)=\sum _{y\in _{k}Y}\hat{\phi }(y)\hat{\phi }(y^{*})=\sum _{y\in _{k}Y}\varphi (y)\rho _{2}(y^{*})$$
$$=\sum _{y^{\prime },y^{\prime \prime }\in _{k}Y}y^{\prime }g_{y^{\prime },y}\hat{g}_{y,y^{\prime \prime }}(y^{\prime \prime })^{*}=\sum _{y^{\prime },y^{\prime \prime }\in _{k}Y}y^{\prime }\delta _{y^{\prime },y^{\prime \prime }}e_{\sigma (y^{\prime \prime })}(y^{\prime \prime })^{*}$$
$$=\sum _{y^{\prime }\in _{k}Y}y^{\prime }(y^{\prime })^{*}.$$

Now observe that $(\hat{\phi }i_{M})|_{M}=(i_{M}\varphi )|_{M}$, therefore
$$\hat{\phi }i_{M}=i_{M}\varphi .$$

Moreover,
$$\hat{\phi }(\overline{e}_{k}M\overline{e}_{k})=\overline{e}_{k}\hat{\phi }(\overline{e}_{k}M\overline{e}_{k})\overline{e}_{k}.$$
and thus $\hat{\phi }(\overline{e}_{k}M\overline{e}_{k})\subseteq \mathrm{Im}i_{\mu _{k}M}$; similarly, $\hat{\phi }(Me_{k}M)\subseteq \mathrm{Im}i_{\mu _{k}M}$.

We have that all the elements $^{*}x$ are in $\mathrm{Im}i_{\mu _{k}M}$ and  $\tau (x)\neq k$ for all $x\in X_{k}$; then, by Proposition \ref{prop7.3}, the elements $\hat{h}_{x^{\prime },x}\in e_{\tau(x^{\prime })}\mathcal{F}_{E}(M)e_{\tau (x)}$ are in $\mathrm{Im}i_{\mu _{k}M}$; therefore $\rho _{1}^{*}(Me_{k})\subseteq \mathrm{Im}i_{\mu _{k}M}$. Similarly, $\rho _{2}((e_{k}M)^{*})$ is also contained in $\mathrm{Im}i_{\mu _{k}M}$.  From here we deduce that there is a morphism:
$\phi :\mu _{k}M\rightarrow \mathcal{F}_{E}(\mu _{k}M)$ such that $\hat{\phi }(i_{\mu _{k}M})|_{M}=i_{\mu _{k}M}\phi $. This map can be extended to an algebra endomorphism of $\mathcal{F}_{E}(\mu _{k}M)$ also denoted by $\phi $; thus
$$\hat{\phi }i_{\mu _{k}M}=i_{\mu _{k}M}\phi .$$

Finally, let us prove that $\phi $ is in fact an isomorphism. Let $m\in \overline{e}_{k}M\overline{e}_{k}\subseteq \mu _{k}M$, then
$$\hat{\phi }(i_{\mu _{k}M}(m))=\varphi (i_{\mu _{k}M}(m))=i_{\mu _{k}M}(m)+\sum _{i}m_{i}e_{k}m_{i}^{\prime }+
\sum _{j,u\ne k}n_{j}e_{u}n^{\prime }_{j}+\nu $$
where $m_{i},m^{\prime }_{i}\in M$, $n_{j},n^{\prime }_{j}\in \overline{e}_{k}M\overline{e}_{k}$, $\nu \in \mathcal{F}_{E}(M)^{\geq 3}$. Then $n_{j},n^{\prime }_{j}$ are in $i_{\mu _{k}M}(\mu _{k}M)$ and $\nu =i_{\mu _{k}}(\nu ^{\prime })$ for some $\nu ^{\prime }\in \mathcal{F}_{E}(\mu _{k}M)^{\geq 2}$. Therefore
$$i_{\mu _{k}M}(\hat{\phi }(m))=\hat{\phi }(i_{\mu _{k}M}(m))=i_{\mu _{k}M}\left(m+\sum _{i}[m_{i}e_{k}m^{\prime }]+r_{0}\right)$$
where $r_{0}\in \mathcal{F}_{E}(\mu _{k}M)^{\geq 2}$. From here we obtain
$$\phi (m)=m+\sum _{i}[m_{i}e_{k}m^{\prime }_{i}]+r_{0}.$$

Now for $f\in Me_{k}M$, take the corresponding element $[f]\in \mu _{k}M$; so $i_{\mu _{k}M}([f])=f$.
Then
$$ i_{\mu _{k}M}\phi ([f])=\hat{\phi }i_{\mu _{k}M}([f])=\hat{\phi }(f)=\varphi (f)=f+\nu _{1}$$
with $\nu _{1}\in \overline{e}_{k}\mathcal{F}_{E}(M)^{\geq 3}$. Therefore $\nu _{1}=i_{\mu _{k}M}(r_{1})$,
with $r_{1}\in \mathcal{F}_{E}(\mu _{k}M)^{\geq 2}$.  Then
$$i_{\mu _{k}M}\phi ([f])=i_{\mu _{k}M}([f]+r_{1})$$
hence
$$\phi ([f])=[f]+r_{1}.$$

For $^{*}x\in ^{*}(Me_{k})\subseteq \mu _{k}M$, we have
$$i_{\mu _{k}M}\phi (^{*}x)=\hat{\phi }i_{\mu _{k}M}(^{*}x)=^{*}x+\sum _{x^{\prime }\in X_{k}}(^{*}x^{\prime })u_{x^{\prime },x},$$
where $(u_{x^{\prime},x})=H-I$. Here $u_{x^{\prime},x}\in e_{\tau (x^{\prime })}\mathcal{F}_{E}(M)^{\geq 1}e_{\tau (x)}$, then $u_{x^{\prime },x}=i_{\mu _{k}M}([u_{x^{\prime },x}])$ with $[u_{x^{\prime },x}]\in \mathcal{F}_{E}(\mu _{k}M)^{\geq 1}$.
Therefore
$$i_{\mu _{k}M}(\phi (^{*}x))=i_{\mu _{k}M}(x^{*}+r_{3})$$
with $r_{3}\in \mathcal{F}_{E}(\mu _{k}M)^{\geq 2}$. From here we obtain
$$\phi (^{*}x)=^{*}x+r_{3}.$$
Similarly, we have 
$$\phi (y^{*})=y^{*}+r_{4}$$
with $r_{4}\in \mathcal{F}_{E}(\mu _{k}M)^{\geq 2}$. From the above we have that $\phi $ is determined by 
a pair of $E$-$E$-bimodule morphisms $(\phi ^{(1)},\phi ^{(2)})$, where 
$$\phi ^{(1)}=\left( \begin{array}{cccc}id&0&0&0\\
h&id&0&0\\
0&0&id&0\\
0&0&0&id \end{array}\right):\mu _{k}M\rightarrow \mu _{k}M$$
with $h:\overline{e}_{k}M\overline{e}_{k}\rightarrow Me_{k}M$. Then $\phi ^{(1)}$ is an isomorphism of $E$-$E$-bimodules which implies that
$\phi$ is an algebra automorphism. This proves our proposition for the case $M=M_{1}$, and when $\varphi$ is unitriangular.

We now prove the general case. First observe that if our proposition holds for the isomorphisms
$$\varphi _{1}:\mathcal{F}_{E}(M)\rightarrow \mathcal{F}_{E}(M_{2}), \varphi _{2}:\mathcal{F}_{E}(M_{2})\rightarrow \mathcal{F}_{E}(M_{3})$$
then the statements of our proposition also hold for the composition $\varphi _{2}\varphi _{1}$.

Now take the algebra isomorphism $\varphi _{1}:\mathcal{F}_{E}(M)\rightarrow \mathcal{F}_{E}(M_{1})$ determined by the pair
$(\varphi ^{(1)},0)$, then the morphism $\varphi _{2}=\varphi _{1}^{-1}\varphi $ is unitriangular and
$\varphi =\varphi _{1}\varphi _{2}$. Therefore, it suffices to prove the result for morphisms
$\varphi :\mathcal{F}_{E}(M)\rightarrow \mathcal{F}_{E}(M_{1})$ determined by a pair $(\varphi ^{(1)},0)$ where $\varphi ^{(1)}=\varphi |_{M}$.

Take the following morphisms of $E$-$E$-bimodules
$$\hat{\phi }^{1}=\left ( \begin{array}{ccc}\varphi ^{(1)}&0&0\\
0&((\varphi ^{-1})^{*})_{1}&0\\
0&0&(^{*}(\varphi ^{-1}))_{1}\end{array}\right):\widehat{M}\rightarrow \widehat{M}_{1},$$

$$\phi ^{(1)}=\left( \begin{array}{ccccc}\varphi _{1} &0&0&0\\
0&\varphi _{2}&0&0\\
0&0&((\varphi ^{-1})^{*}) _{1}&0\\
0&0&0&(^{*}(\varphi ^{-1}))_{1}  \end{array}\right) :\mu _{k}M\rightarrow \mu _{k}M_{1},$$
where $\varphi _{1}=\varphi |_{\overline{e}_{k}M\overline{e}_{k}}$; $\varphi _{2}=\varphi |_{Me_{k}M}$; $(\varphi ^{-1}) ^{*}:M^{*}\rightarrow M_{1}^{*}$ is the right dual map of $\varphi ^{-1}:M_{1}\rightarrow M$; $^{*}(\varphi ^{-1}):^{*}M\rightarrow  ^{*}M_{1}$ is the left dual map of $\varphi ^{-1}$; $((\varphi ^{-1})^{*})_{1}$ is the restriction of $(\varphi ^{-1})^{*}$ to $(e_{k}M)^{*}$;  and
$(^{*}\varphi ^{-1})_{1}$ is the restriction of $^{*}(\varphi ^{-1})$ to $^{*}(Me_{k})$. The pairs $(\hat{\phi }^{(1)},0)$ and $(\phi ^{(1)},0)$ define, respectively, algebra isomorphisms 
$$\hat{\phi }:\mathcal{F}_{E}(\widehat{M})\rightarrow \mathcal{F}_{E}(\widehat{M}_{1}), \phi :\mathcal{F}_{E}(\mu _{k}M)\rightarrow \mathcal{F}_{E}(\mu _{k}M_{1}).$$

One can see that $\hat{\phi }^{(1)}i_{M}|_{M}=i_{M_{1}}|_{M_{1}}\varphi ^{(1)}$,
$\hat{\phi }^{(1)}i_{\mu _{k}M}|_{\mu _{k}M}=i_{\mu _{k}M_{1}}|_{\mu _{k}M_{1}}\phi ^{(1)}$. From this, it follows that
$$\hat{\phi }i_{M}=i_{M_{1}}\varphi $$
and 
$$\hat{\phi }i_{\mu _{k}M}=i_{\mu _{k}M_{1}}\phi  .$$

Observe now that $\varphi (X)$ is a local basis for $_{E}M_{1}$ and $\varphi (Y)$ is a local basis for $(M_{1})_{E}$.
In particular, $\varphi (X_{k})$ is a local basis for $M_{1}e_{k}$ and $\varphi (_{k}Y)$ is a local basis for
$e_{k}M_{1}$. Now if $y^{\prime },y\in Y$, then 
$$(\varphi ^{-1})^{*}(y^{*})(\varphi (y^{\prime }))=y^{*}(\varphi ^{-1}(\varphi (y)))=y^{*}(y^{\prime })=\delta _{y^{\prime },y}e_{\sigma (y)}.$$
This proves that the set $\{\varphi (y), (\varphi ^{-1})^{*}(y^{*})|y \in _{k}Y\}$ is a dual basis for $e_{k}M_{1}$. Therefore 
$$\hat{\phi }\left(\sum _{y \in _{k}Y}yy^{*}\right)=\sum _{y\in _{k}Y}\hat{\phi }(y)\hat{\phi }(y^{*})$$
$$=\sum _{y\in _{k}Y}\varphi (y)(\varphi ^{-1})^{*} (y^{*})=\sum _{y_{1}\in _{k}Y_{1}}y_{1}y_{1}^{*}.$$
Similarly, one has that
$$\hat{\phi }\left(\sum _{x\in X_{k}}(^{*}x)x\right)=\sum _{x_{1}\in (X_{1})_{k}}(^{*}x_{1})x_{1}.$$
This completes the proof.
\end{proof}
$  $

\begin{theorem} \label{theo7.9} Let $\varphi :\mathcal{F}_{E}(M)\rightarrow \mathcal{F}_{E}(M_{1})$ be an algebra isomorphism
such that $\varphi |_{E}=id_{E}$. Then there is an algebra isomorphism $\phi :\mathcal{F}_{E}(\mu _{k}M)\rightarrow \mathcal{F}_{E}(\mu _{k}M_{1})$ such that for each potential $P\in \mathcal{F}_{E}(M)$, with $\overline{e}_{k}P\overline{e}_{k}=0$, we have that $\phi (\mu _{k}P)$ is cyclically equivalent to $\mu _{k}\varphi (P)$.
\end{theorem}
$  $

\begin{proof} Here $\overline{e}_{k}\varphi (P)\overline{e}_{k}=\varphi (\overline{e}_{k}P\overline{e}_{k})=0$; hence $\mu _{k}(\varphi(P))$ is defined. Let $X,X_{1}$ and $Y$ be local bases for $_{E}M,_{E}M_{1}$ and $M_{E}$, respectively. Then we have
$$\mu _{k}(P)=[P]+\sum _{x\in X_{k},y\in _{k}Y}([xy])(y^{*})(^{*}x).$$

By Proposition \ref{prop7.8}, there are algebra isomorphisms $\hat{\phi }:\mathcal{F}_{E}(\widetilde{M})\rightarrow \mathcal{F}_{E}(\widetilde{M_{1}})$ and $\phi :\mathcal{F}_{E}(\mu _{k}M)\rightarrow \mathcal{F}_{E}(\mu _{k}M_{1})$ such that
$$\hat{\phi }i_{M}=i_{M}\varphi ; \ \hat{\phi }i_{\mu _{k}M}=i_{\mu _{k}M_{1}}\phi .$$

We have 
$$\mu _{k}(P)=[P]+\sum _{x\in X_{k},y\in _{k}Y}(^{*}x)[xy](y^{*})+h$$
with $h\in [\mathcal{F}_{E}(\mu _{k}M),\mathcal{F}_{E}(\mu _{k}M)]$. Then 
$$i_{\mu _{k}M}\mu _{k}(P)=i_{\mu _{k}M }([P])+\sum _{x\in X_{k},y\in _{k}Y}(^{*}x)xyy^{*}+i_{\mu _{k}M}(h)$$
$$=i_{M }(P)+\sum _{x\in X_{k},y\in _{k}Y}(^{*}x)xyy^{*}+i_{\mu _{k}M}(h),$$
Applying the morphism $\hat{\phi }$ to the above equality we obtain
$$\hat{\phi }i_{\mu _{k}M}(\mu _{k}P)=i_{\mu _{k}M_{1}}(\phi (\mu _{k}P)) =$$
$$\hat{\phi }i_{M}(P)+\hat{\phi }\left(\sum _{x\in X_{k}}(^{*}x)x\right)\hat{\phi }\left(\sum _{y\in _{k}Y}yy^{*}\right)+\hat{\phi }i_{\mu _{k}M}(h)=$$
$$=i_{M_{1}}(\varphi (P))+\sum _{x_{1}\in (X_{1})_{k},y_{1}\in _{k}Y_{1}}(^{*}x_{1})x_{1}y_{1}y_{1}^{*}+i_{\mu _{k}M_{1}}(\phi (h))$$
$$=i_{\mu _{k}M_{1}}\left([\varphi (P)]+\sum _{x_{1}\in (X_{1})_{k}, y_{1}\in _{k}Y_{1}}(^{*}x_{1})[x_{1}y_{1}]y_{1}^{*}+\phi (h)\right).$$
Therefore
 $$\phi (\mu _{k}P)=[\varphi (P)]+\sum _{x_{1}\in (X_{1})_{k}, y_{1}\in _{k}Y_{1}}(^{*}x_{1})[x_{1}y_{1}]y_{1}^{*}+\phi (h).$$
 From here we deduce that $\phi (\mu _{k}P)$ is cyclically equivalent to $\mu _{k}\varphi (P)$. The proof is now complete.
\end{proof}
In order to define $\mu _{k}$ for an arbitrary potential we first define a continuous $F$-linear map
$$\kappa :\mathcal{F}_{E}(M)_{cyc}\rightarrow \overline{e}_{k}\mathcal{F}_{E}(M)\overline{e}_{k}$$
as follows. For $i\neq k$, we define $\kappa$ as the identity map on $e_{i}\mathcal{F}_{E}(M)e_{i}$, and for $i=k$
$$\kappa : e_{k}\mathcal{F}_{E}(M)e_{k}\rightarrow \overline{e}_{k}\mathcal{F}_{E}(M)\overline{e}_{k}$$  
is defined as $\kappa (f)=\displaystyle \sum _{y\in _{k}Y}y^{*}(f)y$ for all $f\in e_{k}\mathcal{F}_{E}(M)e_{k}$. Observe that $\kappa(f)$ is cyclically equivalent to $f.$  

\begin{definition} \label{def7.10} If $P$ is a potential in $\mathcal{F}_{E}(M)$, we define 
  $$\mu _{k}(P)=\mu _{k}(\kappa (P)).$$
\end{definition}
  
\begin{lemma} \label{lem7.11} Let $f$ and $g$ be directed elements such that $fg\in (\overline{e}_{k}\mathcal{F}_{E}(M)\overline{e}_{k})_{cyc}$, then 
  $$fg-gf=\sum _{i=1}^{l}(\alpha _{i}\beta _{i}-\beta _{i}\alpha _{i})$$
  for some elements $\alpha _{i}, \beta _{i}\in \overline{e}_{k}\mathcal{F}_{E}(M)\overline{e}_{k}$.
  \end{lemma}
  
\begin{proof} Since $f$ and $g$ are directed elements, then there exist idempotents $e_{i}, e_{j}$ such that
  $e_{i}fe_{j}=f$ and $e_{j}ge_{i}=g$. Then if $i\neq k, j\neq k$,
  $$\kappa (fg-gf)=\kappa (e_{i}fge_{i}-e_{j}gfe_{j})=fg-gf$$
  and our claim holds. If $i\neq k, j=k$, we have
  $$\kappa (fg-gf)=fg-\kappa (gf)$$
  $$=\sum _{y\in _{k}Y}fyy^{*}(g)-\sum _{y\in _{k}Y}y^{*}(g)fy=\sum _{y\in _{k}Y}(\alpha _{y}\beta _{y}-\beta _{y}\alpha _{y}),$$
  with $\alpha _{y}=fy=e_{i}fye_{\sigma (y)}, \beta _{y}=y^{*}(g)=e_{\sigma (y)}y^{*}(g)e_{i}$. From here we
 deduce our result for this case. If $i=k$ and $j\neq k$, we proceed in a similar way. Now suppose that $i=j=k$, then
  \begin{align*}
  \kappa (fg-gf)&=\sum _{y\in _{k}Y}y^{*}(fg)y-\sum _{y_{1}\in _{k}Y}y_{1}^{*}(gf)y_{1}\\
  &=\sum _{y,  y_{1}\in Y_{k}}(y^{*}(f)y_{1}y_{1}^{*}(g)y-y_{1}^{*}(g)yy^{*}(f)y_{1})\\
  &=\sum _{y,y_{1}\in _{k}Y}(\alpha _{y,y_{1}}\beta _{y,y_{1}}-\beta _{y,y_{1}}\alpha _{y,y_{1}})
  \end{align*}
  with $\alpha _{y,y_{1}}=y^{*}(f)y_{1}=e_{\sigma (y)}y^{*}(f)y_{1}e_{\sigma (y_{1})}$, and
  $\beta _{y,y_{1}}=y_{1}^{*}(g)y=e_{\sigma (y_{1}}y_{1}^{*}(g)ye_{\sigma (y)}$. Here $\sigma (y)\neq k$ and   $\sigma (y_{1})\neq k$. The result follows.  
  \end{proof}

\begin{prop} \label{prop7.12} If $P$ and $Q$ are cyclically equivalent potentials in $\mathcal{F}_{E}(M)$, then $\mu _{k}P$ is cyclically equivalent to $\mu _{k}Q$.
\end{prop}

\begin{proof} We have $P=Q+h$ with $h\in [\mathcal{F}_{E}(M),\mathcal{F}_{E}(M)]$; here $h=\displaystyle \lim_{l \to \infty }u_{l}$ where each $u_{l}$ is a finite sum of elements of the form $fg-gf$. We have $\kappa (h)=\displaystyle \lim_{l \to \infty }\kappa (u_{l})$. Now, consider the equality
$$i_{M}(\kappa (P))=i_{M}(\kappa (Q))+\displaystyle \lim_{l \to \infty }i_{M}(\kappa (u_{l})).$$
Each $i_{M}\kappa (u_{l})$ is a sum of elements of the form $i_{M}(\alpha )i_{M}(\beta )-i_{M}(\beta )i_{M}(\alpha )$.
Then by Lemma \ref{lem7.11}, each $i_{M}(\kappa (u_{l}))$ is a sum of elements of the form
$i_{\mu _{k}M}([\alpha ][\beta ]-[\beta ][\alpha ])$, so $i_{M}(\kappa (u_{l}))=i_{\mu _{k}M}([\kappa (u_{l})])$, with
$[\kappa (u_{l})]\in [\mathcal{F}_{E}(\mu _{k}M),\mathcal{F}_{E}(\mu _{k}M)]$.
We have

$$i_{\mu _{k}M}([\kappa(P)])=i_{\mu _{k}M}([\kappa(Q)])+i_{\mu _{k}M}(\displaystyle \lim_{l \to \infty }[\kappa (u_{l})]),$$ 
then
$$[\kappa (P)]=[\kappa (Q)]+h^{\prime }$$
with $h^{\prime }\in [\mathcal{F}_{E}(\mu _{k}M), \mathcal{F}_{E}(\mu _{k}M)]$. Therefore
$$\mu _{k}(P)=[\kappa (P)]+\sum _{x\in X_{k},y\in _{k}Y}[xy](y^{*})(^{*}x)=$$
$$[\kappa (Q)]+\sum _{x\in X_{k},y\in _{k}Y}[xy](y^{*})(^{*}x)+h^{\prime }=\mu _{k}(Q)+h^{\prime }.$$
This completes the proof. 
\end{proof}

\begin{rem} \label{rem7.13} Let $M_{1}, M_{2}$ be two $E$-$E$-bimodules such that for $s=1, 2$ whenever $e_{i}M_{s}e_{j}\neq 0$, then this bimodule is an $E_{i}$-free left module and an $E_{j}$-free right module. Suppose that $\overline{e}_{k}M_{2}\overline{e}_{k}=M_{2}$. Let $P$ be a potential in $\mathcal{F}_{E}(M_{1})$ and $Q$ a potential in $\mathcal{F}_{E}(M_{2})$. Then the potential $\mu _{k}(P\oplus Q)$ is right equivalent to $\mu _{k}(P)\oplus Q$. 
\end{rem}

Indeed, take local bases $X$ and $X^{\prime }$ for $_{E}M_{1}$ and $_{E}M_{2}$ respectively, and bases $Y$ and $Y^{\prime }$ for $(M_{1})_{E}$ and $(M_{2})_{E}$, respectively. Let $M=M_{1}\oplus _{M_{2}}$, then $X\cup X^{\prime }$ is a local basis for $_{E}M$ and $Y\cup Y^{\prime }$ is a local basis for $M_{E}$. We have $P\oplus Q\in \mathcal{F}_{E}(M)$,  $\mu _{k}M=\mu _{k}M_{1}\oplus M_{2}$ and  $X^{\prime }_{k}=\emptyset $. Then
\begin{align*}
\mu _{k}(P\oplus Q)&=[P]+[Q]+\sum _{x\in X_{k},y\in _{k}Y}[xy][y^{*}(^{*}x)] \\
&=[P]+\sum _{x\in X_{k},y\in _{k}Y}[xy][y^{*}(^{*}x)]+[Q] \\
&=\mu _{k}(P)+Q
\end{align*}
where $\mu _{k}(P)\in \mathcal{F}_{E}(\mu _{k}M)$ and $Q=[Q]\in \mathcal{F}_{E}(M_{2})$.

\end{section}
\begin{section}{The square of the premutation} \label{sec9}
In this section, the $F$-algebra and the $E$-$E$-bimodules $M$ are as in the previous section, without any assumption on the $F$-bases of $E$.
Let $P$ be a potential in $\mathcal{F}_{E}(M)$ where $e_{i}Me_{k}\neq 0$ implies $e_{k}Me_{i}=0$, and $e_{k}Me_{i}\neq 0 $ implies $e_{i}Me_{k}=0$. We have defined $\mu _{k}M$ and $\mu _{k}P$. From the definition of $\mu _{k}M$, one can see that $e_{i}\mu _{k}Me_{k}\neq 0$ implies $e_{k}\mu _{k}Me_{i}=0$, and that $e_{k}\mu _{k}Me_{i}\neq 0$ implies $e_{i}\mu _{k}Me_{k}=0$; therefore we can define $\mu _{k}^{2}M$ and $\mu _{k}^{2}P$. 
We have the $E$-$E$-bimodule
$$\mu _{k}^{2}M=\overline{e}_{k}(\mu _{k}M)\overline{e}_{k}\oplus (\mu _{k}M)e_{k}(\mu _{k}M)\oplus (e_{k}(\mu _{k}M))^{*}\oplus ^{*}((\mu _{k}M)e_{k})$$

From the definition of $\mu _{k}M$, we obtain
$$\overline{e}_{k}(\mu _{k}M)\overline{e}_{k}=\overline{e}_{k}M\overline{e}_{k}\oplus Me_{k}M,$$
$$(\mu _{k}M)e_{k}(\mu _{k}M)=(e_{k}M)^{*}e_{k}(^{*}(Me_{k})),$$
$$(e_{k}(\mu _{k}M))^{*}=(^{*}(Me_{k}))^{*}; \quad ^{*}((\mu _{k}M)e_{k})=^{*}((e_{k}M)^{*}).$$

By Proposition \ref{prop7.3}, there exists a morphism of topological algebras
$$\tau :\overline{e}_{k}\mathcal{F}_{E}(\mu _{k}M)\overline{e}_{k}\rightarrow \mathcal{F}_{E}(\mu _{k}^{2}M),$$
which induces an isomorphism of topological algebras:
$$\tau :\overline{e}_{k}\mathcal{F}_{E}(\mu _{k}M)\overline{e}_{k}\rightarrow \mathcal{F}_{E_{\overline{k}}}(\overline{e}_{k}(\mu _{k}M)\overline{e}_{k}\oplus (\mu _{k}M)e_{k}(\mu _{k}M)).$$

Consider the $E$-$E$-bimodule $\nu (M)=M\oplus Me_{k}M\oplus (e_{k}M)^{*}e_{k}(^{*}(Me_{k}))$. We will identify $M$ with its inclusion into the first summand of $\nu (M)$; the $E$-$E$-bimodule $Me_{k}M$ will be identified with its inclusion into the second summand of $\nu (M)$. We recall that the elements of $Me_{k}M \subseteq \mathcal{F}_{E}(\mu _{k}M)$ are of the form $[h]$, where $h\in Me_{k}\otimes _{E}e_{k}M \subseteq \mathcal{F}_{E}(M)$. Similarly, the $E$-$E$-bimodule $(e_{k}M)^{*}e_{k}(^{*}(Me_{k}))$ will be identified with its inclusion into the third direct summand of $\nu (M)$; the elements of the latter $E$-$E$-bimodule are of the form $\tau (f)$, where $f\in (e_{k}M)^{*}\otimes _{E}(^{*}(Me_{k}))\subseteq \mathcal{F}_{E}(\mu _{k}M)$.

\begin{theorem} \label{theo8.1} There exists an isomorphism of $E$-$E$-bimodules
$$\mathfrak{f}: \mathcal{F}_{E}(\mu _{k}^{2}M) \rightarrow \mathcal{F}_{E}(M\oplus Me_{k}M\oplus (e_{k}M)^{*}e_{k}^{*}(Me_{k}))$$
such that for any potential $P\in \mathcal{F}_{E}(M)^{\geq 3}$, $\mathfrak{f}(\mu _{k}^{2}(P))$ is cyclically equivalent to the potential
$$P+\sum _{x\in X_{k},y\in _{k}Y}[xy]\tau (y^{*}(^{*}x))+\sum _{x\in X_{k},x_{1}\in _{k}X}[xx_{1}]f_{[xx_{1}]}$$
where $f_{[xx_{1}]}\in \mathcal{F}_{E}(\overline{e}_{k}M\overline{e}_{k}\oplus Me_{k}M)^{\geq 1}$.
\end{theorem}
  
\begin{proof} Let $\nu (M)=M\oplus Me_{k}M\oplus (e_{k}M)^{*}e_{k}(^{*}(Me_{k}))$. We have an isomorphism of $E$-$E$-bimodules, $ev:e_{k}M\rightarrow ^{*}((e_{k}M)^{*})$, where for $m\in e_{k}M$ and $\psi \in (e_{k}M)^{*}$, $ev(m)(\psi )=\psi (m)$. The set consisting of all elements of the form $y^{*}$, where $y\in _{k}Y$, forms a free local basis for $_{E}(e_{k}M)^{*}$. Therefore, the set of elements $^{*}(y^{*})$, with $y\in _{k}Y$, is a free local basis of $^{*}((e_{k}M)^{*})_{E}$. Observe that $ev(y)=^{*}(y^{*})$. Similarly, we have the evaluation isomorphism $ev:Me_{k}\rightarrow (^{*}(Me_{k}))^{*}$, and $ev(x)=(^{*}x)^{*}$ for all $x\in X_{k}$. 

Now we define a morphism of $E$-$E$-bimodules $\mathfrak{g}^{(1)}:\mu _{k}^{2}M\rightarrow \nu (M)$. Let $h\in \overline{e}_{k}(\mu _{k}M)\overline{e}_{k}$, then if $h\in \overline{e}_{k}M\overline{e}_{k}$, $\mathfrak{g}^{(1)}(h)=h\in M\subseteq \nu (M)$; if $h=[h_{1}]\in Me_{k}M$ and $h_{1}\in Me_{k}\otimes _{E}e_{k}M$, $\mathfrak{g}^{(1)}([h_{1}])=[h_{1}]$. Finally, if $g=\tau (g_{1})\in (e_{k}M)^{*}e_{k}(^{*}(Me_{k}))$, with $g_{1}\in (e_{k}M)^{*}\otimes _{E}(^{*}(Me_{k}))$, $\mathfrak{g}^{(1)}(\tau (g_{1})))=\tau (g_{1})$.

Now for the elements in $^{*}((e_{k}M)^{*})$ and in $(^{*}(Me_{k}))^{*}$, we define $\mathfrak{g}^{(1)}$ as the inverse of the evaluation map. Clearly, $\mathfrak{g}^{(1)}$ is a monomorphism and here $M=\overline{e}_{k}M\overline{e}_{k}\oplus Me_{k}\oplus e_{k}M$, so $\mathfrak{g}^{(1)}$ is an isomorphism. Therefore, the pair
$(\mathfrak{g}^{(1)},0)$ induces an isomorphism
$$\mathfrak{g}:\mathcal{F}_{E}(\mu _{k}^{2}M)\rightarrow \mathcal{F}_{E}(\nu (M)).$$

Now let $P$ be any potential in $\mathcal{F}_{E}(M)^{\geq 3}$. We have
$\mu _{k}(P)=\mu _{k}(\kappa (P))$, so we may assume $P=\kappa (P)\in \overline{e}_{k}\mathcal{F}_{E}(M)\overline{e}_{k}$. We have $\mu _{k}(P)=[P]+\displaystyle \sum _{x\in X_{k},y\in _{k}Y}[xy](y^{*})(^{*}x)$, then
$$\mu _{k}^{2}(P)=\tau (\mu _{k}(P))+\sum _{x\in X_{k},y\in _{k}Y}\tau (y^{*}(^{*}x))((^{*}x)^{*})(^{*}(y^{*}))=$$
$$\tau ([P])+\sum _{x\in X_{k},y\in _{k}Y}\tau ([xy])\tau ((y^{*})(^{*}x))+\sum _{x\in X_{k},y\in _{k}Y}\tau (y^{*}(^{*}x))\tau ((^{*}x)^{*})\tau (^{*}(y^{*})).$$

Therefore
$$\mathfrak{g}(\mu _{k}^{2}(P))=\mathfrak{g}(\tau ([P]))+\sum _{x\in X_{k},y\in _{k}Y}[xy]\tau ((y^{*})(^{*}x))+\sum _{x\in X_{k},y\in _{k}Y}\tau (y^{*}(^{*}x))(xy).$$
This implies that the potential $\mathfrak{g}(\mu _{k}^{2}(P))$ is cyclically equivalent to the potential
$$\mathfrak{g}(\tau [P])+\sum _{x\in X_{k},y \in _{k}Y}([xy]+xy)\tau (y^{*}(^{*}x)).$$

Since $M=\overline{e}_{k}M\overline{e}_{k}\oplus Me_{k}\oplus e_{k}M$, we define an automorphism of $E$-$E$-bimodules as the identity map on the first and third summands; and as $-id _{Me_{k}}$ on $Me_{k}$. We now extend this automorphism to an automorphism $\varphi ^{(1)}$ of $M\oplus Me_{k}M\oplus (e_{k}M)^{*}e_{k}(^{*}(Me_{k}))$, defined as the identity map on the third and fourth summands. Let $\varphi $ be the algebra automorphism
of $\mathcal{F}_{E}(\nu (M))$ determined by the pair $(\varphi ^{(1)},0)$. Then
$$\varphi \mathfrak{g}(\mu _{k}^{2}(P))=\mathfrak{g}(\tau ([P]))+\sum _{x\in _{k}X, y\in _{k}Y}([xy]-xy)\tau (y^{*}(^{*}x)).$$
Let 
$$\phi ^{(2)}:M\oplus Me_{k}M\oplus (e_{k}M)^{*}\oplus (^{*}(Me_{k})e_{k}(^{*}(Me_{k}))\rightarrow \mathcal{F}_{E}(\nu (M))^{\geq 2}$$
be the morphism which is the identity on the first and third summands, and on $Me_{k}M$ is the isomorphism $[-]^{-1}:Me_{k}M\rightarrow 
Me_{k}\otimes _{E}e_{k}M $ followed by the composition of inclusions 
$$Me_{k}\otimes _{E}e_{k}M\subseteq \mathcal{F}_{E}(M)^{\geq 2}\subseteq \mathcal{F}_{E}(\nu (M))^{\geq 2}.$$
Now let $\phi$ be the automorphism of $\mathcal{F}_{E}(\nu (M))$ determined by the pair $(id_{\nu (M)},\phi ^{(2)})$. Define $\mathcal{J}$ as the closure of the two-sided ideal of $\mathcal{F}_{E}(\overline{e}_{k}M\overline{e}_{k}\oplus Me_{k}M)$ generated by the elements $[xx_{1}]$ with $x\in X_{k}, x_{1}\in _{k}X$.

$ $
Recall we have 
$$[-]:\overline{e}_{k}\mathcal{F}_{E}(M)\overline{e}_{k}\rightarrow \mathcal{F}_{E_{\overline{k}}}(\overline{e}_{k}M\overline{e}_{k}\oplus Me_{k}M),$$
$$\mathcal{F}_{E_{\overline{k}}}(\overline{e}_{k}M\overline{e}_{k}\oplus Me_{k}M)\subseteq \overline{e}_{k}\mathcal{F}_{E}(\mu _{k}M)\overline{e}_{k}$$
{\bf Claim. } If $h\in \mathcal{F}_{E_{\overline{k}}}(\overline{e}_{k}M\overline{e}_{k}\oplus Me_{k}M)$, then
$$\phi \varphi \mathfrak{g}(\tau (h))=[-]^{-1}(h)+z,$$
with $z\in \mathcal{J}$. 

$ $

{\bf Proof of the Claim } Suppose first that $h=[m]$, with $m\in \overline{e}_{k}M\overline{e}_{k}$. Then $\tau (h)=h$, $\mathfrak{g}(h)=m$, and $\phi \varphi \mathfrak{g}(h)=m=[-]^{-1}(h)$. If $h=[f]$, with $f\in Me_{k}\otimes _{E}e_{k}M$, then $\tau (h)=h$, $\mathfrak{g}(h)=h$ and 
$\phi \varphi \mathfrak{g}(\tau (h))=f+h$.  Since the set of all elements of the form $[xx_{1}]$ with $x\in X_{k}, x_{1}\in _{k}X$, is a set of generators for $_{E}Me_{k}M$, then $h=[f]\in \mathcal{J}$; hence $\phi \varphi \mathfrak{g}(\tau (h))=[-]^{-1}(h)+h$. Therefore, our claim holds for all the generators of the algebra
$T_{E}(\overline{e}_{k}M\overline{e}_{k}\oplus Me_{k}M)$ and hence it holds for the whole algebra.
Now if $h\in \mathcal{F}_{E}(\overline{e}_{k}M\overline{e}_{k}\oplus Me_{k}M)$, then $h=\displaystyle \lim_{l \to \infty }h_{l}$ with $h_{l}\in T_{E}(\overline{e}_{k}M\overline{e}_{k}\oplus Me_{k}M)$. Therefore
$$\phi \varphi \mathfrak{g}(\tau (h))=\displaystyle \lim_{l \to \infty}\phi \varphi \mathfrak{g}(\tau (h_{l}))$$
and $\phi \varphi \mathfrak{g}(\tau (h_{l}))=[-]^{-1}(h_{l})+z_{l}$, with $z_{l}\in \mathcal{J}$. Here the sequences $\{[-]^{-1}(h_{l}+z_{l})\}_{l=1}^{\infty }$ and $\{[-]^{-1}(h_{l})\}_{l=1}^{\infty }$ converge; hence the sequence
$\{z_{l}\}_{l=1}^{\infty }$ converges as well. Therefore
$$\phi \varphi \mathfrak{g}(\tau (h))=[-]^{-1}(h)+\displaystyle \lim_{l \to \infty }z_{l}.$$

This proves our claim.

Now we continue with the proof of our proposition. Since the set of all $[xx_{1}]$ with $x\in X_{k},x_{1}\in _{k}X$ is a set of generators for $_{E}Me_{k}M$, there exist elements $\alpha _{[x,y],x^{\prime }x_{1}}\in E$ such that
$$[xy]-xy=\displaystyle \sum _{x,x^{\prime }\in X_{k},x_{1}\in _{k}X}\alpha _{[xy],x^{\prime }}([x^{\prime }x_{1}]-x^{\prime }x_{1})$$. 

Then
\begin{align*}
\phi \varphi \mathfrak{g}(\mu _{k}^{2}(P))&=P+\sum _{x,x^{\prime }\in X_{k},y\in _{k}Y}\alpha _{[xy],x^{\prime }x_{1}}\phi ([x^{\prime }x_{1}]-x^{\prime }x_{1})\tau (y^{*}(^{*}x))+z\\
&=P+ \sum _{x,x^{\prime }\in X_{k},y\in _{k}Y}\alpha _{[xy],x^{\prime }x_{1}}[x^{\prime }x_{1}]\tau (y^{*}(^{*}x))+z\\
&=P+\sum _{x\in X_{k},y\in _{k}Y}[xy]\tau (y^{*}(^{*}x))+z
\end{align*}
with $z\in \mathcal{J}$. 

We have  $\mathcal{J}=\mathcal{F}_{E}(\overline{e}_{k}M\overline{e}_{k}\oplus Me_{k}M)\mathcal{J}$, then Lemma \ref{lem4.4} yields that $z$ is cyclically equivalent to an element of the form
$\displaystyle \sum _{x\in X_{k},x_{1}\in _{k}X}f_{[xx_{1}]}[xx_{1}]$, where $f_{[xx_{1}]}\in \mathcal{F}_{E}(\overline{e}_{k}M\overline{e}_{k}\oplus Me_{k}M)$. Since $\phi \varphi \mathfrak{g}(\mu _{k}^{2}(P))-P-\displaystyle \sum _{x\in X_{k},y\in _{k}Y}[xy]\tau (y^{*}(^{*}x))$ lies in $\mathcal{F}_{E}(\nu (M))^{\geq 2}$, then 
$f_{[xx_{1}]}\in \mathcal{F}_{E}(\overline{e}_{k}M\overline{e}_{k}\oplus Me_{k}M)^{\geq 1}$. Taking $\mathfrak{f}=\phi \varphi \mathfrak{g}$ yields the desired result.
\end{proof}
We recall that a $n\times n$ matrix $B$ with integer entries is called skew-symmetrizable if there exists a diagonal matrix $D=\operatorname{diag}(d_{1},\dots,d_{n})$, with $d_{1},\dots,d_{n}$ positive integers, such that the matrix $DB$ is skew-symmetric.

\begin{definition} \label{def8.2} An $E$-$E$-bimodule $M$  is said to be $2$-acyclic if $e_{i}Me_{j}\neq 0$ implies $e_{j}Me_{i}=0.$ 
 \end{definition} 
 
To any $E$-$E$-bimodule $M$ with $M_{cyc}=0$ (not necessarily $2$-acyclic) and such that  each $e_{i}Me_{j}\neq 0$ is a free left $E_{i}$-module and also a free right $E_{j}$-module, we associate the matrix $B(M)=(b_{i,j})$, where
 $$b_{i,j}=\mathrm{rank}_{E_{i}}(e_{i}Me_{j})-\mathrm{rank}_{E_{i}}(e_{j}Me_{i}).$$
 
\begin{rem} \label{rem8.3} The matrix $B(M)$ is skew-symmetrizable. Indeed, take $d_{i}=\mathrm{dim}_{F}E_{i}$, then
$$d_{i}b_{i,j}=d_{i}\mathrm{rank}_{E_{i}}(e_{i}Me_{j})-d_{i}\mathrm{rank}_{E_{i}}(e_{j}Me_{i})=\mathrm{dim}_{F}(e_{i}Me_{j})-\mathrm{dim}_{F}(e_{j}Me_{i})$$

Then $d_{j}b_{j,i}=\mathrm{dim}_{F}(e_{j}Me_{i})-\mathrm{dim}_{F}(e_{i}Me_{j})=-d_{i}b_{i,j}$.
\end{rem}
  
\begin{prop} \label{prop8.4} Let $M$ be an $E$-$E$-bimodule as before, and such that $e_{i}Me_{k}\neq 0$ implies
 $e_{k}Me_{i}=0$, and $e_{k}Me_{i}\neq 0$ implies $e_{i}Me_{k}=0$. Then $\mu _{k}M$ is defined and
 $$B(\mu _{k}M)=\mu _{k}B(M),$$
 where $\mu _{k}B(M)$ is the matrix mutation of $B$ in the direction of $k$, in the sense of Fomin-Zelevinsky \cite{5}.
 \end{prop}
  
\begin{proof} Take $(c_{i,j})=B(\mu _{k}M)$. Then for $i\neq k, j\neq k$
 \begin{align*}
 c_{i,j}&=\mathrm{rank}_{E_{i}}(e_{i}\mu _{k}Me_{j})-\mathrm{rank}_{E_{i}}(e_{j}\mu _{k}Me_{i}) \\
&=\mathrm{rank}_{E_{i}}(e_{i}Me_{j})-\mathrm{rank}_{E_{i}}(e_{j}Me_{i})+\mathrm{rank}_{E_{i}}(e_{i}Me_{k}Me_{j})
-\mathrm{rank}_{E_{i}}(e_{j}Me_{k}Me_{i}) \\
&=b_{i,j}+\mathrm{rank}_{E_{i}}(e_{i}Me_{k})\mathrm{rank}_{E_{k}}(e_{k}Me_{j})-\mathrm{rank}_{E_{k}}(e_{j}Me_{k})\mathrm{rank}_{E_{i}}(e_{k}Me_{i})
\end{align*}
Then if $b_{i,k}>0$ and $b_{k,i}>0$, we have
$$c_{i,j}=b_{i,j}+b_{i,k}b_{k,j}.$$
If $b_{i,k}<0$ and $b_{k,j}<0$, then
$$c_{i,j}=b_{i,j}-b_{i,k}b_{k,j}.$$
Now suppose that $b_{i,k}b_{k,j}\leq 0$. If $b_{i,k}\leq 0$, then $b_{k,j}\geq 0$, which implies that $e_{i}Me_{k}=0$ and $e_{j}Me_{k}=0$. If $b_{i,k}\geq 0$, then $b_{k,j}\leq 0$. Thus $e_{k}Me_{i}=0$ and $e_{k}Me_{j}=0$. In any case
$$c_{i,j}=b_{i,j}.$$
Now 
\begin{align*}
c_{k,j}&=\mathrm{rank}_{E_{k}}(^{*}(Me_{k})e_{j})-\mathrm{rank}_{E_{k}}(e_{j}(e_{k}M)^{*}) \\
&=\mathrm{rank}_{E_{k}}\left(^{*}(e_{j}Me_{k})\right)-\mathrm{rank}_{E_{k}}((e_{k}Me_{j})^{*}) \\
&=\mathrm{rank}_{E_{k}}(e_{j}Me_{k})-\mathrm{rank}_{E_{k}}(e_{k}Me_{j})=-b_{k,j}
\end{align*}
Similarly, one can prove the equality $c_{j,k}=-b_{k,j}$. The assertion follows.
\end{proof}
\end{section}

\begin{section}{Realizations} \label{sec10}
In this section we examine a procedure proposed in \cite{6,8} for constructing an $F$-algebra $E$ and an $E$-$E$-bimodule $M$, such that the matrix $B(M)$ of $M$ is skew-symmetrizable. The algebra $E=\displaystyle \prod_{i=1}^{n} E_{i}$ considered in \cite{6,8} has the property that each $F$-algebra $E_{i}$ is a cyclic field extension of $F$. Based on this procedure, and using the construction given in \cite{3}, we give a construction using the same arithmetic given in \cite{6}, but allowing $F$ to be an algebraically closed field.

Let $E=\displaystyle \prod_{i=1}^{n} E_{i}$ and set $\operatorname{dim}_{F} E_{i}=d_{i}$. We put $d_{i,j}=(d_{i},d_{j})$ and $d_{i,j,l}=(d_{i},d_{j},d_{l})$. For each $i,j$ suppose that we have a family $\{M_{\rho}\}_{\rho \in \mathcal{B}(i,j)}$ of $E_{i}$-$E_{j}$-bimodules with the following conditions
\begin{enumerate}
\item For $\rho \in \mathcal{B}(i,j)$, $M_{\rho}$ is free as a left $E_{i}$-module and also free as a right $E_{j}$-module.
\item For $\rho \in \mathcal{B}(i,j)$, $\operatorname{rank}_{E_{i}}M_{\rho}=\frac{d_{j}}{d_{i,j}}$.
\item There exists a bijection $\mathcal{B}(i,j) \rightarrow \mathcal{B}(j,i)$ sending $\rho$ into $\hat{\rho}$ such that $M_{\rho}^{\ast} \cong M_{\hat{\rho}}$.
\item If $\rho_{1} \in \mathcal{B}(i,j)$ and $\rho_{2} \in \mathcal{B}(j,l)$, then
\begin{center}
$M_{\rho_{1}} \otimes_{E_{j}} M_{\rho_{2}} = \displaystyle \left(\bigoplus_{\rho \in \mathcal{B}(\rho_{1},\rho_{2})} M_{\rho} \right)^{\frac{d_{i}d_{i,j,l}}{d_{i,j}d_{j,l}}}$
\end{center}
for some subset $\mathcal{B}(\rho_{1},\rho_{2})$ of $\mathcal{B}(i,l)$. 
\end{enumerate}
In the above situation we set $\mathcal{B}=\displaystyle \bigcup_{i,j} \mathcal{B}(i,j)$.

Let $Q=(Q_{0},Q_{1},h,t)$ be a quiver without loops nor $2$-cycles; also, suppose that given $i,j \in Q_{0}$ there is at most one arrow $\alpha \in Q_{1}$ with $t(\alpha)=i$ and $h(\alpha)=j$. 

\begin{definition} \label{def9.1} Suppose $\{M_{\rho }\}_{\rho \in \mathcal{B}}$ is a family of $E$-$E$-bimodules satisfying the above conditions. Then a modulating function of $Q$ in $\mathcal{B}$ is a function
$$g:Q_{1}\times \mathcal{B}\rightarrow \{0,1,\ldots\}$$
such that if $g(\alpha , \rho )\neq 0$, then $\rho \in \mathcal{B}(h(\alpha ), t(\alpha ))$. For
$\alpha \in Q_{1}$, we put $|g(\alpha )|=\sum _{\rho \in \mathcal{B}}g(\alpha , \rho )$.
\end{definition}
Associated to a modulating function $g$ of $Q$, we define the $E$-$E$-bimodule
$$M_{g}=\bigoplus _{\alpha \in Q_{1}, \rho \in \mathcal{B}}M_{\rho }^{g(\alpha , \rho )}.$$

$ $

Let $B(M_{g})=(b_{i,j}(M_{g}))$ be the corresponding skew-symmetrizable matrix. Then
$$b_{i,j}(M_{g})=(u_{i,j}-u_{j,i})d_{j}/d_{i,j}$$
where $u_{s,t}=|g(\alpha )|$ if there is an arrow $\alpha :t\rightarrow s$, otherwise $u_{s,t}=0$.

Now suppose the quiver $Q$ has no $2$-cycles and $g$ is a modulating function for $Q$ on the family of $E$-$E$-bimodules $\{M_{\rho }\}_{\rho \in \mathcal{B}}$. Then if $k\in Q_{0}$, we define $\mu _{k}Q$ taking $Q$ and reversing the orientation of the arrows $\alpha $ with $h(\alpha )=k$ or $t(\alpha )=k$; that is, we change $\alpha $ by $\tilde{\alpha }$ such that $h(\tilde{\alpha })=t(\alpha )$ and $t(\tilde{\alpha })=h(\alpha )$. Moreover, if $\alpha :i\rightarrow k$ and $\beta :k\rightarrow j$ and there is no arrow $i\rightarrow j$, then we add an arrow $[\beta \alpha ]:i\rightarrow j$. Associated to the modulation function $g$ we define a modulation function $\mu _{k}g$ on
for the quiver $\mu _{k}Q$, on the same family $\{M_{\rho }\}_{\rho \in \mathcal{B}}$, as follows

\begin{itemize}
\item For $\alpha: j \rightarrow i$ with $j \neq k$ and $i \neq k$, such that there is no pair of arrows $\beta: j \rightarrow k$, $\gamma: k \rightarrow i$ we define $\mu_{k}g(\alpha,\rho)=g(\alpha,\rho)$ for all $\rho \in \mathcal{B}$. If there is a pair of arrows $\beta: j \rightarrow k$, and $\gamma: k \rightarrow i$, we define
\begin{center}
$\mu_{k}g(\alpha,\rho)=g(\alpha,\rho)+\displaystyle \sum_{\rho_{1},\rho_{2}| \rho \in \mathcal{B}(\rho_{1},\rho_{2})}g(\gamma,\rho_{1})g(\beta,\rho_{2})d_{k}d_{i,k,j}/d_{i,k}d_{k,j}.$
\end{center}
\item If $\alpha$ is an arrow in $Q$ with $h(\alpha)=k$ or $t(\alpha)=k$, $\mu_{k}g(\hat{\alpha},\hat{\rho})=g(\alpha,\rho)$.
\item Finally, if $\alpha: i \rightarrow k$ and $\beta: k \rightarrow j$ are arrows of $Q$, and there is no arrow in $Q$ from $i$ to $j$ we define
\begin{center}
$g([\beta \alpha],\rho)=\displaystyle \sum_{\rho_{1},\rho_{2}| \rho \in \mathcal{B}(\rho_{1},\rho_{2})} g(\beta,\rho_{1})g(\alpha,\rho_{2})d_{k}d_{i,k,j}/d_{i,k}d_{k,j}.$
\end{center}
\end{itemize}
Let $B \in \mathbb{Z}^{n \times n}$ be any skew-symmetrizable matrix. The quiver of $B$, $Q_{B}$, has vertices $\{1,\ldots,n\}$ and an arrow $\alpha: j \rightarrow i$ if and only if $b_{i,j}>0$. Since $B$ is skew-symmetrizable, $Q_{B}$ has no loops nor $2$-cycles. Now let $E=\displaystyle \prod_{i=1}^{n} E_{i}$ and let $\{M_{\rho}\}_{\rho \in \mathcal{B}}$ be a family of $E$-$E$-bimodules satisfying properties $1$-$4$. Choose a modulating function $g$ for $Q_{B}$ on $\mathcal{B}$ such that for any $\alpha: j \rightarrow i \in (Q_{B})_{1}$, $|g(\alpha)|=b_{i,j}d_{i,j}/d_{j}$, then in case $b_{i,j}>0$, $b_{i,j}(M_{g})=b_{i,j}$; and $b_{j,i}(M_{g})=-|g(\alpha)|d_{i}/d_{j,i}=b_{j,i}$. Therefore $B(M_{g})=B$. 

\begin{prop} \label{prop9.3} If $g$ is a modulating function for the quiver $Q$ on the family $\{M_{\rho}\}_{\rho \in \mathcal{B}}$, then
\begin{center}
$\mu_{k}M_{g} \cong M_{\mu_{k}g}$
\end{center}
\end{prop}

\begin{proof} Let $\alpha: j \rightarrow k$ an arrow in $Q$. Then
$$e_{j}(\mu_{k}M_{g})e_{k}=(e_{k}M_{g}e_{j})^{\ast}=\displaystyle \left(\bigoplus_{\rho \in \mathcal{B}(k,j)}M_{\rho}^{g(\alpha,\rho)}\right)^{\ast} \cong \displaystyle \bigoplus_{\rho \in \mathcal{B}(j,k)}M_{\hat{\rho}}^{g(\alpha,\rho)}=e_{j}M_{\mu_{k}g}e_{k}$$
Similarly, if $\alpha: k \rightarrow j$ is an arrow in $Q$, then
\begin{center}
$e_{k}\mu_{k}M_{g}e_{j} \cong e_{k}M_{\mu_{k}g}e_{j}$.
\end{center}
Now let $i,j$ be vertices of $Q$, both distinct from $k$, and suppose that there is an arrow $\alpha: j \rightarrow i$. Then, if there is no pair of arrows $\beta: j \rightarrow k$, $\gamma: k \rightarrow i$, we have: $e_{i}\mu_{k}M_{g}e_{j}=e_{i}M_{g}e_{j}=\displaystyle \bigoplus_{\rho \in \mathcal{B}(i,j)} M_{\rho}^{g(\alpha,\rho)}=e_{i}M_{\mu_{k}g}e_{j}$. Now if there is a pair of arrows $\beta: j \rightarrow k$, $\gamma: k \rightarrow i$, taking $c=d_{k}d_{i,k,j}/d_{i,k}d_{k,j}$ we have
\begin{align*}
e_{i}\mu_{k}M_{g}e_{j}&=e_{i}M_{g}e_{j} \oplus e_{i}M_{g}e_{k}  \otimes_{E_{k}} e_{k}M_{g}e_{j}\\
&=\displaystyle \bigoplus_{\rho \in \mathcal{B}(i,j)} \left(M_{\rho}^{g(\alpha,\rho)} \displaystyle \bigoplus_{\rho_{1},\rho_{2} | \rho \in \mathcal{B}(\rho_{1},\rho_{2})} M_{\rho}^{g(\gamma,\rho_{1})g(\beta,\rho_{2})c}\right)\\
&=\displaystyle \bigoplus_{\rho \in \mathcal{B}(i,j)} M_{\rho}^{\mu_{k}g(\alpha,\rho)}\\
&=e_{i}M_{\mu_{k}g}e_{j}.
\end{align*}
Similarly, if there are arrows in $Q$, $\beta: j \rightarrow k$, $\gamma: k \rightarrow i$ and no arrow $\alpha: j \rightarrow i$, we have
\begin{center}
$e_{i}\mu_{k}M_{g}e_{j}=e_{i}M_{\mu_{k}g}e_{j}$
\end{center}
This proves our result.
\end{proof}
We now give two examples of families of bimodules satisfying conditions $1$-$4$.
\begin{example} \label{ex9.4}
Let $(d_{1},\ldots,d_{n})$ be a collection of positive integers and let $d=\operatorname{lcm}(d_{1},\ldots,d_{n})$. Let $F$ be a field containing a primitive $d^{th}$-root of unity, and let $G/F$ be a cyclic Galois extension of degree $d$. For each divisor $l$ of $d$, we take the unique intermediate field $F \leq F_{l} \leq G$ having degree $l$. Then take the fields $F_{i}:=F_{d_{i}}$ and the semisimple commutative $F$-algebra $E=\displaystyle \prod_{i=1}^{n} F_{i}$. For $\rho \in \operatorname{Gal}(F_{i,j}/F)$, where $F_{i,j}:=F_{i} \cap F_{j}$, we consider the $F_{i,j}$-bimodule $F_{i}^{\rho}$ given by $F_{i}$ endowed with the left structure given by multiplication in $F_{i}$; the right structure is defined as follows: for $x \in F_{i}$ and $y \in F_{i,j}$, we define $x \ast y=x\rho(y)$. For each $i,j$ consider the family of $F_{i}$-$F_{j}$-bimodules $\{F_{i}^{\rho} \otimes_{F_{i,j}} F_{j}\}_{\rho \in \operatorname{Gal}(F_{i,j}/F)}$. In \cite[pp.12-13]{6} it is proved that the above family satisfies conditions $1$-$4$. 
\end{example}

\begin{example} \label{ex9.5}
As above, let $(d_{1},\ldots,d_{n})$ be a collection of positive integers and let $d=\operatorname{lcm}(d_{1},\ldots,d_{n})$. Let $F$ be a field containing a primitive $d^{th}$-root of unity. Let $G$ be a cyclic group of order $d$; and for each divisor $l$ of $d$, denote by $G_{l}$, the unique subgroup of $G$ of order $l$. Now take the group algebras $E_{i}:=FG_{d_{i}}$ and $E=\displaystyle \prod_{i=1}^{n} E_{i}$. For each subgroup $H \leq G$, denote by $\chi(H)$ the set of all irreducible characters of $H$ over $F$. Since $F$ contains a primitive $d^{th}$-root of unity, all irreducible characters of $H$ over $F$ are one-dimensional. Let $\mathcal{S}$ be the set of all subgroups of $G$ and $\chi(\mathcal{S})=\displaystyle \bigcup_{H \in \mathcal{S}} \chi(H)$. 
Given subgroups $G_{i},G_{j},G_{l}$ of $G$, we set $G_{i,j}:=G_{i} \cap G_{j}$ and $G_{i,j,l}:=G_{i} \cap G_{j} \cap G_{l}$. As in the above example, for a given $\rho \in \chi(G_{i,j})$ we define an $FG_{i}$-$FG_{i,j}$-bimodule $FG_{i}^{\rho}$ by taking $FG_{i}$ endowed with the left structure given by the multiplication of $FG_{i}$; and for $x \in FG_{i}$ and $y \in FG_{i,j}$ ,we set $x \ast y=xy\rho(y)$. Then as we will see below, the following family of $FG_{i}$-$FG_{j}$-bimodules: $\{FG_{i}^{\rho} \otimes_{FG_{i,j}} FG_{j}\}_{\rho \in \chi(G_{i,j})}$ satisfies conditions $1$-$4$. 
\end{example}

For each $\rho \in \chi(G_{i} \cap G_{j})$, consider the map $\underline{\rho}: F(G_{i} \cap G_{j}) \rightarrow F(G_{i} \cap G_{j})$ given by $\underline{\rho}(g)=\rho(g)g$ for all $g \in G_{i} \cap G_{j}$. Given $\rho \in \chi(G_{i}Ê\cap G_{j})$, we consider the element
\begin{center}
$e_{\rho}=(1/d_{i,j}) \displaystyle \sum_{w \in G_{i} \cap G_{j}} \underline{\rho}(w^{-1}) \otimes w \in FG_{i} \otimes_{F} FG_{j}$
\end{center}
Observe that for $g \in G_{i} \cap G_{j}$
\begin{align*}
e_{\rho}g&=(1/d_{i,j})\displaystyle \sum_{w \in G_{i} \cap G_{j}} \rho(w^{-1})w^{-1} \otimes wg \\
&=(1/d_{i,j}) \displaystyle \sum_{w_{1} \in G_{i} \cap G_{j}} \rho(gw_{1}^{-1})gw_{1}^{-1} \otimes w_{1}\\
&=(1/d_{i,j})\rho(g)g\left(\displaystyle \sum_{w_{1} \in G_{i} \cap G_{j}} \rho(w_{1}^{-1})w_{1}^{-1} \otimes w_{1}\right)\\
&=\underline{\rho}(g)e_{\rho}.
\end{align*}
\begin{prop} \label{prop9.6} Suppose $\rho_{1}, \rho_{2} \in \chi(G_{i} \cap G_{j})$. Then
\begin{enumerate}[(a)]
\item $e_{\rho_{1}}e_{\rho_{2}}=\delta_{\rho_{1},\rho_{2}}e_{\rho_{1}}.$
\item $1Ê\otimes 1 = \displaystyle \sum_{\rho \in \chi(G_{i} \cap G_{j})} e_{\rho}.$
\end{enumerate}
\end{prop}
\begin{proof}
Taking into account the above observation we have
\begin{align*}
e_{\rho_{1}}e_{\rho_{2}}&=(1/d_{i,j}) \displaystyle \sum_{w_{1} \in G_{i} \cap G_{j}} \underline{\rho_{1}}(w_{1}^{-1})e_{\rho_{2}}w_{1} \\
&=(1/d_{i,j}) \displaystyle \sum_{w_{1} \in G_{i} \cap G_{j}} \underline{\rho_{1}}(w_{1}^{-1})\underline{\rho_{2}}(w_{1})e_{\rho_{2}} \\
&=(1/d_{i,j}) \displaystyle \sum_{w_{1} \in G_{i} \cap G_{j}} e_{\rho_{2}}\rho_{2}(w_{1})\rho_{1}(w_{1}^{-1}) \\
&=\delta_{\rho_{1},\rho_{2}}e_{\rho_{1}}
\end{align*}
This proves item (a). For item (b), one has
\begin{align*}
\displaystyle \sum_{\rho \in \chi(G_{i} \cap G_{j})} e_{\rho}&=(1/d_{i,j}) \displaystyle \sum_{wÊ\in G_{i} \cap G_{j}} \displaystyle \sum_{\rho \in \chi(G_{i} \cap G_{j})} \rho(w^{-1})w^{-1} \otimes w \\
&=(1/d_{i,j}) \displaystyle \sum_{w \in G_{i} \cap G_{j}} \operatorname{reg}(w^{-1})w^{-1} \otimes w
\end{align*}
where reg is the regular character of $G_{i} \cap G_{j}$. Now recalling that $\operatorname{reg}(w)=0$ if $w \neq 1$ and $\operatorname{reg}(1)=d_{i,j}$, we obtain
\begin{center}
$\displaystyle \sum_{\rho \in \chi(G_{i} \cap G_{j})} e_{\rho}=1 \otimes 1$
\end{center}
This completes the proof.
\end{proof}
As in \cite[p.11]{6} we define the $FG_{i}$-$F(G_{i} \cap G_{j})$-bimodule $FG_{i}^{\rho}$, endowed with the left structure given by $FG_{i}$ and the right structure given by $x \ast g = x\underline{\rho}(g)$ for $g \in G_{i} \cap G_{j}$. We now consider the $FG_{i}$-$FG_{j}$-bimodule $FG_{i}^{\rho} \otimes_{F(G_{i} \cap G_{j})} FG_{j}$.
Observe that there exists an epimorphism of $FG_{i}$-$FG_{j}$-bimodules $FG_{i}^{\rho} \otimes_{F(G_{i} \cap G_{j})} FG_{j} \rightarrow FG_{i}e_{\rho}FG_{j}$ given by $x \otimes y \mapsto xe_{\rho}y$. We now show this map is in fact an isomorphism. \\

\begin{rem} \label{rem9.7} There exists an isomorphism $FG_{i}^{\rho} \otimes_{F(G_{i} \cap G_{j})} FG_{j} \rightarrow FG_{i}e_{\rho}FG_{j}$. 
\end{rem}

\begin{proof} Using (b) of Proposition \ref{prop9.6} we obtain
\begin{center}
$FG_{i} \otimes_{F} FG_{j} = \displaystyle \bigoplus_{\rho \in \chi(G_{i} \cap G_{j})} FG_{i}e_{\rho}FG_{j}$.
\end{center}
Therefore
\begin{align*}
d_{i}d_{j}&=\operatorname{dim}_{F} FG_{i}Ê\otimes_{F} FG_{j} \\
&=\displaystyle \sum_{\rho \in \chi(G_{i} \cap G_{j})} \operatorname{dim}_{F}(FG_{i}e_{\rho}FG_{j}) \\
& \leq \displaystyle \sum_{\rho \in \chi(G_{i} \cap G_{j})} \operatorname{dim}_{F} (FG_{i}^{\rho} \otimes_{F(G_{i} \cap G_{j})} FG_{j}) \\
&=\displaystyle \sum_{\rho \in \chi(G_{i} \cap G_{j})} d_{i}d_{j}/d_{i,j}  \\
&=d_{i}d_{j}
\end{align*}
hence $\operatorname{dim}_{F}(FG_{i}e_{\rho}FG_{j})=\operatorname{dim}_{F}(FG_{i}^{\rho} \otimes_{F(G_{i} \cap G_{j})} FG_{j})$, and the claim follows.
\end{proof}

\begin{lemma} \label{lem9.8} Let $G_{i},G_{j},G_{l}$ be subgroups of the cyclic group $G$. Take $G_{i,j}=G_{i} \cap G_{j}$, $G_{j,l}=G_{j} \cap G_{l}$, $G_{i,j,l}=G_{i} \cap G_{j} \cap G_{l}$. Suppose $\rho_{1}$ is an irreducible character of $G_{i,j}$ and $\rho_{2}$ is an irreducible character of $G_{j,l}$. Then
\begin{center}
$(FG_{i}^{\rho_{1}} \otimes_{FG_{i,j}} FG_{j}) \otimes_{FG_{j}} (FG_{j}^{\rho_{2}} \otimes_{FG_{j,l}} FG_{l}) \cong  (FG_{i}^{\lambda} \otimes_{FG_{i,j,l}} FG_{l})^{[G_{j}: G_{i,j}G_{j,l}]}$
\end{center}
where $\lambda=\rho_{1}|\rho_{2}|$, with $\rho_{1}|$ and $\rho_{2}|$ being the respective restrictions of $\rho_{1}$ and $\rho_{2}$ to $G_{i,j,l}$. 
\end{lemma}

\begin{proof} There exists an isomorphism of $FG_{i}$-$FG_{l}$-bimodules

\begin{center}
$(FG_{i}^{\rho_{1}}Ê\otimes_{FG_{i,j}} FG_{j}) \otimes_{FG_{j}} (FG_{j}^{\rho_{2}} \otimes_{FG_{j,l}} FG_{l}) \cong FG_{i}^{\rho_{1}} \otimes_{FG_{i,j}} FG_{j}^{\rho_{2}} \otimes_{FG_{j,l}} FG_{l}=M$
\end{center}
Let $T$ be a set of representatives of the cosets of $G_{i,j}G_{j,l}$ in $G_{j}$ and $Z$ a set of representatives of the cosets of $G_{i,j,l}$ in $G_{i,j}$. 
For $t \in T$, take $\nu_{t}=\displaystyle \sum_{zÊ\in Z} \rho_{1}(z^{-1}) \otimes zt \otimes 1 \in M$. Define
\begin{center}
$e_{\lambda}=\displaystyle \sum_{wÊ\in G_{i,j,l}} \lambda(w^{-1}) \otimes w \in FG_{i}Ê\otimes_{F} FG_{l}$
\end{center}
Let us show that the set $\{e_{\lambda} \ast \nu_{t}: t \in T\}$ generates the $FG_{i}$-$FG_{l}$-bimodule $M$. Indeed, the set of all elements $1Ê\otimes t \otimes 1$ is a set of generators for $M$ as an $FG_{i}$-$FG_{l}$-bimodule. We have
\begin{align*}
1Ê\otimes t \otimes 1 &=  \displaystyle \sum_{wÊ\in G_{i,j,l}} \rho_{1}(w^{-1}) \otimes wt \otimes 1 \\
&=\displaystyle \sum_{z \in Z,w \in G_{i,j,l}} \rho_{1}(z^{-1})\rho_{1}(w^{-1}) \otimes wzt \otimes 1 \\
&=\displaystyle \sum_{z \in Z} \rho_{1}(z^{-1}) \displaystyle \sum_{w \in G_{i,j,l}} \rho_{1} \rho_{2}(w^{-1}) \otimes \rho_{2}(w)zt \otimes 1 \\
&=\displaystyle \sum_{zÊ\in Z} \rho_{1}(z^{-1}) \displaystyle \sum_{w \in G_{i,j,l}} \rho_{1}\rho_{2}(w^{-1}) \otimes zt \otimes w \\
&=\displaystyle \sum_{w \in G_{i,j,l}} \rho_{1}\rho_{2}(w^{-1}) \left(\displaystyle \sum_{z \in Z} \rho_{1}(z^{-1}) \otimes zt \otimes 1\right)w \\
&=e_{\lambda} \ast \nu_{t}
\end{align*}
Therefore $M=\displaystyle \sum FG_{i}(e_{\lambda} \ast \nu_{t})FG_{l}$. Then
\begin{align*}
d_{i}(d_{j}/d_{i,j})(d_{l}/d_{j,l})&=\operatorname{dim}_{F}M \leq \displaystyle \sum_{t \in T} \operatorname{dim}_{F} (FG_{i}(e_{\lambda} \ast \nu_{t})FG_{l}) \\
& \leq \displaystyle \sum_{t \in T} \operatorname{dim}_{F}(FG_{i}^{\lambda} \otimes_{FG_{i,j,l}} FG_{l}) \\
&=[G_{j}: G_{i,j}G_{j,l}]d_{i}(d_{l}/d_{i,j,l}) \\
&=(d_{j}/d_{i,j}d_{j,l})d_{i,j,l}(d_{i}d_{l}/d_{i,j,l}) \\
&=d_{i}d_{j}d_{l}/d_{i,j}d_{j,l} \\
&=\operatorname{dim}_{F}M
\end{align*}
Consequently
\begin{center}
$\operatorname{dim}_{F}(FG_{i}^{\lambda} \otimes_{FG_{i,j,l}} FG_{l})=\operatorname{dim}_{F} FG_{i}e_{\lambda} \ast \nu_{t}FG_{l}$
\end{center}
this implies that $FG_{i}^{\lambda} \otimes_{FG_{i,j,l}} FG_{l} \cong FG_{i}e_{\lambda} \ast \nu_{t}FG_{l}$ and 
\begin{center}
$M \cong (FG_{i}^{\lambda} \otimes_{FG_{i,j,l}} FG_{l})^{[G_{j}: G_{i,j}G_{j,l}]}$
\end{center}
The proof is now complete.
\end{proof}
\begin{lemma} \label{lem9.9} For $\lambda \in \chi(G_{i,j,l})$ we have
\begin{center}
$FG_{i}^{\lambda} \otimes_{FG_{i,j,l}} FG_{l} \cong \displaystyle \bigoplus_{\rho \in \chi(G_{i,l})_{\lambda}} FG_{i}^{\rho} \otimes_{FG_{i,l}} FG_{l}$
\end{center}
where $\chi(G_{i,l})_{\lambda}$ is the set of characters of $G_{i,l}$ whose restriction to $G_{i,j,l}$ coincides with $\lambda$.
\end{lemma}
\begin{proof} Take $f_{\lambda}=\displaystyle \sum_{\rho \in \chi(G_{i,l})_{\lambda}} e_{\rho} \in FG_{i} \otimes_{F} FG_{l}$. Then for $g \in G_{i,j,l}$ we have
\begin{align*}
f_{\lambda}g&=\displaystyle \sum_{\rho \in \chi(G_{i,l})_{\lambda}} e_{\rho}g \\
&=\displaystyle \sum_{\rho \in \chi(G_{i,l})_{\lambda}} \underline{\rho}(g)e_{\rho}\\
&=\underline{g}f_{\lambda}
\end{align*}
Then we have an epimorphism
\begin{equation} \label{eq9.1}
FG_{i}^{\lambda} \otimes_{FG_{i,j,l}} FG_{l} \rightarrow FG_{i}f_{\lambda}FG_{l}
\end{equation}
Moreover,
\begin{align*}
\operatorname{dim}_{F}(FG_{i}f_{\lambda}FG_{l})&=\displaystyle \sum_{\rho \in \chi(G_{i,l})_{\lambda}} \operatorname{dim}_{F}(FG_{i}^{\rho} \otimes_{FG_{i,l}} FG_{l}) \\
&=(d_{i,l}/d_{i,j,l})d_{i}(d_{l}/d_{i,l}) \\
&=d_{i}d_{l}/d_{i,j,l} \\
&=\operatorname{dim}_{F} FG_{i}^{\lambda} \otimes_{FG_{i,j,l}} FG_{l}
\end{align*}
Hence the epimorphism given in (\ref{eq9.1}) is in fact an isomorphism. This proves our result.
\end{proof}
\begin{prop} \label{prop9.10} If $\rho_{1} \in \chi(G_{i,j})$ and $\rho_{2} \in \chi(G_{j,l})$, then
\begin{center}
$(FG_{i}^{\rho_{1}} \otimes_{FG_{i,j}} FG_{j}) \otimes_{FG_{j}} (FG_{j}^{\rho_{2}} \otimes_{FG_{j,l}} FG_{l}) \cong \displaystyle \bigoplus_{\rho \in \chi(G_{i,l})_{\lambda}} (FG_{i}^{\rho} \otimes_{FG_{i,l}} FG_{l})^{[G_{j}:G_{i,j}G_{j,l}]}$.
\end{center}
\begin{proof} This follows at once by combining Lemmas \ref{lem9.8} and \ref{lem9.9}.
\end{proof}
\end{prop} 
\begin{prop} \label{prop9.11}
For $\rho \in \chi(G_{i,j})$ we have
\begin{center}
$(FG_{i}^{\rho} \otimes_{FG_{i,j}} FG_{j})^{\ast} \cong FG_{j}^{\rho^{-1}} \otimes_{FG_{i,j}} FG_{i}$.
\end{center}
\end{prop}
\begin{proof} Let $e=g_{1},\ldots,g_{l}$ be representatives of the cosets of $G_{i,j}$ in $G_{i}$, then
\begin{center}
$FG_{i}^{\rho} \otimes_{FG_{i,j}} FG_{j}=g_{1} \otimes FG_{j} \oplus \ldots \oplus g_{l} \otimes FG_{j}$
\end{center}
as right $FG_{j}$-modules. Observe that $e=g_{1},\ldots,g_{l}$ is a $FG_{i}$-free basis for $FG_{i}^{\rho} \otimes_{FG_{i,j}} FG_{j}$. 
Thus we have the right $FG_{j}$-morphisms $g_{i}^{\ast}: FG_{i}^{\rho} \otimes_{FG_{i,j}} FG_{j} \rightarrow FG_{j}$, given by $g_{i}^{\ast}(g_{j})=\delta_{i,j}e$. The morphism $e^{\ast}$ generates $(FG_{i}^{\rho} \otimes_{FG_{i,j}} FG_{j})^{\ast}$ as an $FG_{j}$-$FG_{i}$- module. Moreover, for $x \in G_{i,j}$ and every $y \in FG_{i}^{\rho} \otimes_{FG_{i,j}} FG_{j}$, we have:
\begin{center}
$(e^{\ast}x)(y)=e^{\ast}(xy)=e^{\ast}(y\underline{\rho^{-1}}(x))=(\underline{\rho^{-1}}(x)e^{\ast})(y)$
\end{center}
therefore $e^{\ast}x=\underline{\rho^{-1}}(x)e^{\ast}$. Thus we get an epimorphism
\begin{equation} \label{eq9.2}
FG_{j}^{\rho^{-1}} \otimes_{FG_{i,j}} FG_{i} \rightarrow (FG_{i}^{\rho} \otimes_{FG_{i,j}} FG_{j})^{\ast}
\end{equation}
Since $\operatorname{dim}_{F}(FG_{j}^{\rho^{-1}} \otimes_{FG_{i,j}} FG_{i})=\operatorname{dim}_{F}(FG_{i}^{\rho} \otimes_{FG_{i,j}} FG_{j})^{\ast}$. It follows that the epimorphism given in (\ref{eq9.2}) is in fact an isomorphism. This completes the proof. 
\end{proof}

\begin{rem} \label{rem9.12} In Example \ref{ex9.4}, the modules $F_{i}^{\rho} \otimes_{F_{i,j}} F_{j}$ are indecomposable; in Example \ref{ex9.5}, the bimodules $FG_{i}^{\rho} \otimes_{FG_{i,j}} FG_{j}$ are, in general, not indecomposable. In Example \ref{ex9.4}, the field $F$ cannot be algebraically closed; on the other hand, in Example \ref{ex9.5}, the field $F$ can be algebraically closed, so in particular one can take $F=\mathbb{C}$. 
\end{rem}

Let $\mathcal{B}(i,j):=\operatorname{Gal}(F_{i,j}/F)$ and $\mathcal{B}=\displaystyle \bigcup_{i,j} \mathcal{B}(i,j)$. Let $g$ be a modulation of the quiver $Q$ in $\mathcal{B}$ given in Example \ref{ex9.4}. Consider the corresponding $E$-$E$-bimodule
\begin{center}
$M_{g}=\displaystyle \bigoplus_{\alpha \in Q_{1},\rho \in \mathcal{B}} M_{\rho}^{g(\alpha,\rho)}$
\end{center}
For $\rho \in \mathcal{B}(i,j)$ we have $M_{\rho}=F_{i}^{\rho} \otimes_{F_{i,j}} F_{j}$. Take $\sigma_{u}: M_{\rho} \rightarrow M_{\rho}^{g(\alpha,\rho)}$ the $u$-th canonical inclusion where $1 \leq u \leq g(\alpha,\rho)$. Set $a_{u}=\sigma_{u}(1 \otimes 1)$. Denote by $A_{\rho}$ the set $\{a_{u}\}_{1 \leq u \leq g(\alpha,\rho)}$ and $A=\displaystyle \bigcup_{\rho \in \mathcal{B}}A_{\rho}$, which is a local set of generators of the $E$-$E$-bimodule $M_{g}$. For $a \in A_{\rho}$, with $\rho \in \mathcal{B}(i,j)$, we have $EaE=F_{i}aF_{j}Ê\cong F_{i}^{\rho} \otimes_{F_{i,j}} F_{j}$. For $x \in F_{i,j}$ we have $xa=a\rho(x)$. Take $v$ a $d$-th primitive root of unity in $F$, $L_{i}=\{v^{(d/d_{i})m}\}_{0 \leq m \leq (d/d_{i})-1}$ and $L=\displaystyle \bigcup_{i=1}^{n} L_{i}$. The set $B_{i,j}=\{v^{(d_{i}/d_{i,j})m}\}_{0 \leq m \leq (d_{i}/d_{i,j})-1}$ forms a basis of $F_{i}$ over $F_{i,j}$. For $a \in A \cap e_{i}Me_{j}$, $Y_{a}=\{v^{(d/d_{j})m}a\}_{0 \leq m \leq (d_{j}/d_{i,j})-1}$ is an $F_{j}$-basis for $EaE$, and $X=\displaystyle \bigcup_{a \in A} X_{a}$ is a local basis for $_{E}M_{g}$. 
Take now an element in $\mathcal{F}_{E}(M_{g})_{cyc}$ of the form $P=\omega_{0}a_{1}\omega_{1}a_{2}\cdots \omega_{l-1}a_{l-1}\omega_{l}$ with $a_{i} \in e_{h(a_{i})}Ae_{t(a_{i})}$ and $\omega_{i} \in L(\sigma(a_{i+1}))$, $i=0,\ldots,l-1$, $\omega_{l} \in L(\tau(a_{l}))$. We have
\begin{equation} \label{eq9.3}
\begin{split}
\tilde{\delta}_{a^{\ast}}(P) &= \tilde{\delta}(a_{1}\omega_{1} \cdots a_{l}\omega_{l}\omega_{0}) \\
&= \displaystyle \sum_{s \in L} \displaystyle \sum_{k=1}^{l} a^{\ast}s^{-1}(a_{k}\omega_{k} \cdots a_{l} \omega_{l} \omega_{0} a_{1} \omega_{1} \cdots a_{k-1} \omega_{k-1})s \\
&=\displaystyle \sum_{s \in L_{i}} \displaystyle \sum_{k=1}^{l} a^{\ast}s^{-1}(a_{k}w_{k} \cdots a_{l}\omega_{l}\omega_{0}a_{1} \omega_{1} \cdots a_{k-1} \omega_{k-1})s\delta_{a,a_{k}}
\end{split}
\end{equation}
Then
\begin{center}
$s^{-1}a=v^{-r(d/d_{i})}v^{-(d/d_{i})(d_{i}/d_{i,j})m}a=v^{-r(d/d_{i})}a\rho^{-1}(v^{-(d/d_{i,j})m})$
\end{center}
with $0 \leq r \leq d_{i}/d_{i,j}$, $0 \leq m \leq d/d_{i,j}-1$. Now $a^{\ast}(s^{-1}a) \neq 0$ implies $r=0$; then the only non-zero elements in the sum of equation (\ref{eq9.3}) are the summands
\begin{center}
$\rho^{-1}(s^{-1})\omega_{k+1}a_{k+1} \cdots a_{l}\omega_{l}\omega_{0}a_{1} \omega_{1} \cdots a_{k-1} \omega_{k-1}s$
\end{center}
with $s=v^{(d/d_{i,j})m}$, $0 \leq m \leq (d/d_{i,j})-1$. Therefore
\begin{center}
$\tilde{\delta}_{a^{\ast}}(P)=\displaystyle \sum_{s \in B_{i,j}} \displaystyle \sum_{k=1}^{l} \rho^{-1}(s^{-1})a_{k+1}w_{k+1} \cdots a_{l}\omega_{l}\omega_{0}a_{1} \omega_{1} \cdots \omega_{k-1} a_{k-1}s\delta_{a,a_{k}}$.
\end{center}
From the preceding discussion we have the following
\begin{rem} \label{rem9.13} For any potential $P \in \mathcal{F}_{E}(M_{g})$ we have
\begin{center}
$\partial_{a}(P)=(1/d_{i,j})\tilde{\delta}_{a^{\ast}}(P)$
\end{center}
for $a \in A \cap e_{i}M_{g}e_{j}$, where $\partial_{a}$ is the cyclic derivative defined in \cite[p.19]{6}.
\end{rem}

Let $A_{\rho}$ be a set of generators of the $E$-$E$-bimodule $M_{\rho}^{g(\alpha,\rho)}$. For each $a \in e_{i}A_{\rho}e_{j}$, we set $\rho(a)=\rho$. Now consider Example \ref{ex9.4}. For any potential of the form $P=\omega_{0}a_{1}\omega_{1}a_{2}\omega_{2} \cdots \omega_{l-1}a_{l}\omega_{l}$, with $\omega_{i} \in L$ and $a_{j} \in A$, we have the same expression for
\begin{center}
 $\tilde{\delta }_{a^{*}}(P)=\displaystyle \sum _{s\in G_{i,j}}\displaystyle \sum _{k=1}^{l}\rho (a)^{-1}(s^{-1})a_{k+1}\omega _{k+1} \cdots a_{l}\omega _{l}\omega _{0}a_{1}\omega _{1} \cdots \omega _{k-1}a_{k-1})s\delta _{a,a_{k}}.$
\end{center}

\begin{section}{Mutations} \label{sec11}
Let $E=\displaystyle \prod _{i=1}^{n}E_{i}$ be a finite direct product of semisimple algebras, each $E_{i}$ having an $F$-basis $L_{i}$, satisfying conditions 1-2 and (\ref{eq6.1}) of page $21$. Let $e_{i}$ be the image of the unit of $E_{i}$ under the canonical $i$-th inclusion $E_{i} \hookrightarrow E$. We will assume that $E \otimes_{F} E^{op}$ is a semisimple $F$-algebra. In this section we consider $E$-$E$-bimodules $M$ such that $e_{i}Me_{j} \neq 0$, then this bimodule is an $E_{i}$-free left module and an $E_{j}$-free right module. We will present a definition for mutation of potentials in $\mathcal{F}_{E}(M)$. Then we will prove that, up to right-equivalence, the operation of mutation is an involution. 

For $i,j\in [1,n]$, $E_{i} \otimes _{F}E_{j}^{op}$ is a semisimple $F$-algebra. Then as in Section \ref{sec6}, we have a decomposition $1_{E_{i}}\otimes 1_{E_{j}^{op}}=\displaystyle \sum _{x\in U(i,j)}\epsilon _{x}$, into primitive orthogonal idempotents.
We also have a decomposition $1_{E_{i}}\otimes 1_{E_{j}^{op}}=\displaystyle \sum _{h\in V(i,j)}\xi _{h}$ into primitive orthogonal
central idempotents. If $M_{r}$ is a simple $E_{i}$-$E_{j}$-bimodule, then $M_{r}=E_{i}a_{r}E_{j}$, where $a_{r}=\epsilon _{x(r)}*a$ for some $x(r)\in U(i,j)$. Given $M_{r}$, there is a unique $\xi _{h}$ such that $0\neq \xi _{h}*M=M$. In this case we will say that $M_{r}$ belongs to $\xi _{h}$. It is known
that $M_{r}\cong M_{r_{1}}$ if and only if $M_{r}$ and $M_{r_{1}}$ belong to the same $\xi _{h}$. \\

Recall from Section \ref{sec6} there exists an antiisomorphism of algebras $\mathfrak{s}:E_{i}\otimes _{F}E_{j}^{op}\rightarrow E_{j}\otimes _{F}E_{i}^{op}$.

\begin{lemma} \label{lem10.1} Suppose $M_{r}$ is a simple $E_{i}$-$E_{j}$-bimodule and $M_{r_{1}}$ is a simple $E_{j}$-$E_{i}$-bimodule. Then $M_{r_{1}}\cong M_{r}^{*}$ if and only if $M_{r}$ belongs to $\xi _{h}$ and $M_{r_{1}}$ belongs to
$\mathfrak{s}(\xi _{h})$.
\end{lemma}

\begin{proof} It is enough to prove that $M_{r}^{*}$ belongs to $\mathfrak{s}(\xi _{h})$. Let $\{x_{1},\ldots,x_{l}\}$ be an
$F$-basis for $E_{i}$, then 
$$\xi _{h}=\sum _{z=1}^{l}x_{z}\otimes y_{z},$$
for some $y_{z}\in E_{j}$; here $\xi _{h}$ is in the center of $E_{i}\otimes _{F}E_{j}^{op}$. Then for all $a\in E_{i}$, we have
$$(1\otimes a)\xi _{h}=\sum _{z=1}^{l}x_{z}\otimes ay_{z}=\sum _{z=1}^{l}x_{z}\otimes y_{z}a=\xi _{h}(1\otimes a).$$
Therefore for all $a\in E_{j}$ and  $1\leq z \leq l$, $ay_{z}=y_{z}a$; thus, each $y_{z}$ is in the center of $E_{j}$. \\

Let $f\in M_{r}^{*}$, then for $m\in M_{r}$ we have
$$\mathfrak{s}(\xi _{h})(f)(m)=\displaystyle \sum_{z=1}^{l} (y_{z}fx_{z})(m)=\sum _{z=1}^{l}y_{z}f(x_{z}m)=\sum _{z=1}^{l}f(x_{z}m)y_{z}=$$
$$f(\sum _{z=1}^{l}x_{z}my_{z})=f(\xi _{h}*m)=f(m).$$
Then $\mathfrak{s}(\xi _{h})*M_{r}^{*}=M_{r}^{*}$. Therefore $M_{r}^{*}$ belongs to $\mathfrak{s}(\xi _{h})$. The lemma follows.
\end{proof}

\begin{rem} \label{rem10.2} Suppose $M_{r}=E_{i}aE_{j}$ is a simple $E_{i}$-$E_{j}$-bimodule and $M_{r_{1}}=E_{j}bE_{i}$ is a simple
$E_{j}$-$E_{i}$ bimodule. If $a=\epsilon _{x}*a$ and $b=\mathfrak{s}(\epsilon _{x})*b$, then
$M_{r_{1}}\cong M_{r}^{*}$. Indeed, suppose $M_{r}$ belongs to $\xi _{h}$, then
$$\mathfrak{s}(\xi _{h})*E_{j}bE_{i}=E_{j}\mathfrak{s}(\xi _{h})(\mathfrak{s}(\epsilon _{x}))*bE_{i}=$$
$$E_{j}\mathfrak{s}(\epsilon _{x}\xi _{h})*bE_{i}=E_{j}\mathfrak{s}(\epsilon _{x})*bE_{i}=M_{r_{1}}.$$
Therefore $M_{r_{1}}$ belongs to $\mathfrak{s}(\xi _{h})$. This implies that $M_{r_{1}}\cong M_{r}^{*}$.
\end{rem}  
 
In what follows, let $M$ be an $E$-$E$-bimodule of finite dimension over $F$. 
  
\begin{prop} \label{prop10.3} Suppose $M_{r}=E_{i}aE_{j}$ and $M_{r_{1}}=E_{j}bE_{i}$ are simple bimodules  which are direct summands of $M$, with
 $a=\epsilon _{x}*a$ and $b=\mathfrak{s}(\epsilon _{x})*b$. Then
 \begin{enumerate}[(a)]
\item If $P$ is any potential lying in $M_{r}\otimes _{E_{j}}M_{r_{1}}$, and not cyclically equivalent to zero, then there exists an automorphism $\varphi$ of the $E_{j}-E_{i}$-bimodule $M_{r_{1}}$ such that $(id_{M_{r}} \otimes \varphi)P$ is cyclically equivalent to $ab$. 
\item The potential $ab \in \mathcal{F}_{E}(M)$ is not cyclically equivalent to zero, and there exists $\vartheta \in M_{r_{1}}^{*}$ with  $\vartheta \tilde{\delta }(ab)\neq 0.$
 \item The morphism $\tilde{\delta }^{ab}$ induces an isomorphism
 $$\hat{\delta }^{ab}:(M_{r}\oplus M_{r_{1}})^{*}\rightarrow M_{r}\oplus M_{r_{1}}.$$
 \end{enumerate}
\end{prop}

\begin{proof}
\begin{enumerate}[(a)]
\item Here $P=\displaystyle \sum_{u} s_{u}at_{u}b\mu_{u}$, with $s_{u},\mu_{u} \in E_{i}, \mu_{u} \in E_{j}$. Therefore $P$ is cyclically equivalent to $a\left(\displaystyle \sum_{u} t_{u}b\mu_{u}s_{u}\right)=ab_{1}$.  Since $P$ is not cyclically equivalent to zero, then $b_{1}\neq 0$. Thus $P$ is cyclically equivalent to $ab_{1}=\epsilon _{x}*ab_{1}$, which is cyclically equivalent to $(\epsilon _{x}*a)\mathfrak{s}(\epsilon _{x})*b_{1}$. Since $\mathfrak{s}(\epsilon _{x})*b_{1}\neq 0$, there exists an automorphism $\varphi $ of $M_{r_{1}}$ such that $\varphi (\mathfrak{s}(\epsilon _{x})*b_{1})=b$. This proves (a).
\item Let $\{x_{u},x_{u}^{*}\}_{u\in T}$ be a finite dual basis of $(M_{r})_{E_{j}}$; here $M_{r_{1}}\cong M_{r}^{*}$. Let $\varrho :M_{r_{1}}\rightarrow M_{r}^{*}$ be an isomorphism of $E_{j}$-$E_{i}$-bimodules and $y_{u}\in M_{r_{1}}$ with $\varrho (y_{u})=x_{u}^{*}$. Thus if $s\in L(j)$, 
\begin{align*}
\left(\varrho \otimes id_{M_{r}}\right)\left(s \left(\displaystyle \sum_{u} y_{u} \otimes x_{u} \right)\right) &= s\left(\displaystyle \sum_{u} x_{u}^{\ast} \otimes x_{u}\right) \\
&=\displaystyle \sum_{u} (x_{u}^{\ast} \otimes x_{u})s \\
&=\left(\varrho \otimes id_{M_{r}}\right) \left(\displaystyle \sum_{u} y_{u} \otimes x_{u}\right)s \\
&=(\varrho \otimes id_{M_{r}}) \left( \left(\displaystyle \sum_{u} y_{u} \otimes x_{u}\right)s\right)
\end{align*}
Then if $Q=\displaystyle \sum _{u}y_{u}\otimes x_{u}$, we have
$$\tilde{\delta }(Q)=\sum _{s\in L(j)}s^{-1}Qs+\sum _{t\in L_{i}}t^{-1}\left(\sum _{u}x_{u}\otimes y_{u}\right)t=
c(j)Q+Q^{\prime },$$
where $c(j)=\mathrm{dim}_{F}E_{j}$, and $Q^{\prime }\in M_{r}\otimes M_{r_{1}}$. Since $c(j)\neq 0$, then
$\tilde{\delta }(Q)\neq 0$. Take $\vartheta \in M_{r_{1}}^{*}$ with $\vartheta Q\neq 0$, then $\vartheta \hat{\delta }(Q)\neq 0$. By (a), there exists an automorphism $\varphi $ of $ M_{r_{1}}$, such that $a\varphi (b)$ is cyclically equivalent to $Q$. Then $\tilde{\delta }(a\varphi (b))=\tilde{\delta }(Q)$. Therefore
$$0\neq \vartheta (\tilde{\delta }(a\varphi (b)))=\vartheta\left(\displaystyle \sum_{s \in L(j)} s^{-1} \varphi(b)as\right)=\vartheta\left(\displaystyle \sum_{s \in L(j)} \varphi(s^{-1}b)as\right)=\vartheta \varphi (\tilde{\delta }(ab)).$$
This proves (b).
\item By (b), $\tilde{\delta }^{ab}$ induces a non-zero morphism of $E_{i}$-$E_{j}$-bimodules from $M_{r_{1}}^{*}$ to $M_{r}$. Since $M_{r}$ is simple, then $\tilde{\delta }^{ab}$ induces an isomorphism from $M_{r_{1}}^{*}$ to $M_{r}$. Now $\epsilon _{x}=\mathfrak{s}(\mathfrak{s}(\epsilon _{x}))$, then applying (b) to $ba$, yields a map $\psi \in M_{r}^{*}$ such that $\tilde{\delta }^{ab}(\psi )=\psi (\tilde{\delta }(ab))\neq 0.$ Then 
$\tilde{\delta }^{ab}$ induces an isomorphism from $M_{r}^{*}$ to $M_{r_{1}}$. We conclude that
$\tilde{\delta }:(M_{r}\oplus M_{r_{1}})^{*}\rightarrow M_{r}\oplus M_{r_{1}}$ is an isomorphism. This proves (c).
\end{enumerate}
\end{proof}

Now let $L$ be an $E$-$E$-bimodule and $L=L_{1}\oplus L_{2}\oplus \ldots \oplus L_{l} $ where $L_{i}$ are $E$-$E$-subbimodules of $L$. We will identify
$L_{i}^{*}$ with those $h\in L^{*}$ such that $h(L_{j})=0$ for $j\neq i$. Under this identification,
$L^{*}=L_{1}^{*}\oplus \ldots \oplus L_{l}^{*}$.

\begin{theorem} \label{theo10.4} Suppose $P\in \mathcal{F}_{E}(M)$ is a potential. Then $P$ is right-equivalent to a potential $Q$, which is the direct sum of a trivial potential $W$ and $Q^{\geq 3}$.
\end{theorem}

\begin{proof} Let $\mathcal{A}$ be a polarization of $M$, we denote  $M_{\mathcal{A}}=\displaystyle \bigoplus _{(i,j)\in \mathcal{A}}e_{i}Me_{j}$,
$M_{\mathcal{A}^{\prime }}=\displaystyle \bigoplus _{(i,j)\in \mathcal{A}}e_{j}Me_{i}$.

We have $\tilde{\delta }^{P^{(2)}}:M_{\mathcal{A}}^{*}\rightarrow M_{\mathcal{A}^{\prime }}$. Let $\mathcal{L}=\mathrm{Ker}(\tilde{\delta }^{P^{(2)}})$, here $\tilde{\delta }^{P^{(2)}}$ is a morphism of $E$-$E$-bimodules, then $\mathcal{L}$ is an $E$-$E$-subbimodule of $M_{\mathcal{A}}^{*}$. Since we are assuming that $E\otimes _{F}E^{op}$ is semisimple, then there exists an $E$-$E$-subbimodule $\mathcal{L}_{1}$ of $M_{\mathcal{A}}^{*}$ such that
$$M_{\mathcal{A}}^{*}=\mathcal{L}\oplus \mathcal{L}_{1}.$$

Let $N_{\mathcal{A}}=\{m\in M_{\mathcal{A}}|h(m)=0 \textrm{ for all } h\in \mathcal{L}_{1}\}$, and $N^{\prime}$ be the subspace consisting of all elements $m\in M_{\mathcal{A}}$ such that $h(m)=0$ for all $h\in \mathcal{L}$. Then we have a decomposition of $E$-$E$-bimodules 
$$M_{\mathcal{A}}=N_{\mathcal{A}}\oplus N^{\prime }.$$

Moreover, $N_{\mathcal{A}}^{*}=\mathcal{L}$ and $(N^{\prime })^{*}=\mathcal{L}_{1}$. Now let $N_{\mathcal{A}^{\prime }}=\mathrm{Im}(\tilde{\delta }^{P})$, we have a decomposition into
$E$-$E$-bimodules
$$M_{\mathcal{A}^{\prime }}=N_{\mathcal{A}^{\prime }}\oplus M^{\prime }.$$

We now decompose $N_{\mathcal{A}}$ and $N^{\prime }$ into a direct sum of simple $E$-$E$-bimodules
$$N_{\mathcal{A}}=\displaystyle \bigoplus _{u\in I}E_{i(u)}a_{u}E_{j(i)}, \quad N^{\prime }=\displaystyle \bigoplus _{v\in J}E_{i(v)}a_{v}E_{j(v)},$$
and $a_{u}=\epsilon _{x(u)}*a_{u}, a_{v}=\epsilon _{x(v)}*a_{v}$.

The potential $P^{(2)}$ is cyclically equivalent to a potential
$$P^{\prime }=\sum _{u\in I}a_{u}b_{u}+\sum _{v\in J}a_{v}b_{v}.$$

We claim that for $v_{0}\in J$, $a_{v_{0}}b_{v_{0}}$ is cyclically equivalent to zero. Indeed, otherwise
$a_{v_{0}}b_{v_{0}}$ is cyclically equivalent to $(\epsilon _{x(v_{0})}*a_{v_{0}})(\mathfrak{s}(\epsilon _{x(v_{0})})*b_{v_{0}})$ with $\mathfrak{s}(\epsilon _{v_{0}(x)})*b_{v_{0}}\neq 0$. By (b) of Proposition \ref{prop10.3}, there exists
$\vartheta \in E_{i(v_{0})}a_{v_{0}}E_{j(v_{0})}^{*}$ such that $\tilde{\delta}^{a_{v_{0}}b_{v_{0}}}(\vartheta)=\vartheta \tilde{\delta }(a_{v_{0}}b_{v_{0}})\neq 0$. 

We can now extend $\vartheta $ to some  funtion $\underline{\vartheta }\in M_{\mathcal{A}}^{*}$ such that
$\underline{\vartheta }(N_{\mathcal{A}})=0$ and $\underline{\vartheta }(E_{i(v)}a_{v}E_{j(v)})=0$ for
$v\neq v_{0}$. Then $\underline{\vartheta }\in \mathcal{L}$ and
$$\tilde{\delta}^{P^{(2)}}(\underline{\vartheta})=\tilde{\delta}^{a_{v_{0}}b_{v_{0}}}(\underline{\vartheta}) \neq 0,$$
a contradiction. Therefore we may assume
$$P^{\prime }=\sum _{u\in I}a_{u}b_{u},$$
with $\epsilon _{x(u)}*a_{u}=a_{u}, \mathfrak{s}(\epsilon _{x(u)})*b_{u}=b_{u}$. 

We have 
$$N_{\mathcal{A}}^{*}=\displaystyle \bigoplus _{u\in I}\left(E_{i(u)}a_{u}E_{j(u)}\right)^{*},$$
 $$\tilde{\delta }^{P^{(2)}}\left((E_{i(u)}a_{u}E_{j(u)})^{*}\right)=\tilde{\delta }^{a_{u}b_{u}}((E_{i(u)}a_{u}E_{j(u)})^{*})=E_{j(u)}b_{u}E_{i(u)}.$$
 Then
 $$N_{\mathcal{A}^{\prime }}=\hat{\delta }^{P^{(2)}}(N_{\mathcal{A}}^{*})=\displaystyle \bigoplus _{u\in I}E_{j(u)}b_{u}E_{i(u)}.$$
 
 By (c) of Proposition \ref{prop10.3} we have
 $$\tilde{\delta }^{P^{2}}((N_{\mathcal{A}}\oplus N_{\mathcal{A}^{\prime }})^{*})=
\displaystyle \bigoplus _{u\in I}\tilde{\delta }^{a_{u}b_{u}}(\left(E_{i(u)}a_{u}E_{j(u)} \oplus E_{j(u)}bE_{I(u)}\right)^{*})=$$
$$\displaystyle \bigoplus _{u\in I}(E_{i(u)}a_{u}E_{j(u)}\oplus E_{j(u)}b_{u}E_{i(u)})=N_{\mathcal{A}}\oplus N_{\mathcal{A}^{\prime }}.$$

Therefore the potential $P^{\prime }=\displaystyle \sum _{u\in I}a_{u}b_{u}$ is trivial. Then $P$ is cyclically equivalent
to a potential $Q_{1}=\displaystyle \sum _{u\in I}a_{u}b_{u}+Q_{1}^{\geq 3}$. By Corollary \ref{coro5.6}, the potential
$\displaystyle \sum _{u\in I}a_{u}b_{u}$ has the splitting property. Therefore, 
there exists a decomposition of $E$-$E$-bimodules 
$$M=N_{\mathcal{A}}\oplus N_{\mathcal{A}^{\prime }}\oplus M^{\prime \prime }$$
and $Q_{1}$ is right-equivalent to a potential
$Q=\displaystyle \sum _{u\in I}a_{u}b_{u}+Q^{\geq 3}$ with $Q^{\geq 3}\in \mathcal{F}_{E}(M^{\prime \prime })$.
This completes the proof.
\end{proof}

\begin{definition} \label{def10.5} Let $P\in \mathcal{F}(M)$ be a potential and $k$ an integer in $\{1,\ldots,n\}$. Suppose that there are no two-cycles passing through $k$. Then $\mu _{k}(P)$, the pre-mutation of $P$ in the direction $k$, is defined. By Theorem \ref{theo10.4}, $\mu _{k}(P)$ is right-equivalent to the direct sum of a trivial potential $W$ and a potential $Q$ in $\mathcal{F}_{E}(M)^{\geq 3}$. Following \cite{4}, we define the mutation of $P$ in the direction $k$, as $\overline{\mu }_{k}(P)=Q$.
\end{definition}
In order to see that $\overline{\mu }_{k}(P)$ is well-defined up to right-equivalence we need to prove that Proposition \ref{prop4.5} holds in our case.

\begin{prop} \label{prop10.6} Suppose $P$ and $P^{\prime }$ are potentials in $\mathcal{F}_{E}(M)$ such that $P^{\prime }-P\in J(P)^{2}$. Then there exists an algebra automorphism $\varphi $ of $\mathcal{F}_{E}(M)$, with $\varphi |_{E}=id_{E}$, such that $\varphi (P)$ is cyclically equivalent to $P^{\prime }$. Moreover, $\varphi (f)-f\in J(P)$ for all $f\in \mathcal{F}_{E}(M)$.
\end{prop}

\begin{proof} We have $M=\displaystyle \bigoplus _{u\in  U}M_{u}$, where each $M_{u}$ is a simple $E$-$E$-bimodule and $M_{u}= E_{i(u)}a_{u}E_{j(u)}$, with $a_{u}=\epsilon _{x(u)}*a_{u}$. For each $uÊ\in U$, choose an embedding of $E_{i(u)}-E_{j(u)}$-bimodules $M_{u} \rightarrow \hat{M}_{u}=E_{i(u)}\hat{a}_{u}E_{j(u)}$, such that $a_{u}=\epsilon _{x(u)}*\hat{a}_{u}$, and there is an isomorphism from $\hat{M}_{u}$ into $E_{i(u)} \otimes_{F} E_{j(u)}$, sending $\hat{a}_{u}$ into $1_{E_{i(u)}} \otimes 1_{E_{j(u)}}$. Then $M$ is an $E$-$E$-subbimodule of $\hat{M}=\displaystyle \bigoplus_{u \in U} \hat{M}_{u}$ and we have an inclusion $\mathcal{F}_{E}(M)\subseteq \mathcal{F}_{E}(\hat{M})$, where $\hat{M}$ is a $Z$-free $E$-$E$ bimodule. Using the notation introduced in Section \ref{sec4}, we have $J(P)=R(P)\cap \mathcal{F}_{E}(M)$.

Here $P^{\prime }-P\in J(P)^{2}\subseteq R(P)^{2}$, then by Proposition \ref{prop4.5}, there exists an algebra automorphism
$\hat{\varphi }$ of $\mathcal{F}_{E}(\hat{M})$ such that $\hat{\varphi }(P)$ is cyclically equivalent to $P^{\prime }$.

We will show that $\hat{\varphi }$ can be chosen in such a way that it restricts to an automorphism $\varphi $ of $\mathcal{F}_{E}(M)$ and $\varphi (P)-P^{\prime }\in [\mathcal{F}_{E}(M),\mathcal{F}_{E}(M)]$. 

For this, we will prove that the morphisms: $\varphi _{l}^{(2)}:\hat{M}\rightarrow \mathcal{F}_{E}(\hat{M})$ of the Claim in the proof of Proposition \ref{prop4.5}, can be chosen in such a way that $\varphi _{l}^{(2)}(\hat{a}_{u})\in \mathcal{F}_{E}(M)$ for all $u\in U$.  In this case, taking $\varphi_{l}$ the automorphism of $\mathcal{F}_{E}(\hat{M})$ determined by the pair $(id_{\hat{M}}, \varphi_{l}^{2})$, we have for $u\in U$
$$\varphi_{l}^{(2)}(a_{u})=\varphi_{l}(\epsilon _{x(u)}*\hat{a}_{u})=\epsilon _{x(u)}*\varphi_{l}(\hat{a}_{u})=\epsilon _{x(u)}*(\hat{a}_{u}+\varphi_{l} ^{(2)}(\hat{a}_{u}))=a_{u}+\epsilon _{x(u)}*\varphi_{l}^{(2)}(\hat{a}_{u})\in \mathcal{F}_{E}(M).$$

Therefore $\varphi_{l}$ restricts to an automorphism of $\mathcal{F}_{E}(M)$. We will prove by induction on $l$ that the morphisms $\varphi _{l}^{(2)}$ can be chosen such that $\varphi _{l}^{(2)}(a_{u})\in \mathcal{F}_{E}(M)$ for all $u\in U$.

Since $E\otimes _{F}E^{op}$ is semisimple, then $\hat{M}=M\oplus L$ for some $E$-$E$-subbimodule $L$ of $\hat{M}$. Let $\mathcal{L}$ denote the closure of the two-sided ideal generated by $L$ in $\mathcal{F}_{E}(\hat{M})$.
We have 
$$\mathcal{F}_{E}(\hat{M})=\mathcal{F}_{E}(M)\oplus \mathcal{L}.$$
If $w\in [\mathcal{F}_{E}(\hat{M}),\mathcal{F}_{E}(\hat{M})]$, then $w=w_{1}+w_{2}$ with $w_{1}\in [\mathcal{F}_{E}(M),\mathcal{F}_{E}(M)]$ and $w_{2}\in \mathcal{L}$. Moreover, note that for $u\in U$, $X_{\hat{a}_{u}^{*}}(P)\in \mathcal{F}_{E}(M)$; then
$X_{\hat{a}_{u}^{*}}(P)\in J(P)=R(P)\cap \mathcal{F}_{E}(M)$. Now if $\psi \in M^{*}$, let $\hat{\psi }\in \hat{M}^{*}$ be the extension of $\psi$ with $\hat{\psi }(L)=0$. Then
$$\hat{\delta }^{P}(\psi )=\hat{\delta }^{P}(\hat{\psi }).$$

We have $\hat{\psi }=\displaystyle \sum _{u\in U} \displaystyle \sum _{i\in Z(u)}\mu _{i}\hat{a}_{u}\nu _{i}$, with
$\mu _{i}, \nu _{i}\in E$. From this we obtain
$$\hat{\delta }^{P}(\psi )=\sum _{u\in U}\sum _{i\in Z(u)}\mu _{i}\hat{\delta }^{P}(\hat{a}_{u})\nu _{i}.$$

Therefore, $J(P)$ is the closure of the two-sided ideal in $\mathcal{F}_{E}(M)$ generated by the elements $X_{\hat{a}_{u}^{*}}=\hat{\delta }^{P}(\hat{a}_{u})$. 

We first choose $\varphi _{1}^{(2)}$. Applying Lemma \ref{lem4.4} to $\mathcal{F}_{E}(M)$ we have $$P^{\prime }-P=\sum _{u\in U}f(\hat{a}_{u})X_{\hat{a}_{u}^{*}}+w$$ with $w\in [\mathcal{F}_{E}(M),\mathcal{F}_{E}(M)]$ and $f(\hat{a}_{u})\in J(P)$. Since $P\in \mathcal{F}_{E}(M)^{\geq 3}$, then $J(P)\subseteq \mathcal{F}_{E}(M)^{\geq 2}$. Therefore $f(\hat{a}_{u})\in \mathcal{F}_{E}(M)^{\geq 2}\cap J(P)$.
Now choose $\varphi _{1}^{(2)}$ with $\varphi ^{(2)}(\hat{a}_{u})=f(\hat{a}_{u})$. Suppose we have already chosen $\phi_{1}^{(2)},\ldots,\phi _{l}^{(2)}$ with $\phi _{t}^{(2)}(\hat{a}_{u})\in \mathcal{F}_{E}(M)$ for all $t=1,\ldots,l$ and $u\in U$ satisfying conditions (i), (ii) in the Claim of the proof of Proposition \ref{prop4.5}.

We have
$$\varphi _{l}(P)-P-\sum _{u\in U}\varphi _{l}^{(2)}(\hat{a}_{u})X_{\hat{a}_{u}^{*}}(P)+w\in (\mathcal{F}_{E}(\hat{M})^{\leq l+2}\cap R(P))R(P)$$
with $w\in [\mathcal{F}_{E}(\hat{M}),\mathcal{F}_{E}(\hat{M})]$, and $R(\varphi _{l}(P))=R(P)$. Here $R(P)\subseteq J(P)+\mathcal{L}$, and $R(P)\cap \mathcal{F}_{E}(\hat{M})^{\geq l+2}\subseteq J(P)\cap \mathcal{F}_{E}(M)^{\geq l+2}+\mathcal{L}$. Moreover, $w=w_{1}+w_{2}$ with $w_{1}\in [\mathcal{F}_{E}(M),\mathcal{F}_{E}(M)]$ and $w_{2}\in \mathcal{L}$. Therefore
$$\varphi _{l}(P)-P-\sum _{u\in U}\varphi _{l}^{(2)}(\hat{a}_{u})X_{\hat{a}_{u}^{*}}(P)+w_{1}\in (J(P)\cap \mathcal{F}_{E}(M))J(P).$$
Also,
$$ \varphi _{l}(J(P))=J(\varphi _{l}(P))=R(\varphi _{l}(P))\cap \mathcal{F}_{E}(M)=R(P)\cap \mathcal{F}_{E}(M)=J(P).$$
As in the proof of Proposition \ref{prop4.5}, there exists a $z\in (\mathcal{F}_{E}(M)^{\geq l+2}\cap J(P))J(P)$, such that
$$\varphi _{l}(P)-P-\sum _{u\in U}\varphi _{l}^{(2)}(\hat{a}_{u})X_{\hat{a}_{u}^{*}}(P)+w_{1}=\varphi _{l}(z).$$

Applying Lemma \ref{lem4.4} to $\mathcal{F}_{E}(M)$ yields
$$z=\sum _{u\in U}-h(\hat{a}_{u})X_{\hat{a}_{u}^{*}}(P)+w_{2}$$
with $w_{2}\in [\mathcal{F}_{E}(M),\mathcal{F}_{E}(M)]$ and $h(\hat{a}_{u})\in \mathcal{F}_{E}(M)^{\geq l+2}J(P)$. Then we define $\varphi _{l+1}^{2}:\hat{M} \rightarrow \mathcal{F}_{E}(\hat{M})^{\geq 2}$ such that $\varphi _{l+1}(\hat{a}_{u})=h(\hat{a}_{u})$. Now, by Proposition \ref{prop1.10}, there exists an algebra automorphism $\varrho $ of $\mathcal{F}_{E}(\hat{M})$ such that
\begin{enumerate}
\item[(a)] $P^{\prime }-\varrho (P)\in [\mathcal{F}_{E}(\hat{M}),\mathcal{F}_{E}(\hat{M})]$.
\item [(b)] For $f\in \mathcal{F}_{E}(\hat{M})$, $\varrho (f)=\displaystyle \lim_{l \to \infty } \varrho _{l} \cdots \varrho _{0}(P)$, where $\varrho _{l}=\varphi ^{-1}_{l}$.
\item[(c)] For $m\in\hat{ M}$, $\varrho (m)-m\in \mathcal{F}_{E}(M)^{\geq 2}$.
\end{enumerate}

Since $\varphi _{l}$ restricts to an automorphism of $\mathcal{F}_{E}(M)$, then $\varphi _{l}^{-1}$ restricts to an automorphism of $\mathcal{F}_{E}(M)$. Using (b), we obtain the inclusion $\varrho (\mathcal{F}_{E}(M))\subseteq \mathcal{F}_{E}(M)$. Moreover, by (c), $\varrho (m)-m\in \mathcal{F}_{E}(M)^{\geq 2}$ for $m\in \hat{M}$. Therefore $\varrho $ restricts to an automorphism of $\mathcal{F}_{E}(M)$. Finally, using (a) we get
$$P^{\prime }-\varrho (P)\in [\mathcal{F}_{E}(\hat{M}), \mathcal{F}_{E}(\hat{M})]\cap \mathcal{F}_{E}(M)=[\mathcal{F}_{E}(M), \mathcal{F}_{E}(M)].$$
This completes the proof.
\end{proof}

\begin{definition} \label{def10.7} For $P$ a potential in $\mathcal{F}_{E}(M)$ we define $\Theta (P)=\tilde{\delta }^{P^{(2)}}(M^{*})\subseteq M$.
\end{definition}

\begin{prop} \label{prop10.8} Let $P$ be a potential in $\mathcal{F}_{E}(M)$ and $\varphi: \mathcal{F}_{E}(M)\rightarrow \mathcal{F}_{E}(M_{1})$ be an algebra isomorphism with $\varphi |_{E}=id_{E}$, determined by the pair $(\varphi ^{(1)},\varphi ^{(2)})$. Then $ \Theta (\varphi (P))=\varphi ^{(1)}\left(\Theta (P)\right).$
\end{prop}

\begin{proof} Since any isomorphism $\varphi$ with $\varphi |_{E}=id_{E}$ is the composition of an automorphism $\psi$ with $\psi (M)=M$ and an uni-triangular automorphism, it is enough to prove our claim for the case in which $\varphi $ is uni-triangular and when $\varphi (M)=M$. In the first case, $\varphi (P)=P^{(2)}+\varphi (P)^{\geq 3}$; then $\varphi (P)^{(2)}=P^{(2)}$ and our result follows. In the second case, 
\begin{align*}
&\Theta (\varphi (P ))=\Theta (\varphi (P)^{(2)})=\Theta (\varphi (P^{(2)}))=J(\varphi (P^{(2)}))\cap M_{1} \\
&=\varphi (J(P^{(2)}))\cap \varphi (M)=\varphi (J(P^{(2)})\cap M)=\varphi (\Theta (P))
\end{align*}
This proves our result.
\end{proof}

\begin{prop} \label{prop10.9} Suppose $P+W\in \mathcal{F}_{E}(M\oplus N)$ and $P^{\prime }+W^{\prime }\in \mathcal{F}_{E}(M_{1}\oplus N_{1})$ are potentials  with $P\in \mathcal{F}_{E}(M)^{\geq 3}$,
$P^{\prime }\in \mathcal{F}_{E}(M_{1})^{\geq 3}$, and  $W\in \mathcal{F}_{E}(N)$, $ W^{\prime }\in \mathcal{F}_{E}(N_{1})$ are trivial potentials. Then $P+W$ right-equivalent to $P^{\prime }+W^{\prime }$ implies that $P$ is right-equivalent to $P^{\prime}$. 
\end{prop}

\begin{proof} Let $\mathcal{N}$ be the closure of the two-sided ideal in $\mathcal{F}_{E}(M)$ generated by $N$. Similarly, let $\mathcal{N}_{1}$, be the closure of the two-sided ideal generated by $N_{1}$.
We have
$$\mathcal{F}_{E}(M\oplus N)=\mathcal{F}_{E}(M)\oplus \mathcal{N},\quad \mathcal{F}_{E}(M_{1}\oplus N_{1})=\mathcal{F}_{E}(M_{1})\oplus \mathcal{N}_{1}.$$

Suppose there exists an isomorphism $\varphi :\mathcal{F}_{E}(M\oplus N)\rightarrow \mathcal{F}_{E}(M_{1}\oplus N_{1})$ such that $\varphi |_{E}=id_{E}$. Then $\varphi $ is determined by a pair of morphisms of $E$-$E$-bimodules, $(\varphi ^{(1)},\varphi ^{(2)})$. By  Proposition \ref{prop10.8}, we have 
$$N_{1}=\Theta (P^{\prime}+W^{\prime })=\Theta \varphi (P+W)=\varphi ^{(1)}\Theta (P+W)=\varphi ^{(1)}(N).$$

Then 
$$\varphi ^{(1)}=\left( \begin{array}{cc}\varphi _{1}^{(1)}&0\\ \varphi ^{(1)}_{3}& \varphi ^{(1)}_{2}\end{array}\right):
M\oplus N \rightarrow M_{1} \oplus N_{1}.$$
Thus $\varphi ^{(1)}_{1}: M\rightarrow M_{1}$ and $\varphi ^{(1)}_{2}: N\rightarrow N_{1}$ are isomorphisms.

Let $\psi =p\varphi \iota: \mathcal{F}_{E}(M)\rightarrow \mathcal{F}_{E}(M_{1})$ be the morphism $\psi $ determined by the pair $(\varphi _{1}^{(1)}, p\varphi ^{(2)}\iota )$; hence $\psi $ is an automorphism. Moreover,

$$\varphi =\left( \begin{array}{cc}\psi & \varphi _{2}\\ \varphi _{3}&\varphi _{4}\end{array}\right) :\mathcal{F}_{E}(M)\oplus \mathcal{N}\rightarrow \mathcal{F}_{E}(M_{1})\oplus \mathcal{N}_{1}.$$

We have
$$\varphi (J(P)+\mathcal{N})=\varphi J(P+W)=J(\varphi (P+W))=J(P^{\prime }+W_{1})=J(P^{\prime })+\mathcal{N}_{1}.$$
Then $\varphi (P+W)=\psi (P)+\varphi _{2}(W)+\varphi _{3}(P)+\varphi _{4}(W)$ with
$\psi (P)+\varphi _{2}(W)\in \mathcal{F}_{E}(M_{1})$ and $\varphi _{3}(P)+\varphi _{4}(W)\in \mathcal{N}_{1}.$ Similarly, for $n\in N$, $\varphi (n)=\varphi _{2}(n)+\varphi _{4}(n)$ with $\varphi _{2}(n)\in J(P^{\prime })$ and $\varphi _{4}(n)\in \mathcal{N}_{1}$. Then $\phi _{2}(W)\in J(P^{\prime })^{2}$. Consequently,
$P^{\prime }-\psi (P)=\varphi _{2}(W)\in J(P^{\prime })^{2}$. Then $\psi ^{-1}(P^{\prime })-P\in J(\psi ^{-1}(P^{\prime }))^{2}$. By Proposition \ref{prop10.6}, there exists an automorphism $\tau $ of $\mathcal{F}_{E}(M)$ with, $\tau |_{E}=id_{E}$, such that $\tau \psi ^{-1}(P^{\prime })$ is cyclically equivalent to $P$. Therefore $\psi \tau ^{-1}(P^{\prime })$ is cyclically equivalent to $P$. The proof is now complete. 
\end{proof}

\begin{prop} \label{prop10.10} If $P\in \mathcal{F}_{E}(M)$ is right-equivalent to $P^{\prime }\in \mathcal{F}_{E}(M_{1})$, then $\overline{\mu }_{k}P$ is right-equivalent to $\overline{\mu }_{k}P^{\prime }$.
\end{prop}

\begin{proof} Let $\varphi :\mathcal{F}_{E}(M)\rightarrow \mathcal{F}_{E}(M_{1})$ be an algebra isomorphism such that $\varphi (P)$ is cyclically equivalent to $P^{\prime }$. We have that $\varphi (\kappa P)$ is cyclically equivalent to $\varphi (P)$. Then by Proposition \ref{prop7.12}, $\mu _{k}\varphi (\kappa P)$ is cyclically equivalent to $\mu _{k}P^{\prime }$. Since $\overline{e}_{k}\kappa P\overline{e}_{k}=0$, by Theorem \ref{theo7.9}, we have that $\varphi (\mu _{k}P)=\varphi (\mu _{k}\kappa P)$ is cyclically equivalent to $\mu _{k}\varphi (\kappa P)$. Therefore $\varphi (\mu _{k}P)$ is cyclically equivalent to $\mu _{k}P^{\prime }$. It follows that $\mu _{k}P$ is right-equivalent to $\mu _{k}P^{\prime }$. Applying Proposition \ref{prop10.9}, yields $\overline{\mu }_{k}P$ is right-equivalent to $\overline{\mu }_{k}P^{\prime }$. This completes the proof.
\end{proof}

\begin{theorem} \label{theo10.11} Let $M$ be a $2$-acyclic $E$-$E$-bimodule, such that if $e_{i}Me_{j}\neq 0$, then $e_{i}Me_{j}$ is an $E_{i}$-free left module and an $E_{j}$-right free module.  Suppose $P$ is a potential in $\mathcal{F}_{E}(M)$ and $k$ an integer in $[1,n]$. Then $\overline{\mu }_{k}^{2}(P)$ is right-equivalent to $P$.
\end{theorem}

\begin{proof}
In order to prove our theorem we first prove the following.
\begin{prop} \label{prop10.12} Let $M$ be a $2$-acyclic $E$-$E$-bimodule satisfying the conditions of Theorem \ref{theo10.11}. Let $P$ be a  potential in $\mathcal{F}_{E}(M)$, then $\mu _{k}^{2}(P)$ is right-equivalent to the direct sum of $P$ and a trivial quadratic potential.
\end{prop}

\begin{proof} Let $M\rightarrow \underline{M}$ be an inclusion, where $\underline {M}$ is a $Z$-free bimodule.  We have $\underline{M}=M\oplus L$ for some $E$-$E$-bimodule $L$. Then $\underline{M}e_{k}\underline{M}=Me_{k}M\oplus H_{1}$ for $H_{1}=Me_{k}L\oplus Le_{k}M\oplus Le_{k}L$. Moreover, $^{*}(\underline{M}e_{k})=^{*}(Me_{k})\oplus ^{*}(Le_{k})$ and $(e_{k}\underline{M})^{*}=(e_{k}M)^{*}\oplus (e_{k}L)^{*}$. Then
$$(e_{k}\underline{M})^{*}e_{k}(^{*}(\underline{M}e_{k}))=(e_{k}M)^{*}e_{k}(^{*}(Me_{k}))\oplus H_{2},$$
where $H_{2}=(e_{k}L)^{*}e_{k}(^{*}(Me_{k}))\oplus (e_{k}M)^{*}e_{k}(^{*}L)\oplus L^{*}e_{k}(^{*}L)$. We set $\nu _{k}(\underline{M})=\underline{M}\oplus \underline{M}e_{k}\underline{M}\oplus (e_{k}\underline{M})^{*}e_{k}(^{*}(\underline{M}e_{k}))$ and $\nu _{k}(M)=M\oplus Me_{k}M\oplus (e_{k}M)^{*}e_{k}(^{*}(Me_{k}))$.

We have:
$$\nu _{k}(\underline{M})=\nu _{k}(M)\oplus H$$
where  $H=L\oplus H_{1}\oplus H_{2}$, and 
$$\mathcal{F}_{E}(\nu _{k}(\hat{M}))=\mathcal{F}_{E}(\nu _{k}(M))\oplus \mathcal{H},$$
where $\mathcal{H}$ is the closure of the two-sided ideal generated by $H$. 

Now for $i,j\in [1,n]$, $e_{i}\underline{M}e_{j}=e_{i}Me_{j}\oplus e_{i}Le_{j}$. If $e_{i}Le_{j}\neq 0$, then  $e_{i}\underline{M}e_{j}$ is an $E_{i}$-free left module and an $E_{j}$-right module. If $e_{i}Me_{j}\neq 0$, then $e_{i}Me_{j}$ is an $E_{i}$-free left and an $E_{j}$-free right module; the same holds for $e_{i}Le_{j}$.

In case $e_{i}Me_{j} \neq 0$, let $X_{i,j}$ be a free basis for $e_{i}Me_{j}$; otherwise, take $X_{i,j}=\emptyset$. If $e_{i}Le_{j} \neq 0$, then let $X^{\prime }_{i,j}$ be an $E_{i}$-free basis for $e_{i}Le_{j}$; otherwise, let $X^{\prime }_{i,j}=\emptyset $. Now let $\hat{X}_{i,j}=X_{i,j}\cup X_{i,j}^{\prime }$. Similarly, we define $Y_{i,j}$ and $Y^{\prime }_{i,j}$ right $E_{j}$-free bases for $e_{i}Me_{j}\neq 0 $ and $e_{i}Le_{j}\neq 0$, respectively. We take $\hat{Y}_{i,j}=Y_{i,j}\cup Y^{\prime }_{i,j}$, $X=\cup _{i,j}X_{i,j}$, $\hat{X}=\cup _{i,j}\hat{X}_{i,j}$, $X^{\prime }=\cup _{i,j}X_{i,j}^{\prime }$. Similar definitions for $Y$, $Y^{\prime }$ and $\hat{Y}$. If $Z$ is one of the above sets, we let $_{k}Z=\cup _{j}Z_{k,j}$, $Z_{k}=\cup _{i}Z(i,k)$.

We recall from Section \ref{sec8} the morphisms 
$$i_{M}:\mathcal{F}_{E}(M)\rightarrow \mathcal{F}_{E}(\widehat{M}), j_{\mu _{k}M}:\mathcal{F}_{E}(\mu _{k}M)\rightarrow \mathcal{F}_{E}(\widehat{\mu _{k}M}),$$
$$j_{\mu _{k}^{2}M}:\mathcal{F}_{E}(\mu _{k}^{2}M)\rightarrow \mathcal{F}_{E}(\widehat{\mu _{k}M})$$
and the corresponding maps $i_{\underline{M}}, j_{\mu _{k}\underline{M}}, j_{\mu _{k}^{2}\underline{M}}$.

We have inclusions $\sigma _{1}:\mathcal{F}_{E}(\mu _{k}M)\rightarrow \mathcal{F}_{E}(\mu _{k}\underline{M})$,
$\sigma _{2}:\mathcal{F}_{E}(\mu _{k}^{2}M)\rightarrow \mathcal{F}_{E}(\mu _{k}^{2}\underline{M})$ and
$\sigma _{3}:\mathcal{F}_{E}(\widehat{\mu _{k}M})\rightarrow \mathcal{F}_{E}(\widehat{\mu _{k}\underline{M}}),$
and the relations $\sigma_{3}j_{\mu _{k}^{2}M}=j_{\mu _{k}^{2}\underline{M}}\sigma _{2}$, $\sigma _{3}j_{\mu _{k}M}=j_{\mu _{k}\underline{M}}\sigma _{1}$.
Consider the following  morphisms defined before Theorem \ref{theo8.1}:
$$\tau : \overline{e}_{k}\mathcal{F}_{E}(\mu _{k}M)\overline{e}_{k}\rightarrow \mathcal{F}_{E}(\mu _{k}^{2}M),$$  
$$\underline{\tau }: \overline{e}_{k}\mathcal{F}_{E}(\mu _{k}\underline{M})\overline{e}_{k}\rightarrow \mathcal{F}_{E}(\mu _{k}^{2}\underline{M}),$$
they satisfy the following equalities: $j_{\mu _{k}^{2}M}\tau =j_{\mu _{k}M}, j_{\mu _{k}^{2}\underline{M}}(\underline{\tau })=j_{\mu _{k}\underline{M}}.$

We claim that the restriction of $\underline{\tau }$ coincides with $\tau $, that is $\underline{\tau }\sigma _{1}=\sigma _{2}\tau $. Indeed, we have the equalities
$$ j_{\mu _{k}^{2}\underline{M}}\sigma _{2}\tau =\sigma _{3}j_{\mu _{k}^{2}M}\tau =\sigma _{3}j_{\mu _{k}M}=j_{\mu _{k}\underline{M}}\sigma _{1},$$
$$i_{\mu _{k}^{2}\underline{M}}\underline{\tau }\sigma _{1}=j_{\mu _{k}\underline{M}}\sigma _{1}.$$
and therefore $\sigma _{2}\tau =\underline{\tau }\sigma _{1}$. Now, consider the automorphism 
$\underline{\varphi }$ of $\mathcal{F}_{E}(\nu _{k}\underline{M})$ such that is the identity on $\overline{e}_{k}\underline{M}\overline{e}_{k}$, $e_{k}\underline{M}$ , $^{*}(e_{k}\underline{M})\otimes e_{k}(^{*}(\underline{M}e_{k}))$; $-id _{\underline{M}e_{k}}$ on $\underline{M}e_{k}$; and a similar automorphism $\varphi $ on $\mathcal{F}_{E}(\nu _{k}M)$.  Clearly, the restriction of $\underline{\varphi }$ to $\mathcal{F}_{E}(\nu _{k}M)$ coincides with $\varphi $. Now define an automorphism $\underline{\phi }$ of $\mathcal{F}_{E}(\nu _{k}\underline{M})$ as the identity map on $\underline{M}$ and $(e_{k}\underline{M})^{*}e_{k}(^{*}(\underline{M}e_{k}))$; and for $x\in \underline{M}e_{k}, x_{1}\in e_{k}\underline{M}$, $\underline{\phi }([xx_{1}])=[xx_{1}]+xx_{1}$. A similar definition for the automorphism $\phi $ of $\mathcal{F}_{E}(\nu _{k}M)$. It is clear that the restriction of $\underline{\phi }$ to $\mathcal{F}_{E}(\nu _{k}M)$ coincides with $\phi$.

Let $P$ be a potential in $\mathcal{F}_{E}(M)$. In what follows, we denote by $\underline{\mu }_{k}^{2}(P)$ the square of the premutation $\underline{\mu}_{k}$ in $\mathcal{F}_{E}(\underline{M})$, and by $\mu _{k}^{2}(P)$ the corresponding square in $\mathcal{F}_{E}(M)$. By Theorem \ref{theo8.1}, there exist isomorphisms of $E$-$E$-bimodules $\mathfrak{f}^{1}:\mu _{k}^{2}M\rightarrow \nu _{k}M$ and $\underline{\mathfrak{f}}^{1}:\mu _{k}^{2}\underline{M}\rightarrow \nu _{k}\underline{M}$. It is easy to see that the restriction of $\underline{\mathfrak{f}}^{1}$ to $\mu _{k}^{2}M$ coincides with $\mathfrak{f}^{1}$. The above isomorphisms induce isomorphisms of topological algebras: $\underline{\mathfrak{f}}:\mathcal{F}_{E}(\mu _{k}^{2}\underline{M})\rightarrow \mathcal{F}_{E}(\nu _{k}\underline{M})$, $\mathfrak{g}:\mathcal{F}_{E}(\mu _{k}^{2}M)\rightarrow \mathcal{F}_{E}(\nu _{k}M)$. Clearly, the second morphism is the restriction of the first one.

Denote by $\mathcal{W}$ the complement of $(X_{k}\otimes _{k}Y)$ in $\hat{X}_{k}\times _{k}\hat{Y}$. We have (see proof of Theorem \ref{theo8.1})
 \begin{align*}
 &\underline{\varphi }\underline{\mathfrak{f}}(\underline{\mu }_{k}^{2}(P))=\underline{\mathfrak{f}}(\underline{\tau }([P]))+\sum _{x\in \hat{X}_{k},y\in _{k}\hat{Y}}([xy]-xy)\underline{\tau }(y^{*}(^{*}x))  \\
 &=\mathfrak{g}(\tau [P])+\sum _{x\in X_{k},y\in _{k}Y}([xy]-xy)\tau (y^{*}(^{*}x))+\sum _{(x,y)\in \mathcal{W}}([xy]-xy)\underline{\tau }(y^{*}(^{*}x)) \\
&= \varphi \mathfrak{g}(\mu _{k}^{2}(P))+\sum _{(x,y)\in \mathcal{W}}([xy]-xy)\underline{\tau }(y^{*}(^{*}x)).
 \end{align*}
Therefore
$$\underline{\phi }\underline{\varphi }\underline{\mathfrak{f}}(\underline{\mu }_{k}^{2}(P))=\phi \varphi \mathfrak{g}(\mu _{k}^{2}(P))+\sum _{(x,y)\in \mathcal{W}}\underline{\phi }([xy]-xy)\underline{\tau }(y^{*}(^{*}x)).$$ Applying Theorem \ref{theo8.1} to $\mathcal{F}_{E}(M)$ yields
$$\underline{\phi }\underline{\varphi }\underline{\mathfrak{f}}(\underline{\mu }_{k}^{2}(P))=P+\sum _{x\in \hat{X}_{k},y\in _{k}\hat{Y}}[xy]\underline{\tau }(y^{*}(^{*}x))+\sum _{x\in X_{k}, x_{1}\in _{k}X}[x,x_{1}]f_{[x,x_{1}]}+z_{1},$$
where $f_{[x,x_{1}]}\in \mathcal{F}_{E}(\underline{e}_{k}M\underline{e}_{k}\oplus Me_{k}M)^{\geq 1}$ and $z_{1}\in [\mathcal{F}_{E}(\nu _{k}M),\mathcal{F}_{E}(\nu _{k}M)]$.
Let $A$ be a set of $Z$-free oriented generators of $\underline{M}$, we set $A_{k}=A\cap \underline{M}e_{k}, _{k}A=A\cap e_{k}\underline{M}$. The set of elements of the form $a_{1}ta_{2}$ with $a_{1}\in A_{k}, a_{2}\in _{k}A, t\in L_{k}$ is a set of $Z$-free generators of the $E$-$E$-bimodule $\underline{M}e_{k}\underline{M}$. Likewise, the set $U _{k}$  formed by the elements $a_{1}t$ with $a_{1}\in A_{k}, t\in L_{k}$ is a local basis for $_{E}\underline{M}e_{k}$; and the set $_{k}V$, formed by the elements of the form $ta$ with $a\in _{k}A$ and $t\in L_{k}$, is a local basis for $e_{k}\underline{M}_{E}$. 
  
In what follows, we set $N=Me_{k}M\oplus (e_{k}M)^{*}\otimes _{E_{k}}(^{*}(Me_{k}))$ and $\underline{N}=\underline{M}e_{k}\underline{M}\oplus (e_{k}\underline{M})^{*}\otimes _{E_{k}}(^{*}(\underline{M}e_{k}))$. We have $\underline{N}=N\oplus H_{3}$, with $H_{3}=H_{1}\oplus H_{2}$.
  
As in the observation after Definition \ref{def7.6} we have
$$\sum _{x\in \hat{X}_{k}, y\in _{k}\hat{Y}}[xy]\underline{\tau }(y^{*}(^{*}x))=\sum _{x\in U_{k}, y\in_{k} V}[xy]\underline{\tau }(y^{*}(^{*}x))+z_{2}$$
with $z_{2}\in [\mathcal{F}_{E}(\underline{N}), \mathcal{F}_{E}(\underline{N}))$.  
Let $a_{1}\in A_{k},a_{2}\in _{k}A$, using the equalities of Proposition \ref{prop6.9}, we obtain
 
 \begin{align*}
\sum _{t_{i},t_{j}\in L_{k}}[a_{1}t_{i}t_{j}a_{2}]\underline{\tau }(a_{2}^{*}t_{j}^{-1}t_{i}^{-1}(^{*}a_{1}))&=\sum _{r_{1},r_{2},t_{i},t_{j}\in L_{k}}r_{1}^{*}(t_{i}t_{j})(r_{2}^{-1})^{*}(t_{j}^{-1}t_{i}^{-1})[a_{1}r_{1}a_{2}]\underline{\tau }(a_{2}^{*}r_{2}^{-1}(^{*}a_{1})) \\
&=c_{k}\sum _{r\in L_{k}}[a_{1}ra_{2}]\underline{\tau }(a_{2}^{*}r^{-1}(^{*}a_{1}))
\end{align*}
where $c_{k}=\mathrm{Card}(L_{k})=\mathrm{dim}_{F}(E_{k})$. Therefore
$$\sum _{x\in U_{k},y\in _{k}V}[xy]\underline{\tau }(y^{*}(^{*}x))=c_{k}\sum _{a_{1}\in A_{k},a_{2}\in _{k}A, r\in L_{k} }[a_{1}ra_{2}]\underline{\tau }(a_{2}^{*}r^{-1}(^{*}a_{1})).$$

Now for $x\in X_{k}$ and $x_{1}\in _{k}X$ 
$$[xx_{1}]=\sum _{a_{1}\in A_{k},a_{2}\in _{k}A,r\in L_{k},t_{s}\in L}\alpha ^{a_{1}ra_{2},t_{s}} _{[xx_{1}]}[a_{1}ra_{2}]t_{s},$$
where $\alpha ^{a_{1}ra_{2},t_{s}} _{[xx_{1}]}\in E$. Then
$$\sum _{x\in X_{k},x_{1}\in _{k}X}[xx_{1}]f_{[xx_{1}]}=\sum _{a_{1}\in A_{k},a_{2}\in _{k}A, r\in L_{k}}[a_{1}ra_{2}]g_{a_{1}ra_{2}}+z_{3},$$
where $g_{a_{1}ra_{2}}=\displaystyle \sum _{x\in X_{k},x_{1}\in _{k}X, t\in L}tf_{[xx_{1}]}\alpha _{[xx_{1}]}^{a_{1}ra_{2},t}\in \mathcal{F}_{E}(\underline{e}_{k}M\underline{e}_{k}\oplus Me_{k}M)^{\geq 1}$ and $z_{3}\in [\mathcal{F}_{E}(\underline{N}),\mathcal{F}_{E}(\underline{N})]$. From this, we obtain

\begin{center}
$\underline{\phi }\underline{\varphi }\underline{\mathfrak{f}}(\underline{\mu }_{k}^{2}(P))=P+c_{k}\displaystyle \sum _{a_{1}\in A_{k},a_{2}\in _{k}A,r \in L_{k}}[a_{1}ra_{2}](\underline{\tau }(a_{2}^{*}r^{-1}(^{*}a_{1}))+c_{k}^{-1}g_{a_{1}ra_{2}})+z_{1}+z_{2}+z_{3}.$
\end{center}

Define an $E$-$E$-morphism $\tilde{\varrho }:\nu _{k}\underline{M}\rightarrow \mathcal{F}_{E}(\underline{e}_{k}M\underline{e}_{k}\oplus Me_{k}M)^{\geq 1}$ as the zero map on $\underline{M}\oplus \underline{M}e_{k}\underline{M}$, and on $(e_{k}\underline{M})^{*}e_{k}(^{*}(\underline{M}e_{k}))$ is defined as $\tilde{\varrho }(\underline{\tau }(a_{2}^{*}r^{-1}(^{*}a_{1}))=-c_{k}^{-1}g_{a_{1}ra_{2}}$. The morphism $\tilde{\varrho }=\tilde{\varrho }^{(1)}+\tilde{\varrho }^{(2)}$, where $\tilde{\varrho }^{(1)}:\nu _{k}\underline{M}\rightarrow \nu _{k}\underline{M}$ and $\tilde{\varrho }^{(2)}:\nu _{k}\underline{M}\rightarrow \mathcal{F}_{E}(\underline{e}_{k}M\underline{e}_{k}\oplus Me_{k}M)^{\geq 2}$.  Now we consider the endomorphism of $\mathcal{F}_{E}(\nu _{k}\underline{M})$ determined by the pair $(id_{\nu _{k}\underline{M}}+\tilde{\varrho }^{1}, \tilde {\varrho }^{(2)})$.

Observe that the matrix of $id_{\nu _{k}\underline{M}}+\tilde{\varrho }^{(1)}$ with respect to the decomposition $\nu _{k}\underline{M}=L\oplus (e_{k}\underline{M})^{*}e_{k}(^{*}(\underline{M}e_{k}))$, where $L=\underline{M}\oplus \underline{M}e_{k}\underline{M}$, has the form $\left( \begin{array}{cc}E_{1}&Z\\0&E_{2} \end{array} \right)$ with $E_{1}=id_{L}$ and $E_{2}=id_{(e_{k}\underline{M})^{*}e_{k}(^{*}(\underline{M}e_{k}))}$. Therefore, $id_{\nu _{k}\underline{M}}+\tilde{\varrho }^{(1)}$ is an automorphism of $\nu _{k}\underline{M}$. Using Proposition \ref{prop1.7} we get that $\varrho $ is an automorphism of $\mathcal{F}_{E}(\nu _{k}\underline{M})$.

Then:
$$\varrho \underline{\phi }\underline{\varphi }\underline{\mathfrak{f}}(\underline{\mu }_{k}^{2}(P))=P+c_{k}\sum _{a_{1}\in A_{k},a_{2}\in _{k}A,r\in L_{k}}[a_{1}ra_{2}]\underline{\tau }(a_{2}^{*}r^{-1}(^{*}a_{1}))+w$$
with $w=\varrho (z_{1}+z_{2}+z_{3})$. From this equality we obtain
$$\varrho \left(P+\sum _{x\in X_{k},y\in _{k}Y}[xy]\tau (y^{*}(^{*}x))+\sum _{x\in X_{k},x_{1}\in _{k}X}[xx_{1}]f_{[xx_{1}]}+z_{1}\right)$$
$$+\sum _{(x,y)\in \mathcal{W}}[xy]\varrho (\underline{\tau }(y^{*}(^{*}x)))=P+c_{k}\sum _{a_{1}\in A_{k},a_{2}\in _{k}A, r\in L_{k}}[a_{1}ra_{2}]\underline{\tau }(a_{2}^{*}r^{-1}(^{*}a_{1}))+w.$$

We have $[a_{1}ra_{2}]=w_{a_{1}ra_{2}}+w_{a_{1}ra_{2}}^{\prime }$, with $w_{a_{1}ra_{2}}\in Me_{k}M, w_{a_{1}ra_{2}}^{\prime }\in H_{1}$, $ \underline{\tau }(a_{2}^{*}r^{-1}(^{*}a_{1}))=\hat{w}_{a_{1}ra_{2}}+\hat{w}_{a_{1}ra_{2}}^{\prime }$, with $\hat{w}_{a_{1}ra_{2}}\in (e_{k}M)^{*}e_{k}(^{*}(Me_{k}))$, $\hat{w}_{a_{1}ra_{2}}^{\prime }$ lying in
$H_{2}$.  As we saw before, the restriction of $\varrho $ to $\mathcal{F}_{E}(\nu _{k}M)$ induces an automorphism $\rho $ of $\mathcal{F}_{E}(\nu _{k}M)$. Moreover,  $\displaystyle \sum _{(x,y)\in \mathcal{W}}[xy]\varrho (\underline{\tau }(y^{*}(^{*}x)))$ lies in $\mathcal{H}$. From the above equality, we deduce that $\rho \phi \varphi \mathfrak{g}(\mu _{k}^{2}(P))$ is cyclically equivalent
to $$P+c_{k}\sum _{a_{1}\in A_{k},a_{2}\in _{k}A, r\in L_{k}}w_{a_{1}ra_{2}}\hat{w}_{a_{1}ra_{2}}.$$

Finally, let us prove that $W=c_{k}\displaystyle \sum _{a_{1}\in A_{k},a_{2}\in _{k}A, r\in L_{k}}w_{a_{1}ra_{2}}\hat{w}_{a_{1}ra_{2}}\in \mathcal{F}_{E}(N)$ is a trivial potential. Consider the potential 
$$Q=c_{k}\sum _{a_{1}\in A_{k},a_{2}\in _{k}A, r\in L_{k}}[a_{1}ra_{2}]\underline{\tau }(a_{2}^{*}r^{-1}(^{*}a_{1})).$$
Here, the set of elements $[a_{1}ra_{2}], \underline{\tau }(a_{2}^{*}r^{-1}(^{*}a_{1}))$ with $a_{1}\in A_{k}, a_{2}\in _{k}A$ and $r\in L_{k}$, is a $Z$-free set of free generators for $\underline{N}$. By Proposition \ref{prop3.4}, $\tilde{\delta }_{[a_{1}ra_{2}]^{*}}(Q)=\underline{\tau }(a_{2}^{*}r^{-1}(^{*}a_{1}))$ and $\tilde{\delta }_{\underline{\tau }(a_{2}^{*}r^{-1}(^{*}a_{1}))^{*}}(Q)=[a_{1}ra_{2}]$. 

The potential $Q$ is cyclically equivalent to the potential $Q_{1}+Q_{2}$, where $Q_{1}=\displaystyle \sum _{x\in X_{k},y\in _{k}Y}[xy]\tau(y^{*}(^{*}x))$ and $Q_{2}=\displaystyle \sum _{(x,y)\in \mathcal{W}}[xy]\underline{\tau }(y^{*}(^{*}x))$, $Q_{1}\in \mathcal{F}_{E}(N)$, $Q_{2}\in \mathcal{F}_{E}(H)$.

Then $\tilde{\delta }^{Q}=\tilde{\delta }^{Q_{1}}+\tilde{\delta }^{Q_{2}}$, $\tilde{\delta }^{Q_{1}}(N^{*})\subseteq N$, $\tilde{\delta }^{Q_{1}}(H^{*})=0$, $\tilde{\delta }^{Q_{2}}(H^{*})\subseteq H$ and $\tilde{\delta }^{Q_{2}}(N)=0$. Since $\tilde{\delta }^{Q}(N^{*}\oplus H^{*})=N\oplus H$, we deduce that $\tilde{\delta }^{Q_{1}}(N^{*})=N$. Therefore $Q_{1}$ is a trivial potential and this potential is cyclically equivalent to $W$; thus, $W$ is trivial. The proof is now complete.
\end{proof}
Now we continue with the proof of Theorem \ref{theo10.11}. We have that $\mu _{k}(P)$ is right-equivalent to $\overline{\mu }_{k}(P)\oplus W_{1}$, $\mu _{k}M=L_{1}\oplus L_{2}$, $\overline{\mu }_{k}(P)\in \mathcal{F}_{E}(L_{1})$, $W_{1}\in \mathcal{F}_{E}(L_{2})$, and $W_{1}$ is a trivial potential. Since $\mu _{k}M$ has no two-cycles passing through $k$, we have that $W_{1}=\displaystyle \sum _{s}w_{s}w^{\prime}_{s}$, where $w_{s}\in e_{u}L_{2}e_{v}, w^{\prime }_{s}\in e_{v}L_{2}e_{u}$, so $u\neq k, v\neq k$.  From here we deduce that $\overline{e}_{k}L_{2}\overline{e}_{k}=L_{2}$. Then by Remark \ref{rem7.13}, $\mu _{k}(\overline {\mu }_{k}P\oplus W_{1})$ is right-equivalent to $\mu _{k}\overline{\mu }_{k}(P)\oplus W_{1}$. This last potential is right-equivalent to $\overline{\mu }_{k}^{2}(P)\oplus W_{2}\oplus W_{1}$, with $W_{2}$ a trivial potential. By Theorem \ref{theo7.9}, $\mu _{k}^{2}(P)$ is right-equivalent to $\overline{\mu } _{k}^{2}(P)\oplus W_{2}\oplus W_{1}$. Applying Proposition \ref{prop10.12}, yields that $\mu _{k}^{2}(P)$ is right-equivalent to $P\oplus W$, with $W$ a trivial potential. Therefore $P\oplus W$ is right-equivalent to $\overline{\mu}_{k}^{2}(P)\oplus W_{2}\oplus W_{1}$. Applying Proposition \ref{prop10.9} we get that $\overline{\mu}_{k}^{2}(P)$ is right-equivalent to $P$. This completes the proof. 
\end{proof}

\begin{definition} \label{def10.13} A polarization $\mathcal{A}$ of the $E$-$E$-bimodule $M$ is called dominant if for any $(i,j)\in \mathcal{A}$, $\mathrm{dim}_{F}(e_{j}Me_{i})\leq \mathrm{dim}_{F}(e_{i}Me_{j})$.
\end{definition}
Given a dominant polarization $\mathcal{A}$ of $M$, we set $A^{\sigma}=\{(i,j): (j,i)Ê\in \mathcal{A}, e_{i}Me_{j} \neq 0\}$. Also, we define $M_{\mathcal{A}}=\displaystyle \sum_{(i,j) \in \mathcal{A}} e_{i}Me_{j}$ and $M_{A^{\sigma}}=\displaystyle \sum_{(i,j) \in \mathcal{A}^{\sigma}} e_{i}Me_{j}$.

\begin{definition} \label{def10.14} A potential $P$ in $\mathcal{F}_{E}(M)$ is said to be fully reducible if $P$ is right-equivalent to the direct sum of a trivial potential and a potential $Q \in \mathcal{F}_{E}(N)^{\geq 3}$, where $N$ is a $2$-acyclic $E$-$E$-bimodule.
\end{definition}

\begin{prop} \label{prop10.15} Let $P$ be a potential in $\mathcal{F}_{E}(M)$. Then the following conditions are equivalent
\begin{enumerate}[(i)]
\item The potential $P$ is fully reducible. 
\item One has $\mathrm{dim}_{F}(\Theta (P))=2\mathrm{dim}_{F}M_{\mathcal{A}^{\sigma }},$  for any dominant polarization $\mathcal{A}$ of $M$.
\item There exists a dominant polarization $\mathcal{A}$ of $M$ such that
$\delta ^{P^{(2)}}:(M_{\mathcal{A}^{\sigma }})^{*}\rightarrow M_{\mathcal{A}}$
is a monomorphism.
\end{enumerate}
\end{prop}

\begin{proof} \begin{enumerate} 
\item (i) implies (ii).  Suppose $P$ is right-equivalent to $W+Q$, where $W\in \mathcal{F}_{E}(M_{1})$ is a trivial potential, $Q\in \mathcal{F}_{E}(M_{2})^{\geq 3}$ and $M=M_{1}\oplus M_{2}$; where $M_{2}$ is $2$-acyclic. Let $\mathcal{A}$ be a dominant polarization of $M$ and consider $\mathcal{C}=\{(i,j)\in \mathcal{A} | e_{i}M_{1}e_{j}\neq 0 \}$. The set $\mathcal{C}$ is a polarization of $M_{1}$. Moreover, $\mathcal{C}^{\sigma }=\mathcal{A}^{\sigma }$ and $(M_{1})_{\mathcal{C}^{\sigma }}=M_{\mathcal{A}^{\sigma }}$. Since $W$ is a trivial potential in $\mathcal{F}_{E}(M_{1})$, then $\Theta (W)=M_{1}$, $M_{1}=(M_{1})_{\mathcal{C}^{\sigma }}\oplus M_{\mathcal{C}}$ and $\mathrm{dim}_{F}(M_{\mathcal{C}^{\sigma }})=\mathrm{dim}_{F}(M_{\mathcal{C}})$. Therefore
$$\mathrm{dim}_{F}(\Theta (P))=\mathrm{dim}_{F}(\Theta (W))=2\mathrm{dim}_{F}((M_{1})_{\mathcal{C}^{\sigma }})=2\mathrm{dim}_{F}(M_{\mathcal{A}^{\sigma }}).$$
\item (ii) implies (iii). We have $\delta ^{P^{2}}(M_{\mathcal{A}^{\sigma }}^{*}) \subseteq M_{\mathcal{A}}$ and $\delta ^{P^{(2)}}(M_{\mathcal{A}}^{*}) \subseteq M_{\mathcal{A}^{\sigma }}$.
Then 
\begin{align*}
\mathrm{dim}_{F}(\Theta (P))&=\mathrm{dim}_{F}(\delta ^{P^{2}})(M^{*}) \\
&=\mathrm{dim}_{F}\delta ^{P^{2}}(M_{\mathcal{A}^{\sigma }}^{*})+\mathrm{dim}_{F}(\delta ^{P^{(2)}}(M_{\mathcal{A}}^{*})) \\
&\leq 2\mathrm{dim}_{F}(M_{\mathcal{A}^{\sigma }})
\end{align*}
Then by our hypothesis on $P$, we have $\mathrm{dim}_{F}\delta ^{P^{2}}((M_{\mathcal{A}^{\sigma }})^{*})=\mathrm{dim}_{F}(M_{\mathcal{A}^{\sigma }})$. This implies (iii).
\item (iii) implies (i). Since $\tilde{\delta }^{P^{2}}:(M_{\mathcal{A}^{\sigma }})^{*}\rightarrow M_{\mathcal{A}}$ is a monomorphism, then for $(i,j)\in \mathcal{A}^{\sigma }$
 $$\mathrm{dim}_{F}(e_{j}\Theta (P)e_{i})=\mathrm{dim}_{F}\tilde{\delta }^{P^{(2)}}((e_{i}Me_{j})^{*})=\mathrm{dim}_{F}(e_{i}Me_{j}).$$
 
 By Theorem \ref{theo10.4}, there exists an automorphism of topological algebras $\varphi :\mathcal{F}_{E}(M)\rightarrow \mathcal{F}_{E}(M_{1}\oplus M_{2})$, with $\varphi |_{E}=id_{E}$, such that $\varphi (P)$ is cyclically equivalent to $Q\oplus W$ with $Q\in \mathcal{F}_{E}(M_{1})^{\geq 3}$ and $W$ is  a quadratic trivial potential in $\mathcal{F}_{E}(W)$. We have   $M\cong  M_{1}\oplus M_{2}$. Then $\Theta (\varphi (P))=\Theta (Q\oplus W)=\Theta (W)=M_{2}.$ By Proposition \ref{prop10.8}, 
 $\Theta (\varphi (P))=\varphi ^{(1)}\Theta (P)$, where $\varphi $ is determined by the pair $(\varphi ^{(1)}, \varphi ^{(2)})$ . Therefore, for $(i,j)\in \mathcal{A}^{\sigma }$, we have
 $$\mathrm{dim}_{F}(e_{i}M_{2}e_{j})=\mathrm{dim}_{F}(e_{j}M_{2}e_{i})=\mathrm{dim}_{F}e_{j}\Theta (\varphi (P))e_{i}=\mathrm{dim}_{F}(e_{i}Me_{j}).$$
and then $\mathrm{dim}_{F}M_{\mathcal{A}^{\sigma }}=\mathrm{dim}_{F}(M_{2})_{\mathcal{A}^{\sigma }}.$
This implies $\mathrm{dim}_{F}(M_{1})_{\mathcal{A}^{\sigma }}=0$, so $M_{1}$ is $2$-acyclic and thus $P$ is fully reducible. This proves (i).
 \end{enumerate}   
\end{proof}
 
Now suppose $Q$ is a quiver without loops nor $2$-cycles and such that for each pair of vertices $i,j$ in $Q_{0}$, there is at most one arrow from $i$ to $j$. Let $\{M_{\rho }\}_{\rho \in \mathcal{B}}$ be  a family of $E$-$E$-bimodules satisfying the conditions of Section \ref{sec10} and with the additional condition that $\mathrm{Hom}_{E-E}(M_{\rho _{1}}, M_{\rho _{2}})\neq 0$ implies $\rho _{1}=\rho _{2}$ (this condition holds for examples 9.4 and 9.5).  Let $g$ be a modulation of $Q$ on the above family of $E$-$E$-bimodules. Consider $P$ a potential in $\mathcal{F}_{E}(M_{g})^{\geq 3}$, with $\mu _{k}P$  fully reducible. Then there exists a dominant polarization $\mathcal{A}$ of $M_{\mu _{k}g}$ such that for $\alpha: j\rightarrow i$ with $(i,j)\in \mathcal{A}$, either there is no arrow from $i$ to $j$ or there is an arrow $\alpha ^{\sigma }:i\rightarrow j$ such that 
$$\mu _{k}g(\alpha ^{\sigma }, \overline{\rho })\leq \mu _{k}g(\alpha , \rho )$$
for all $\rho \in \mathcal{B}$.
In this case, define $\tilde{\mu }_{k}Q$ eliminating from $\mu _{k}Q$ the arrows $i\rightarrow j$ with $(i,j)\in \mathcal{A}$. Now define $\tilde{\mu }_{k}g$, the modulation of $\tilde{\mu }_{k}Q$, as follows: for $(i,j)\in \mathcal{A}$, $\tilde{\mu}_{k}g(\alpha , \rho )=\mu _{k}g(\alpha , \rho )$ in case there is no arrow in $\mu _{k}Q$ from $i$ to $j$; otherwise, $\tilde{\mu }_{k}g(\alpha , \rho )=\mu _{k}g(\alpha , \rho )-\mu _{k}g(\alpha ^{\sigma }, \overline{\rho })$. Then $\tilde{\mu }_{k}M_{g}=M_{\tilde{\mu }_{k}g}$.

\end{section}

\begin{section}{Non-degeneracy} \label{sec12}
In this section we assume $E=\displaystyle \prod_{i=1}^{n} E_{i}$ as in section $11$, we also adopt the notation used there. Throughout this section, we assume $F$ is an infinite field. Here we consider the case in which $M$ is a $Z$-free $E$-$E$-bimodule. Our main purpose is to prove that for each sequence of integers $k_{l},\ldots,k_{1}$ of $[1,n]$ there exist potentials which are $(k_{l},\ldots,k_{1})$-non-degenerate in the sense of Definition $7.2$ of \cite{4}. For this, we use the main ideas given in \cite{4}. This result does not hold in general for $E$-$E$-bimodules which are not $Z$-free, as has been observed in section $12$ of \cite{3} and in example $4.1$ of \cite{7}.

Let $B$ be a non-empty set. We denote by $F^{B}$ the $F$-vector space consisting of all functions $B \rightarrow F$ and by $F[x]_{x \in B}$ the free commutative $F$-algebra on the set $B(T)$.
\begin{definition} \label{def11.1} A function $g: F^{B} \rightarrow F$ is polynomial if there exists a polynomial $G_{g}=\displaystyle \sum_{i=0, \ldots, i_{l}=0}^{n_{1},\ldots,n_{l}} a_{i_{1},\ldots, i_{l}} x_{1}^{i_{1}} \cdots x_{l}^{i_{l}}$ in $F[x]_{x \in B}$, with $a_{i_{1}, \ldots, i_{l}} \in F$, such that
\begin{center}
$G_{g}(u)=\displaystyle \sum_{i_{1}=0,\ldots,i_{l}=0}^{n_{1},\ldots,n_{l}} a_{i_{1},\ldots,i_{l}} u(x_{1})^{i_{1}} \cdots u(x_{l})^{i_{l}}$
\end{center}
for all $u \in F^{B}$.
\end{definition}

\begin{definition} \label{def11.2} Let $B$ and $B'$ be non-empty sets. We say that a function $h: F^{B} \rightarrow F^{B'}$ is polynomial if for each $x \in B'$, the function $h_{x}: F^{B} \rightarrow F$ given by $h_{x}(u)=h(u)(x)$ for all $u \in F^{B}$, is a polynomial function.
\end{definition}

\begin{rem} \label{rem11.3} If $B'$, $B''$ are non-empty sets and $h_{1}: F^{B} \rightarrow F^{B'}$, $h_{2}: F^{B'} \rightarrow F^{B''}$ are polynomial functions, then $h_{2}h_{1}: F^{B} \rightarrow F^{B''}$ is a polynomial function.
\end{rem}

Let $M$ be an $E$-$E$-bimodule with $A$ a $Z$-local-free basis. We put $A(i,j)=A \cap e_{i}Me_{j}$. Let $B_{1}$ be the set of all elements $z \in M$ of the form $z=t_{1}(z)a(z)t_{2}(z)$ with $a(z) \in A$, $t_{1}(z) \in L(\tau(a(z)))$ and $t_{2}(z) \in L(\sigma(a(z)))$. Clearly, $B_{1}$ is an $F$-basis for $M$. Now for $m \geq 2$ take $B_{m}$, the set of all elements $x \in M^{\otimes m}_{cyc}$ of the form $x=t_{1}(x)a_{1}(x)t_{2}(x)a_{2}(x) \cdots t_{m}a_{m}(x)t_{m+1}(x)$, with $a_{1}(x), \ldots, a_{m}(x) \in A$, $t_{1}(x) \in L(\tau(a_{1}(x)))$, $t_{2}(x) \in L(\tau(a_{2}(x)))$,$\ldots$, $t_{m}(x) \in L(\sigma(a_{m}(x)))$. Then $B_{m}$ is an $F$-basis for $M^{\otimes m}_{cyc}$ and we set $B = \displaystyle \bigcup_{m=2}^{\infty} B_{m}$. \\
If $P$ is a potential in $\mathcal{F}_{E}(M)$, then $P=\displaystyle \sum_{m=2}^{\infty} \displaystyle \sum_{x \in B_{m}} c_{x}x$ with $c_{x} \in F$. Therefore to each potential $P$ corresponds a unique element $\underline{c}(P) \in F^{B}$ such that $\underline{c}(P)(x)=c_{x}$. Conversely, if $u \in F^{B}$, we have the potential $P_{u}=\displaystyle \sum_{m=2}^{\infty} \displaystyle \sum_{x \in B_{m}} u(x)x$, with $\underline{c}(P_{u})=u$.  \\

Let $M'$ be an $E$-$E$-bimodule with $A'$ a $Z$-local-free basis. Let $B_{m}'$ be the corresponding $F$-basis for $(M')_{cyc}^{\otimes m}$. Set $B'= \displaystyle \bigcup_{m=2}^{\infty} B_{m}'$. 

\begin{lemma} \label{lem11.4} Let $\phi: \mathcal{F}_{E}(M)_{cyc} \rightarrow \mathcal{F}_{E}(M')_{cyc}$ be an $F$-linear map such that $\phi(\mathcal{F}_{E}(M)^{\geq m}_{cyc}) \subseteq \mathcal{F}_{E}(M')^{\geq m}_{cyc}$ for all $m \geq 2$. Then $\phi$ is a continuous map, and there exists a polynomial function

\begin{center}
$\underline{\phi}: F^{B} \rightarrow F^{B'}$
\end{center}
such that for any potential $P \in \mathcal{F}_{E}(M)$
\begin{center}
$\underline{c}(\phi(P))=\underline{\phi}(\underline{c}(P))$.
\end{center}
\end{lemma}

\begin{proof} In what follows, we set $B^{\leq m}=\displaystyle \bigcup_{l=2}^{m} B_{l}$ and $B^{\geq m}=\displaystyle \bigcup_{l=m}^{\infty} B_{l}$, and analogous definitions for $(B')^{\leq m}$ and $(B')^{\geq m}$. For $x \in B_{m}$, we have $\phi(x)=\displaystyle \sum_{y \in (B')^{\geq m}} \alpha_{x,y}y$. Let $\underline{\phi}: F^{B} \rightarrow F^{B'}$ be defined as follows. For $u \in F^{B}$ and $y \in B'_{m}$, we define

\begin{center}
$\underline{\phi}(u)(y)=\displaystyle \sum_{x \in B^{\leq m}}u(x)\alpha_{x,y}$
\end{center}
Take $u=\underline{c}(P)$, then $P=\displaystyle \sum_{m=2}^{\infty} \left(\displaystyle \sum_{x \in B_{m}} u(x)x \right)$. Therefore
\begin{align*}
\phi(P)&=\displaystyle \sum_{m=2}^{\infty} \left( \displaystyle \sum_{x \in B_{m}} u(x)\phi(x) \right) \\
&=\displaystyle \sum_{m=2}^{\infty} \displaystyle \sum_{x \in B_{m}} \left( \displaystyle \sum_{y \in (B')^{\geq m}} \alpha_{x,y}yÊ\right) \\
&=\displaystyle \sum_{m=2}^{\infty} \displaystyle \sum_{y \in B_{m}'} \left( \displaystyle \sum_{x \in B^{\leq m}} u(x)\alpha_{x,y}\right)y \\
&=\displaystyle \sum_{yÊ\in B'} \underline{\phi}(u)(y)y
\end{align*}
\end{proof}

\begin{definition} \label{def11.5} For a polynomial $T \in F[x]_{x \in B}$, we denote by $\mathcal{Z}(T)$ the set of all functions $u \in F^{B}$ such that $T(u) \neq 0$. A function $h: \mathcal{Z}(T) \rightarrow F$ is called regular if there exists a polynomial function $G_{h}: \mathcal{Z}(T) \rightarrow F$ and a non-negative integer $m$ such that for all $u \in \mathcal{Z}$, $h(u)=G_{h}(u)/T(u)^{m}$. If $B'$ is a non-empty set, a function $H: \mathcal{Z}(T) \rightarrow F^{B'}$ is called regular if for all $x \in B'$, the function $H_{x}: F^{B} \rightarrow F$, sending every $u \in F^{B}$ to $H(u)(x)$, is regular.
\end{definition}

\begin{rem} \label{rem11.6}
Suppose $B'$ and $B''$ are non-empty sets and $T$ is a non-zero polynomial in $F[x]_{x \in B}$. If $h: \mathcal{Z}(T) \rightarrow F^{B'}$ is a regular function and $h_{1}: F^{B'} \rightarrow F^{B''}$ is a polynomial function, then $h_{1}h$ is a regular function.
\end{rem}

\begin{rem} \label{rem11.7}
Let $T_{1} \in F[x]_{x \in B}$ and $T_{2} \in F[x]_{x \in B'}$ be non-zero polynomials and $g_{1}: \mathcal{Z}(T_{1}) \rightarrow F^{B'}$, $g_{2}: \mathcal{Z}(T_{2}) \rightarrow F^{B''}$ be regular functions. Suppose there exists $u_{0} \in \mathcal{Z}(T_{1})$ with $g_{1}(u_{0}) \in \mathcal{Z}(T_{2})$, then there exists a non-zero polynomial $H$, with $\mathcal{Z}(H) \subseteq \mathcal{Z}(T_{1})$, such that $g_{1}(\mathcal{Z}(H)) \subseteq \mathcal{Z}(T_{2})$. Moreover, the composition
\begin{center}
$g_{2}(g_{1})_{|\mathcal{Z}(H)}: \mathcal{Z}(H)Ê\rightarrow F^{B''}$
\end{center}
is regular.
\end{rem}

\begin{proof} Note that $T_{2}$ induces a polynomial function $\underline{T_{2}}: F^{B'}Ê\rightarrow F$, so $\underline{T_{2}}g_{1}$ is a regular function. Therefore, there exists a polynomial $G \in F[x]_{x \in B}$ and a non-negative integer $m$ such that $\underline{T_{2}}g_{1}(u)=G(u)/T_{1}(u)^{m}$ for all $u \in \mathcal{Z}(T_{1})$. It follows that $\underline{T_{2}}g_{1}(u_{0})=T_{2}(g_{1}(u_{0})) \neq 0$; therefore $G(u_{0}) \neq 0$ and hence $G \neq 0$. Now consider the product $H=T_{1}G$. We have that $H \neq 0$ and $\mathcal{Z}(H) \subseteq \mathcal{Z}(T_{1})$. One can now verify that the composition $g_{2}(g_{1})|_{\mathcal{Z}(H)}$ is regular.
\end{proof}

\begin{lemma} \label{lem11.8} Let $T$ be a polynomial in $F[x]_{x \in B}$ and suppose that for any $u \in \mathcal{Z}(T)$ we have an algebra automorphism $\varphi^{u}$ of $\mathcal{F}_{E}(M)$ such that $\varphi^{u}|_{E}=id_{E}$, and for any $x \in B_{1}$
\begin{center}
$\varphi^{u}(x)=\displaystyle \sum_{m=1}^{\infty} \displaystyle \sum_{z \in e_{\tau(x)}B_{m}e_{\sigma(x)}} G_{x,z}(u)z$
\end{center}
where the functions that send any $u \in \mathcal{Z}(T)$ to $G_{x,z}(u)$ are regular. Then there exists a regular function $h: \mathcal{Z}(T) \rightarrow F^{B}$ such that for a potential $P$, with $\underline{c}(P) \in \mathcal{Z}(T)$,
\begin{center}
$\underline{c}(\varphi^{\underline{c}(P)})=h(\underline{c}(P))$.
\end{center}
\end{lemma}

\begin{proof} Take $P$ a potential with $u=\underline{c}(P)Ê\in \mathcal{Z}(T)$. Then we have $P=\displaystyle \sum_{m=2}^{\infty} \displaystyle \sum_{x \in B_{m}} u(x)x$ and $\varphi^{u}(P)=\displaystyle \sum_{m=2}^{\infty} \displaystyle \sum_{x \in B_{m}} u(x)\varphi^{u}(x)$. \\
If $x \in B_{m}$, $x=x_{1}x_{2} \cdots x_{m}$, with $x_{i} \in B_{1}$, then
\begin{align*}
\varphi^{u}(x)&=\displaystyle \prod_{j=1}^{m} \varphi^{u}(x_{j}) \\
&=\displaystyle \sum_{z_{i_{1}},\ldots,z_{i_{m}} \in B} U_{z_{i_{1}},\ldots,z_{i_{m}}}(u)z_{i_{1}} \cdots z_{i_{m}} \\
&=\displaystyle \sum_{l=m}^{\infty} \displaystyle \sum_{w \in B_{l}} U_{z_{i_{1}},\ldots,z_{i_{m}}}(u)\lambda_{w}^{z_{i_{1}},\ldots,z_{i_{m}}}w
\end{align*}
Then
\begin{center}
$\varphi^{u}(P)=\displaystyle \sum_{m=2}^{\infty} \displaystyle \sum_{w \in B_{m}} \displaystyle \sum_{x \in B^{\leq m}} u(x)G_{x,w}(u)w$
\end{center}
Take now $h: \mathcal{Z}(T) \rightarrow F^{B}$ such that for $u \in \mathcal{Z}(T)$ and $w \in B_{l}$,
\begin{center}
$h(u)(w)=\displaystyle \sum_{x \in B^{\leq l}} u(x)G_{x,w}(u)=\overline{G}_{w}(u)$
\end{center}
Clearly, $\overline{G}_{w}(u)$ is a regular function in $u$. This proves the lemma.
\end{proof}

\begin{prop} \label{prop11.9} Let $\mathcal{A}$ be a dominant polarization for $M$, with $i: M_{\mathcal{A}^{\sigma}} \rightarrow M_{\mathcal{A}}$ the associated inclusion. Suppose that we have a decomposition of $E$-$E$-bimodules
\begin{center}
$M=\mathcal{L}_{1} \oplus \mathcal{L}_{2}$
\end{center}
such that $\mathcal{L}_{1} \cong (M_{A^{\sigma}})^{\ast}$. Then there exists a polynomial $T \in F[x]_{x \in B}$ such that
\begin{enumerate}
\item If $P$ is a potential in $\mathcal{F}_{E}(M)$ with $T(\underline{c}(P)) \neq 0$, then the composition of $E$-$E$-morphisms
\begin{center}
$(M_{A^{\sigma}})^{\ast}\stackrel{\tilde{\delta}^{P^{(2)}}}{\longrightarrow}M_{\mathcal{A}} \stackrel{\pi}{\longrightarrow}\mathcal{L}_{1}$
\end{center}
is an isomorphism, where $\pi$ is the projection onto $\mathcal{L}_{1}$.
\item There exists a regular function $\underline{\xi}: \mathcal{Z}(T) \rightarrow F^{B}$, such that for any potential $P \in \mathcal{F}_{E}(M)$, with $T(\underline{c}(P)) \neq 0$, there is a potential $\xi(P)=\displaystyle \sum_{aÊ\in \mathcal{A}^{\sigma}} ai(a)+\xi(P)^{\geq 3}$ with $\underline{\xi}(\underline{c}(P))=\underline{c}(\xi(P))$.
\end{enumerate}
\end{prop}

\begin{proof} Let $Y$ be the set of all non-zero elements $y$ of $(M_{A^{\sigma}})^{\ast}$ of the form $y=t_{1}(y)(a(y)^{\ast})t_{2}(y)$ with $t_{1}(y)$, $t_{2}(y) \in L$, and $a(y) \in M_{A^{\sigma}}$. The set $Y$ is an $F$-basis for $(M_{A^{\sigma}})^{\ast}$. There exists an isomorphism $(M_{A^{\sigma}})^{\ast} \rightarrow \mathcal{L}_{1}$. Take $W_{1}$ to be the image of $Y$ under this isomorphism. We have $M_{\mathcal{A}_{1}} \cong (M_{A^{\sigma}})^{\ast}$, thus $\mathcal{L}_{2} \cong M_{\mathcal{A}_{2}}$. Therefore there exists an isomorphism $\phi_{2}: \mathcal{L}_{2} \rightarrow M_{\mathcal{A}_{2}}$. Now define $W_{2}=\phi_{2}^{-1}(B_{1} \cap M_{\mathcal{A}_{2}})$. 
We also have an isomorphism $\phi: \mathcal{L}_{1} \rightarrow M_{\mathcal{A}}$ such that $\phi_{1}(W_{1})=B_{1} \cap M_{\mathcal{A}_{1}}$. We take $W=W_{1} \cup W_{2}$ an $F$-basis for $M_{\mathcal{A}}$. For $a \in M_{\mathcal{A}^{\sigma}}$ and $x \in B_{2}$, one has
\begin{center}
$\tilde{\delta}_{a^{\ast}}(x)=\displaystyle \sum_{w \in W} c_{a,w}(x)w$
\end{center}
For a potential $P$, with $u=\underline{c}(P)$, we have $P=\displaystyle \sum_{m=2}^{\infty} \displaystyle \sum_{x \in B_{m}} u(x)x$. Then
\begin{center}
$\tilde{\delta}^{P^{(2)}}(a^{\ast})=\tilde{\delta}_{a^{\ast}}(P^{(2)})=\displaystyle \sum_{x \in B_{2}} u(x)\tilde{\delta}_{a^{\ast}}(x)=\displaystyle \sum_{x \in B_{2},w \in W}u(x)c_{a,w}(x)w$.
\end{center}
Then for $y \in Y$ we have
\begin{align*}
\tilde{\delta}^{P^{(2)}}(y)&=\displaystyle \sum_{xÊ\in B_{2}} u(x)\tilde{\delta}_{a^{\ast}}(x) \\
&=\displaystyle \sum_{x \in B_{2},w \in W} u(x)c_{a(y),w}(x)t_{1}(y)wt_{2}(y) \\
&=\displaystyle \sum_{x \in B_{2},w,w_{1} \in W} u(x)c_{a(y),w}(x)\lambda_{w_{1}}^{y,w}w_{1} \\
&=\displaystyle \sum_{w_{1} \in W} k_{y,w_{1}}(u)w_{1}
\end{align*}
where $t_{1}(y)wt_{2}(y)=\displaystyle \sum_{w_{1} \in W} \lambda_{w_{1}}^{y,w}w_{1}$, with $\lambda_{w_{1}}^{y,w} \in F$ and $k_{y,w_{1}}$ is the following polynomial in $F[x]_{xÊ\in B}$
\begin{center}
$k_{y,w_{1}}=\displaystyle \sum_{xÊ\in B_{2},w \in W}xc_{a(y),w}(x)\lambda_{w_{1}}^{y,w}$
\end{center}
We have the square matrix $A=(k_{y,w})_{yÊ\in Y, wÊ\in W_{1}}$ and the polynomial $T$ in $F[x]_{x \in B}$, given by $T=\operatorname{det}(A)$.
\begin{enumerate}[(i)]
\item The composition $\pi \tilde{\delta}^{P^{(2)}}: (M_{\mathcal{A}^{\sigma}})^{\ast} \rightarrow \mathcal{L}_{1}$ is an isomorphism if and only if $T(\underline{c}(P))Ê\neq 0$. This proves $1$.
\item Let $u \in \mathcal{Z}(T)$ and consider the potential $P_{u}$. Define $Y^{u}$ as the set consisting of the elements $\tilde{\delta}^{P_{u}^{(2)}}(y),w_{2},y \in Y, w \in W_{2}$. The elements of this set can be expressed in terms of the $F$-basis $W$ of $M_{\mathcal{A}}$ by the matrix
$$
\begin{pmatrix}
A(u) & 0 \\
B(u) & I \\
\end{pmatrix}
$$
where $A(u)=(k_{y,w}(u))_{y \in Y, w \in W_{1}}$. Since $T(u)=\operatorname{det}A(u) \neq 0$, then the matrix $A(u)$ is invertible; thus the elements of $W$ can be expressed in terms of the elements $Y^{u}$, by means of the matrix
$$
\begin{pmatrix}
A(u)^{-1} & 0 \\
-B(u)A(u)^{-1} & I
\end{pmatrix}
$$
Therefore for $w \in W$, there exist elements $\lambda_{w,y}(u)$, $\lambda_{w,w_{2}}(u) \in F$, with
\begin{center}
$w=\displaystyle \sum_{y \in Y} \lambda_{w,y}(u)\tilde{\delta}^{P_{u}}(y)+\displaystyle \sum_{w_{2} \in W_{2}} \lambda_{w,w_{2}}(u)w_{2}$
\end{center}
\end{enumerate}
Moreover, the functions sending $u \in \mathcal{Z}(T)$ into $\lambda_{w,y}(u)$, $\lambda_{w,w_{2}}(u)$ are regular functions. The elements $\tilde{\delta}^{P_{u}^{(2)}}(a^{\ast})$ and $\phi_{2}^{-1}(b)$, with $a \in M_{\mathcal{A}^{\sigma}}$ and $b \in M_{\mathcal{A}_{2}}$, form a set of $Z$-free generators for $M_{\mathcal{A}}$. Therefore, there exists an isomorphism $\varphi^{u}: M \rightarrow M$ such that $\varphi^{u}(a)=a$ for $a \in M_{\mathcal{A}^{\sigma}}$; $\varphi^{u}(\tilde{\delta}^{P_{u}^{(2)}}(a^{\ast}))=i(a)$ for $a \in M_{\mathcal{A}^{\sigma}}$; and $\varphi^{u}(\phi_{2}^{-1}(b))=b$ for $b \in M_{\mathcal{A}_{2}}$. Now let $x \in B_{1}$, then $x=\displaystyle \sum_{w \in W} \mu_{x,w}w$ with $\mu_{x,w} \in F$. Therefore
\begin{align*}
\varphi^{u}(x)&=\displaystyle \sum_{wÊ\in W, yÊ\in Y} \mu_{x,w} \lambda_{w,y}(u)\psi(y) + \displaystyle \sum_{w \in W, w_{2} \in W_{2}} \mu_{x,w} \lambda_{w,w_{2}}(u)\phi_{2}(w_{2}) \\
&=\displaystyle \sum_{x_{1} \in B_{1} \cap M_{\mathcal{A}}} \underline{\nu}_{x,x_{1}}(u)x_{1}
\end{align*}
where $\underline{\nu}_{x,x_{1}}(u)$ is a regular function in $u$. By Lemma \ref{lem11.8}, there exists a regular function $h: \mathcal{Z}(T) \rightarrow F^{B}$ such that $\underline{c}(\varphi^{u}(P_{u}))=h(u)$. Now consider the $F$-linear map $\zeta: \mathcal{F}_{E}(M)_{cyc} \rightarrow \mathcal{F}_{E}(M)_{cyc}$ such that for $x \in B_{m}$, with $m>2$, $\zeta(x)=x$; for $x=t_{1}(x)a_{1}(x)t_{2}(x)a_{2}(x)t_{3}(x) \in B_{2}$, we define $\zeta(x)=a_{1}(x)t_{2}(x)a_{2}(x)t_{3}(x)t_{1}(x)$ in case $a_{1}(x) \in M_{\mathcal{A}^{\sigma}}$; otherwise $a_{2}(x) \in M_{\mathcal{A}^{\sigma}}$, and in this case we define $\zeta(x)=a_{2}(x)t_{3}(x)t_{1}(x)a_{1}(x)t_{2}(x)$.
Observe that if $v=\underline{c}\zeta(P)$, then $P_{v}=\zeta(P_{u})$, so $P_{v}$ and $P_{u}$ are cyclically equivalent, hence $\tilde{\delta}^{P_{v}^{(2)}}=\tilde{\delta}^{P_{u}^{(2)}}$. Therefore $v \in \mathcal{Z}(T)$ and $\varphi^{v}(\tilde{\delta}^{P_{v}^{(2)}}(a^{\ast}))=i(a)$; also, $\varphi^{u}(\tilde{\delta}^{P_{v}^{(2)}}(a^{\ast}))=\varphi^{u}(\tilde{\delta}^{P_{u}^{(2)}}(a^{\ast}))=i(a)$. From here we obtain the equality $\varphi^{v}=\varphi^{u}$. By Proposition \ref{prop6.19}, we have $\zeta(P)=\displaystyle \sum_{a \in M_{\mathcal{A}^{\sigma}}} a\tilde{\delta}^{P_{u}^{(2)}}(a^{\ast})+\zeta(P)^{\geq 3}$. Therefore
\begin{center}
$\varphi^{u}(\zeta(P))=\displaystyle \sum_{a \in M_{\mathcal{A}^{\sigma}}}ai(a)+\varphi^{u}(\zeta(P)^{\geq 3})=\displaystyle \sum_{a \in M_{\mathcal{A}^{\sigma}}}ai(a)+\varphi^{u}(\zeta(P))^{\geq 3}.$
\end{center}
By the above, there exists a regular function $h:Ê\mathcal{Z}(T) \rightarrow F^{B}$ with $\underline{c}(\varphi^{u}(\zeta(P)))=h(\underline{c}(\zeta(P)))=h\underline{\zeta}(u)$, where $\underline{\zeta}: F^{B} \rightarrow F^{B}$ is a polynomial function. Taking $\underline{\xi}=h\underline{\zeta}$ and $\xi(P)=\varphi^{\underline{c}(P)}(\zeta(P))$ yields item $2$.
\end{proof}

\begin{prop} \label{prop11.10} There exists a polynomial $T \in F[x]_{x \in B}$ and a regular function $\underline{\chi}: \mathcal{Z}(T) \rightarrow F^{B}$ such that for any potential $P$ with $\underline{c}(P) \in \mathcal{Z}(T)$, there exists a potential $\chi(P)=\displaystyle \sum_{a \in M_{\mathcal{A}^{\sigma}}} ai(a)+\chi(P)^{\geq 3}$ such that
\begin{enumerate}
\item $\chi(P)$ is right-equivalent to $P$.
\item $\underline{c}(\chi(P))=\underline{\chi}(\underline{c}(P))$.
\item $\chi(P)^{\geq 3} \in \mathcal{F}_{E}(M_{\mathcal{A}_{2}}) \subseteq \mathcal{F}_{E}(M)$.
\end{enumerate}
\end{prop}

\begin{proof} We define the $F$-linear map $\rho: \mathcal{F}_{E}(M)_{cyc} \rightarrow \mathcal{F}_{E}(M)_{cyc}$ such that for $x=t_{1}(x)a_{1}(x) \cdots t_{m}(x)a_{m}(x)t_{m+1}(x) \in B_{m}$, $\rho(x)=x$ if all the $a_{i}(x) \in M_{\mathcal{A}_{2}}$. Otherwise, some $a_{i}(x) \in M_{\mathcal{A}^{\sigma}}$ or there is no $a_{i}(x) \in M_{\mathcal{A}^{\sigma}}$, but there is some $a_{i}(x) \in M_{\mathcal{A}_{1}}$. In the first case, we take the minimal $i$ such that $a_{i}(x) \in M_{\mathcal{A}^{\sigma}}$; then if $i=1$, $\rho(x)=a_{1}(x)t_{2}(x) \cdots a_{m}(x)t_{m+1}(x)t_{1}(x)$; if $i>1$, $\rho(x)=a_{i}(x)t_{i+1}(x) \cdots a_{m}(x)t_{m}(x)t_{1}(x)a_{1}(x) \cdots a_{i-1}(x)t_{i}(x)$. In the second case, we take $j$ maximal such that $a_{j}(x) \in M_{\mathcal{A}_{1}}$. If $j=m$, $\rho(x)=t_{m+1}(x)t_{1}(x)a_{1}(x) \cdots a_{m}(x)$; if $j<m$, then $\rho(x)=t_{j+1}(x)a_{j+1}(x) \cdots a_{j}(x)$. \\
Now for $u \in F^{B}$, we define the automorphism $\psi^{u}$ of $\mathcal{F}_{E}(M)$ such that for $a \in M_{\mathcal{A}_{2}}$, $\psi^{u}(a)=a$; and for $a \in M_{\mathcal{A}^{\sigma}}$
\begin{align*}
\psi^{u}(a)&=a-\displaystyle \sum_{l=2}^{\infty} \displaystyle \sum_{x \in e_{\tau(a)}B_{l}e_{\sigma(a)}} u(xi(a))x \\
\psi^{u}(i(a))&=i(a) - \displaystyle \sum_{l=2}^{\infty} \displaystyle \sum_{x \in e_{\sigma(a)}B_{l}e_{\tau(a)}} u(ax)x.
\end{align*}
For any potential $P$, we define $\vartheta(P)=\psi^{\underline{c}(P)}(P)$. By Lemma \ref{lem11.8}, there exists a polynomial function $\underline{\vartheta}: F^{B} \rightarrow F^{B}$ such that for any potential $P$, we have
\begin{center}
$\underline{c}(\vartheta(P))=\underline{\vartheta}(\underline{c}(P))$.
\end{center}
By $2$ of Proposition \ref{prop11.9}, there exists a polynomial $T \in F[x]_{x \in B}$ and a regular function $\underline{\xi}: \mathcal{Z}(T) \rightarrow F^{B}$ such that for any potential $P$ with $T(\underline{c}(P)) \neq 0$, there exists a potential $\xi(P)$, which is right-equivalent to $P$, such that $\chi(P)=\displaystyle \sum_{a \in M_{\mathcal{A}^{\sigma}}} ai(a) + \chi(P)^{\geq 3}$ and $\underline{c}(\chi(P))=\underline{\chi}(\underline{c}(P))$. \\
Since $P$ is a potential with $T(\underline{c}(P)) \neq 0$, then
\begin{center}
$\displaystyle \lim_{l \to \infty} (\vartheta \rho)^{l} \zeta(P)=\chi(P)=\displaystyle \sum_{a \in M_{\mathcal{A^{\sigma}}}} ai(a)+\chi(P)^{\geq 3}$
\end{center}
with $(\chi(P))^{\geq 3} \in \mathcal{F}_{E}(M_{\mathcal{A}_{2}})$. 
Finally, let us prove that for any $x \in B_{m}$, the function which sends any $u \in \mathcal{Z}(T)$ into $\underline{c}(\chi(P_{u}))(x)$ is a regular function. By (ii) of Claim $2$ of the proof of the Theorem \ref{theo5.3}, we have that for each potential $\xi(P_{u})$, $\varphi^{l}(\xi(P_{u}))-\varphi^{m}(\xi(P_{u})) \in \mathcal{F}_{E}(M)^{\geq m+1}$ for all $l>m$. Now choose $l>m$ such that $\chi(P_{u})-\varphi^{l}(\xi(P_{u})) \in \mathcal{F}_{E}(M)^{\geq m+1}$. Thus $\chi(\xi(P_{u}))-\varphi^{m}(\xi(P_{u})) \in \mathcal{F}_{E}(M)^{\geq m+1}$. Then
\begin{center}
$\underline{c}(\chi(P_{u}))(x)=(\underline{\vartheta}\rho)^{m}\underline{\xi}(u)$
\end{center}
Here $\underline{\xi}: \mathcal{Z}(T) \rightarrow F^{B}$ is a regular function and $(\underline{\vartheta}\rho)^{m}$ is a polynomial function. Then the function sending each $u \in \mathcal{Z}(T)$ into $\underline{c}(\chi(P_{u}))(x)$ is a regular function. The proof is now complete. 
\end{proof}
Take $k \in [1,n]$ and let $\mu_{k}M$ be the premutation of a $Z$-free $E$-$E$-bimodule $M$ with $Z$-free basis $A$. Take $\mu_{k}A$ the $Z$-free basis of $\mu_{k}M$ consisting of the elements $A \cap \overline{e}_{k}M\overline{e}_{k}$; the elements $B_{2} \cap Me_{k}M$; the elements $a^{\ast} \in (e_{k}M)^{\ast}$ with $a \in A \cap e_{k}M$, and the elements $^{\ast}a \in ^{\ast}(Me_{k})$ with $a \in Me_{k}$. Take $\tilde{B}_{l}^{(k)}$ the corresponding $F$-basis for $((\mu_{k}M)^{\otimes l})_{cyc}$ and $\tilde{B}^{(k)}=\displaystyle \bigcup_{l=2}^{\infty} \tilde{B}_{l}^{(k)}$. We recall (see Proposition \ref{prop7.3}) that there exists an isomorphism
\begin{center}
$[-]: \bar{e}_{k}\mathcal{F}_{E}(M)\bar{e}_{k} \rightarrow \mathcal{F}_{E}(\bar{e}_{k}M\bar{e}_{k} \oplus Me_{k}M)$
\end{center}
that sends elements of $B$ into elements of $\tilde{B}^{(k)}$. Therefore, there exists a polynomial function
\begin{center}
$\underline{\lambda}: F^{\bar{e}_{k}B\bar{e}_{k}} \rightarrow F^{\tilde{B}^{(k)}}$
\end{center}
such that for any $P \in \bar{e}_{k}\mathcal{F}_{E}(M)\bar{e}_{k}$, we have $\underline{c}([P])=\underline{\lambda}(\underline{c}(P))$. Moreover, there exists a polynomial function $\underline{\kappa}: F^{B} \rightarrow F^{\bar{e}_{k}B\bar{e}_{k}}$ such that $\underline{c}(\kappa(P))=\underline{\kappa}(\underline{c}(P))$. We have the polynomial function $\underline{\mu}_{k}=\underline{\kappa}\underline{\lambda}+\underline{c}(\Delta_{k})$ such that
\begin{center}
$\underline{c}(\mu_{k}P)=\underline{\mu}_{k}(\underline{c}(P))$
\end{center}
for any potential $P \in \mathcal{F}_{E}(M)$. With this notation we have the following
\end{section}
\begin{prop} \label{prop11.11} Suppose $P_{0}$ is a potential in $\mathcal{F}_{E}(M)$ such that the potential $\mu_{k}P_{0}$ is fully reducible. Let $\mathcal{A}$ be a dominant polarization of $\mu_{k}M$, and let $B^{(k)}=\mathcal{F}_{E}(M_{\mathcal{A}_{2}}) \cap \tilde{B}^{(k)}$. Then there exists a polynomial $T \in F[x]_{x \in B}$ and a regular function $\underline{\tilde{\mu}}_{k}: \mathcal{Z}(T) \rightarrow F^{B^{(k)}}$ such that for any potential $P$ in $\mathcal{F}_{E}(M)$ with $\underline{c}(P) \in \mathcal{Z}(T)$, there is a potential $\tilde{\mu}_{k}(P) \in \mathcal{F}_{E}(M_{\mathcal{A}_{2}})$ such that $P$ is right-equivalent to
\begin{center}
$\displaystyle \sum_{\alpha \in M_{\mathcal{A}^{\sigma}}} \alpha i(\alpha)+\tilde{\mu}_{k}(P)$
\end{center}
and $\underline{c}(\tilde{\mu}_{k}(P))=\underline{\tilde{\mu}}_{k}(\underline{c}(P)).$
\end{prop}
\begin{proof} Since $\mu_{k}P_{0}$ is fully reducible, then by Proposition \ref{prop11.10}, there exists a polynomial $T_{1} \in F[x]_{\tilde{B}}$ and a regular function $\underline{\chi}: F^{\tilde{B}} \rightarrow F^{\tilde{B}}$ such that $\underline{c}(\mu_{k}P_{0}) \in \mathcal{Z}(T_{1})$ and for all potentials $Q$ with $\underline{c}(Q) \in \mathcal{Z}(T_{1})$, $Q$ is right-equivalent to a potential $\chi(Q)=\displaystyle \sum_{\alpha \in M_{\mathcal{A}^{\sigma}}} \alpha i(\alpha)+\chi(Q)^{\geq 3}$, with $\chi(Q)^{\geq 3} \in \mathcal{F}_{E}(M_{\mathcal{A}_{2}})$. Here $\underline{\mu}_{k}(P_{0}) \in \mathcal{Z}(T_{1})$, then by Remark \ref{rem11.7}, there exists a non-zero polynomial $T \in F[x]_{x \in B}$ such that $\underline{\mu}_{k}(\mathcal{Z}(T)) \subseteq \mathcal{Z}(T_{1})$ and $\underline{c}(P_{0}) \in \mathcal{Z}(T)$. Now take the regular function
\begin{center}
$\underline{\tilde{\mu}}_{k}=\underline{r} \underline{\chi} \underline{\mu}_{k}: \mathcal{Z}(T) \rightarrow F^{B^{(k)}}$
\end{center}
where $\underline{r}: F^{\tilde{B}} \rightarrow F^{B^{(k)}}$ is given by the restriction map. Then setting $\tilde{\mu}_{k}(P)=\chi(\mu_{k}P)^{\geq 3}$ yields the desired result.
\end{proof}

\begin{prop} \label{prop11.12} Let $k_{1},\ldots,k_{l}$ be a sequence of elements of $\{1,\ldots,n\}$. Suppose that there exists a potential $P_{0} \in \mathcal{F}_{E}(M)$ such that $\overline{\mu}_{k_{l}} \overline{\mu}_{k_{l-1}} \cdots \overline{\mu}_{k_{1}}P_{0}$ is defined. Then there is a non-zero polynomial $T \in F[x]_{x \in B}$ such that $\overline{\mu}_{k_{l}} \overline{\mu}_{k_{l-1}} \cdots \overline{\mu}_{k_{1}}P$ is defined for every potential $P$ with $T(\underline{c}(P)) \neq 0$. 
\end{prop}

\begin{proof} We prove our proposition by induction on $l$. The case $l=1$ follows by Proposition $\ref{prop11.11}$. Hence we may assume that there are non-zero polynomials $T_{1} \in F[x]_{x \in B}$, $T_{2} \in F[x]_{x \in B^{(k_{1})}}$, and regular functions $\underline{\phi}_{1}: F^{B} \rightarrow F^{B^{(k_{1})}}$, $\underline{\phi}_{2}: F^{B^{(k_{1})}} \rightarrow F^{B'}$, where $B'=(\ldots(B^{(k_{1})})^{(k_{2})} \ldots)^{(k_{l})}$ such that $\underline{c}(P_{0}) \in \mathcal{Z}(T_{1})$, $\underline{c}(\overline{\mu}_{k_{1}}P_{0}) \in \mathcal{Z}(T_{2})$. Moreover, for any potential $Q$ in $\mathcal{F}_{E}(M)$, with $\underline{c}(Q) \in \mathcal{Z}(T_{1})$, $\overline{\mu}_{k}(Q)$ is defined and $\underline{c}(\overline{\mu}_{k}(Q))=\underline{\phi}_{1}(\underline{c}(Q))$; also, for any potential $Q' \in \mathcal{F}_{E}(\mu_{k_{1}}M)$, with $\underline{c}(Q') \in \mathcal{Z}(T_{2})$, $\overline{\mu}_{k_{l}} \cdots \overline{\mu}_{k_{2}}Q$ is defined and $\underline{c}(\overline{\mu}_{k_{l}} \cdots \overline{\mu}_{k_{2}}Q)=\underline{\phi}_{2}(\underline{c}(Q))$. 
We have $\underline{c}(\overline{\mu}_{k}P_{0})=\underline{\phi}_{1}(\underline{c}(P_{0})) \in \mathcal{Z}(T_{2})$. By Remark \ref{rem11.7}, there exists a polynomial $T \in F[x]_{x \in B}$ such that $\mathcal{Z}(T) \subseteq \mathcal{Z}(T_{1})$ and $\underline{\phi}_{1}(\mathcal{Z}(T)) \subseteq \mathcal{Z}(T_{2})$. Consider the regular function $\underline{\phi}=\underline{\phi}_{2}(\underline{\phi}_{1})|_{\mathcal{Z}(T)}: \mathcal{Z}(T) \rightarrow F^{B'}$. Take $P$ a potential in $\mathcal{F}_{E}(M)$ with $\underline{c}(P) \in \mathcal{Z}(T) \subseteq \mathcal{Z}(T_{1})$, then $\overline{\mu}_{k}(P)$ is defined and $\underline{c}(\overline{\mu}_{k}P)=\underline{\phi}_{1}(\underline{c}(P)) \in \mathcal{Z}(T_{2})$, so $\overline{\mu}_{k_{l}} \cdots \overline{\mu}_{k_{1}}(P)$ is defined and $\underline{\phi}(\underline{c}(P))=\underline{c}(\overline{\mu}_{k_{l}} \cdots \overline{\mu}_{k_{1}}(P))$. The result follows.  
\end{proof}

\begin{lemma} \label{lem11.13} For any $k \in [1,n]$, there exists a potential $P \in \mathcal{F}_{E}(M)$ such that the mutation $\overline{\mu}_{k}(P)$ is defined. 
\end{lemma}

\begin{proof} If $\mu_{k}M$ is $2$-acyclic, then for any potential $P \in \mathcal{F}_{E}(M)$, one has $\mu_{k}P=\overline{\mu}_{k}P$; so we may assume that $\mu_{k}M$ has $2$-cycles. Let $\mathcal{A}$ be a dominant polarization of $\mu_{k}M$. Then we have an injective function $i: M_{\mathcal{A}^{\sigma}}Ê\rightarrow M_{\mathcal{A}}$. We have the quadratic potential $Q=\displaystyle \sum_{\alpha \in M_{\mathcal{A}^{\sigma}}} \alpha i(\alpha)$, with $\Theta(Q)=\mu_{k}M_{\mathcal{A}^{\sigma}} \oplus \mu_{k}M_{\mathcal{A}_{1}}$. One can see that
\begin{center}
$\mathcal{F}_{E}(\bar{e}_{k}M\bar{e}_{k} \oplus Me_{k}M) \cong \bar{e}_{k}\mathcal{F}_{E}(M)\bar{e}_{k}$.
\end{center}
Therefore there exists a potential $P \in \bar{e}_{k}\mathcal{F}_{E}(M)\bar{e}_{k}$ with $[P]=Q$. Consequently, $\mu_{k}(P)=Q+\Delta_{k}$ and hence $\overline{\mu}_{k}(P)=\Delta_{k}$. This completes the proof.
\end{proof}

\begin{theorem} \label{theo11.14} Let $k_{1},\ldots,k_{l}$ be any sequence of elements of $\{1,\ldots,n\}$. Then there exists a potential $P \in \mathcal{F}_{E}(M)$ such that $\overline{\mu}_{k_{l}} \cdots \overline{\mu}_{k_{1}}(P)$ is defined.
\end{theorem}

\begin{proof} We argue by induction on $l$. The case $l=1$ follows by Lemma \ref{lem11.13}. Suppose that there exists a potential $Q_{0} \in \mathcal{F}_{E}(\mu_{k}M)$ such that $\overline{\mu}_{k_{l}} \cdots \overline{\mu}_{k_{2}}(Q_{0})$ is defined. We also have a potential $Q_{1} \in \mathcal{F}_{E}(\mu_{k_{1}}M)$ such that $\overline{\mu}_{k_{1}}(Q_{1})$ is defined. Then there are non-zero polynomials $T_{1},T_{2} \in F[x]_{x \in B^{(k_{1})}}$ such that for $Q$ with $\underline{c}(Q) \in \mathcal{Z}(T_{1})$, $\overline{\mu}_{k_{l}} \cdots \overline{\mu}_{k_{2}}(Q)$ is defined; and for $Q'$ with $\underline{c}(Q') \in \mathcal{Z}(T_{2})$, $\overline{\mu}_{k}(Q')$ is also defined. Here $T_{1}T_{2}$ is a non-zero polynomial in $F[x]_{x \in B^{(k_{1})}}$. Since $F$ is an infinite field, then let $u \in \mathcal{Z}(T_{1}T_{2})$ and take $Q_{u}=\displaystyle \sum_{x \in B^{(k_{1})}} u(x)x$. Then $\underline{c}(Q_{u})=u \in \mathcal{Z}(T_{1}T_{2}) \subseteq \mathcal{Z}(T_{1}) \cap \mathcal{Z}(T_{2})$. Therefore $\overline{\mu}_{k_{l}} \cdots \overline{\mu}_{k_{2}}(Q_{u})$ is defined, and likewise, $\overline{\mu}_{k_{1}}(Q_{u})$ is defined. We have $\overline{\mu}_{k_{1}}(Q_{u}) \in \mathcal{F}_{E}(\overline{\mu}_{k_{1}} \overline{\mu}_{k_{1}}M) \cong \mathcal{F}_{E}(M)$, then there exists a potential $P \in \mathcal{F}_{E}(M)$ such that $\overline{\mu}_{k_{1}}Q_{u}$ is right-equivalent to $P$. Therefore, $\overline{\mu}_{k_{1}}(P)$ is defined and it is right-equivalent to $Q_{u}$. It follows that $\overline{\mu}_{k_{l}} \cdots \overline{\mu}_{k_{1}}(P)$ is defined. This completes the proof.
\end{proof}
\end{section}



\begin{thebibliography}{12}
\bibitem{1} R. Bautista, D. L\'{o}pez-Aguayo. \textit{Tensor algebras and decorated representations}. \href{https://arxiv.org/abs/1606.01974}{arXiv:1606.01974}. 

\bibitem{2} R.Bautista, L. Salmer\'{o}n, R. Zuazua, \textit{Differential tensor algebras and their module categories}. London Mathematical Society Lecture Note Series, 362. Cambridge University Press, Cambridge, (2009).

\bibitem{3} L. Demonet. \textit{Mutations of group species with potentials and their representations. Applications to cluster algebras}. \href{https://arxiv.org/abs/1003.5078}{arXiv:1003.5078}.

\bibitem{4} H. Derksen, J. Weyman and A. Zelevinsky. \textit{Quivers with potentials and their representations I: Mutations}. Selecta Math. 
\textbf{14} (2008), no. 1, 59-119. \href{https://arxiv.org/abs/0704.0649}{arXiv:0704.0649}.

\bibitem{5} S. Fomin and A. Zelevinsky. \textit{Cluster algebras I. Foundations}. J. Amer. Math. Soc. \textbf{15} (2002), no.2, 497-529. \href{https://arxiv.org/abs/math/0104151}{arXiv:0104151}.

\bibitem{6} J. Geuenich, D. Labardini-Fragoso, \textit{Species with potential arising from surfaces with orbifold points of order 2, Part I: one choice of weights}. Mathematische Zeitschrift, \textbf{286} (2017), no. 3-4, 1065-1143. \href{http://arxiv.org/abs/1507.04304}{arXiv:1507.04304}.

\bibitem{7} J. Geuenich, D. Labardini-Fragoso, \textit{Species with potential arising from surfaces with orbifold points of order 2, Part II: arbitrary weights}. International Mathematics Research Notices. Advanced online publication. DOI: 10.1093/imrn/rny090 \href{https://arxiv.org/abs/1611.08301}{arXiv:1611.08301}.

\bibitem{8} D. Labardini-Fragoso and A. Zelevinsky. \textit{Strongly primitive species with potentials I: Mutations}. Bolet\'{i}n de la Sociedad Matem\'{a}tica Mexicana (Third series), Vol. \textbf{22}, (2016), Issue 1, 47-115. \href{http://arxiv.org/abs/1306.3495}{arXiv:1306.3495}.

\bibitem{9} D. L\'{o}pez-Aguayo (2018). \textit{A note on species realizations and nondegeneracy of potentials}. Journal of Algebra and Its Applications. Advanced online publication. DOI: 10.1142/S0219498819500245.

\bibitem{10} B.Nguefack. \textit{Potentials and Jacobian algebras for tensor algebras of bimodules}. \href{http://arxiv.org/abs/1004.2213}{arXiv:1004.2213}.

\bibitem{11} S.Roman, \textit{Field Theory}. Graduate Texts in Mathematics, 158. Springer-Verlag New York, 2006.

\bibitem{12} G.-C.Rota, B.Sagan and P.R.Stein. \textit{A cyclic derivative in noncommutative algebra}. Journal of Algebra. \textbf{64} (1980) 54-75.

\end{thebibliography}
\end{document}